\documentclass[12pt,a4paper,english]{article}
\usepackage{amsfonts,amssymb,latexsym,xspace,textcomp}
\usepackage{amsmath,enumerate}
\usepackage{amsthm}
\usepackage{graphicx}
\usepackage{color}
\usepackage[english]{babel}
\usepackage[latin1]{inputenc}
\usepackage[linktocpage]{hyperref}
\numberwithin{equation}{subsection}
\usepackage{titlesec}
\titleformat{\subsubsection}[runin]
{\normalfont\large\bfseries}{\thesubsubsection}{1em}{}
\theoremstyle{definition}
\newtheorem{theorem}[subsubsection]{Theorem}
\newtheorem{definition}[subsubsection]{Definition}

\newtheorem{remark}[subsubsection]{Remark}
\newtheorem{proposition}[subsubsection]{Proposition}
\newtheorem{lemma}[subsubsection]{Lemma}
\newtheorem{corollary}[subsubsection]{Corollary}
\begin{document}

\begin{center}
\textbf{{\Large Singularities of Fredholm maps with one-dimensional kernels, II: \\ \vspace{10pt} 
Local behaviour and pointwise conditions}}\footnote{This research was partially supported by MIUR project ``Elliptic and Hamiltonian Differential Problems and their applications''.}  \\
\vspace{50pt}
\textbf{F. Balboni} \\ \vspace{8pt} 
\textit{\footnotesize Dipartimento di Matematica, Università degli Studi ``Tor Vergata'', Via della Ricerca Scientifica, 00133, Roma, Italy}\\ 
{\footnotesize E-mail: balboni@mat.uniroma2.it} \\ \vspace{18pt} 
and \\  \vspace{18pt} 
\textbf{F. Donati   } \textbf{\textdied}\\
\vspace{8pt} 
\textit{\footnotesize Dipartimento di Matematica, Università degli Studi ``Tor Vergata'', Via della Ricerca Scientifica, 00133, Roma, Italy}\\ 
\vspace{18pt} 
\end{center}
\textbf{\footnotesize Abstract}\bigskip  \\
\indent \footnotesize In analogy to what happens in finite dimensions we state the Normal Form Theorem for $k$-singularities, introduced in the previous paper of the series. By means of that we study the local behaviour near a singularity i.e. we deduce local results of existence and multiplicity of solutions for the equation $F(x)=y$ where $F$ is a 0-Fredholm map and $x$ belongs to a suitable neighbourhood of a singular point $x_{o}$, once $x_{o}$ is identified as one of the three kinds of singularities defined in the first paper. To this end we also start to seek alternative strategies for the determination of the type of a given singularity according to our classification. Here we give a pointwise approach for lower-order singularities which is coherent with the Thom-Boardman classification in finite dimensions. We conclude by applying the pointwise condition to a differential problem where, under suitable hypotheses, we determine a swallow's tail singularity (3-singularity). 

\normalsize
\section*{{\large INTRODUCTION}}
\quad This paper is the natural continuation of \cite{B-D 1}, where a complete classification of simple singularities of Fredholm maps with index 0 was given. In this view we refer to \cite{B-D 1} for a thorough comparison with the existing literature and a more extensive bibliography. Here we are interested in analysing the local behaviour of a smooth 0-Fredholm map near a simple singularity, i.e. at a point where the Fréchet derivative has a one-dimensional kernel. Our aim is to obtain information, as precise as possible, on the number of solutions to the given equation near such a singular point. More precisely we study the equation $F(u) = h$ where $F$ is the given map, $u$ and $h$ belong to suitable neighbourhoods of $u_{o}$ and $h_{o}$, for $u_{o}$ a fixed singularity and $h_{o}=F(u_{o})$. This study is quite successful when dealing with the $k$-singularities described in \cite{B-D 1}, i.e. the infinite-dimensional analogues of Morin singularities. Indeed, it is possible to state, in the Banach space setting, a Normal Form Theorem for $k$-singularities which generalizes the well-known Normal Form Theorem for Morin singularities in the finite-dimensional case (cf. Theorem \ref{Teo111} and following remark). This means that the map to be studied is equivalent, up to local diffeomorphisms, to a suitable polynomial map. As a consequence, quite sharp results of existence and multiplicity of solutions to the equation $F(u) = h$ can be deduced on a neighbourhood of a $k$-singularity. Moreover, as we will show, similar information can also be obtained for the maximal $k$-transverse singularities introduced in \cite{B-D 1}. This kind of results, in conjunction with other tools developed in \cite{B-D 3}, are useful to obtain multiple solutions to some nonlinear differential problems formulated as operator equations in suitable Banach spaces, \cite{B-D 5}. An anticipation of this approach to the study of differential problems is given in this paper, see Section \ref{ss25}. It may be worthwhile to recall that the use of Singularity Theory in studying nonlinear differential problems is a process made up of three steps, i.e. classification (cf. \cite{B-D 1}), knowledge of the local behaviour (here) and identification of the singular point to be considered (here and in \cite{B-D 3}). Hence we start here to present a possible approach to the third step. In order to identify the considered singularity it would be desirable to have quite easy conditions to verify. However if we use the very definition of  $k$-singularity or maximal $k$-transverse singularity then we can be faced with quite cumbersome computations, and this is all the more true as the value of $k$ is larger. Nevertheless for ``low-order'' singularities, i.e. for $k=1,2,3,4$, such a direct approach is reasonable and, in fact, we will show that it is also possible to pass from \textit{local conditions}, i.e. on a neighbourhood, to \textit{pointwise conditions}, i.e. at the singular point, with the advantage of easier computations. Of course one needs a completely different approach when dealing with larger values of the integer $k$, as can happen with solution curves, and for these cases we refer to  \cite{B-D 3}. A comparison between the pointwise conditions presented here and the existing literature can be found, for the interested reader, in Remark \ref{Rem216}.\\
\indent This paper is organized as follows. In Section 1.1 we state the Normal Form Theorem for $k$-singularities, Theorem \ref{Teo111}, and then, after a technical result on the local behaviour of a polynomial map, we state and prove the existence and multiplicity of solutions to the equation  $F(u)=h$,  for  $u$ near a $k$-singularity  $u_{o}$, in Theorem \ref{Teo115} (Local Multiplicity Theorem). Roughly speaking, for this case we obtain the existence of  \textit{at  most} $k+1$ solutions near $u_{o}$ with some  $\tilde{h}$  near $F(u_{o})$  having exactly $k+1$ solutions. A slightly less precise result is then proved for $u$ near a $(k+1)$-transverse singularity. Namely, in this case we can only show the existence of \textit{at least} $k+1$ solutions. As a consequence we can still study what happens near a maximal $k$-transverse or an $\infty$-transverse singularity, even though for these two kinds of singularities there is no possibility to provide a Normal Form Theorem (cf. Section \ref{ss12}). We then consider, in Section \ref{ss13}, an important feature of $\infty$-transverse singularities of analytic maps: there is always a solution curve passing through an $\infty$-transverse singularity. As some counterexamples show, in general this property does not hold for non-analytic, $\infty$-transverse singularities and quite different behaviours can occur. In Chapter \ref{s2} we give, for $k=1,2,3,4$, the explicit form of the pointwise conditions. Particularly, in Section \ref{ss21} we state such conditions and we provide the terminology necessary to deal with and prove them. Their proofs rely on the study of suitable equations on the strata of singularities for all integers $k$ and on their relationship with the local conditions. This study is carried through in Sections \ref{ss22} and \ref{ss23}, while Section \ref{ss24} is devoted to a technical result, of combinatorial nature, we need to make sure that the pointwise conditions are indipendent from the vectors we choose to satisfy them. Finally, in Section \ref{ss25} we provide a simple application of pointwise conditions to the Neumann problem of a Liénard equation   
\begin{equation*}
\text{(P)}\; \begin{cases}
u''+f(u)u'+g(u)=h \qquad \text{in  }(0,\pi)\\
u'(0)= 0
=u'(\pi).
\end{cases}
\end{equation*}
Under suitable assumptions for $ f $ and $ g $ we show that $u\equiv 0$ is a swallow's tail (3-singularity) for the map naturally associated with the problem. Actually, the proof relies on pointwise conditions suitably simplified because of symmetry of the problem. This seems to be one of the first explicit examples of such singularities and we refer to Section 2.5 for a detailed discussion. \\
\indent For the reader's convenience, a table of contents is found below.  

\tableofcontents
\vspace{10 pt}
\section{Local Behaviour near a Singularity}\label{s1}
\subsection{The Normal Form Theorem for  $k$-singularities and its Consequences}\label{ss11}
\quad Let $F:U\subseteq X \rightarrow V \subseteq Y$ be a smooth, 0-Fredholm map between open subsets $U,V$ of the $B$-spaces $X,Y$. Let $u_{o}\in U$ be a simple singularity for $F$ which satisfies, for some integer $k$, one of the conditions given in \cite{B-D 1} (see Definition 2.1.1). For each of these conditions we shall provide quite precise information on the number of solutions to the equation $F(u)=h$, for $u$ near $u_{o}$ and $h$ near $h_{o}:=F(u_{o})$. First, when $u_{o}$ is a $k$-singularity for $F$, in order to study the local behaviour of $F$ near $u_{o}$ we can adopt the same strategy used for the Morin singularities in the finite-dimensional case. We thus need to establish a \textquotedblleft normal form\textquotedblright $ \, $  theorem, i.e. to show that $F$ is locally equivalent, near $u_{o}$, to a suitable simpler map, and this also allows knowing the number of solutions to $F(u)=h$. This is accomplished with Theorem \ref{Teo111}, the Normal Form Theorem for $k$-singularities, which states that a map, near a $k$-singularity, is locally equivalent to another one having a polynomial form that we call here generalized Whitney map (cf \cite{B-D 1}, Example 2.6.8). We incidentally note that the $k$-singularities coincide, in the finite-dimensional case, with Morin singularities as will be sketched in the next section. Hence, at least in order to get analogous results for the other kinds of singularities, considered in Definition 2.1.1 of \cite{B-D 1}, we will just combine the information related to the stratification of singularities near $u_{o}$ with the normal form theorem.\\
\indent However, we should note that the knowledge of the multiplicity of solutions is only a step towards significant applications to nonlinear differential problems. In fact, we also need an operative way to determine the kind of the considered singularity. This problem is studied in the case of low-order singularities in Chapter \ref{s2}, while for the general case we refer to \cite{B-D 3}.\\
\par
\begin{theorem}\label{Teo111}(Normal Form Theorem). \textit{Let} $U,V$ \textit{be open subsets of the} $B$-spaces $X,Y$\textit{ respectively, and }$F:U\subseteq X \rightarrow V \subseteq Y$\textit{ a }$C^{\infty}$ 0-\textit{Fredholm map. Let }$u_{o}\in U$\textit{ and }$k\geq 1$\textit{ an integer. Then the following conditions are equivalent:} \medskip \\ \medskip
NF1)  $u_{o}$\textit{ is a } $k$-singularity\textit{ for }$F$ ;\\
NF2)  \textit{there exists a  local commutative diagram of class} $C^{\infty}$ ($C^{\infty}$ l.c.d.)\textit{ given by} 
\qquad \begin{center}
\begin{tabular}{ccc}
&$F \quad $& \\
$u_o \in U \subseteq X \quad $&$ \rightarrow \quad $&$ V \subseteq Y $ \bigskip  \\
$\alpha \downarrow $&&$ \downarrow \beta$ \bigskip \\ 
$\underline{0}\in \mathbb{R}^{k} \times Z \quad $&$ \rightarrow  \quad $&$ \mathbb{R}^{k} \times Z$ ,\\
&$w_{k,Z} \quad $&\\
\end{tabular}
\end{center}
\textit{where }$Z$\textit{ is a suitable }$B$\textit{-space}, $\underline{0}$\textit{ is the origin of } $\mathbb{R}^{k} \times Z$\textit{ and }$w_{k,Z}$\textit{ is the generalized Whitney map } $w_{k,Z}(t,t_{1},\ldots,t_{k-1},z)=(t^{k+1}+t_{k-1}t^{k-1}+\ldots+t_{2}t^{2}+t_{1}t,t_{1},\ldots,t_{k-1},z).$
\end{theorem}
\begin{remark} \label{Rem112} The definition of a $C^{\infty}$ l.c.d. can be found in Subsection 1.2.1. of \cite{B-D 1}. The rather long proof that NF1) $\Rightarrow$ NF2) is given in \cite{Ba2}; here we shall only add that it follows a strategy similar to the one adopted in the finite-dimensional case (e.g. see \cite{G-G}, theorem 4.1, chapter VII). The main tool for this kind of proofs is the so-called Preparation Theorem which also exists for Banach spaces (cf.\cite{ Ba1}, \cite{B-C-T}) though with a weaker statement.\\
\indent On the other hand we can easily prove that NF2) $\Rightarrow$ NF1) by using the Invariance Theorem (cf. Theorem 2.5.5 in \cite{B-D 1}), which essentially states that the kind of singularity is not affected by changes of coordinates. Namely, in Example 2.6.8 of \cite{B-D 1} we showed that the origin $\underline{0}$ is a $k$-singularity for $w_{k,Z}$. Hence, from the Invariance Theorem we immediately obtain that $u_{o}$ is a $k$-singularity for $F$.
\end{remark}
From the above theorem it follows that, in order to study the equation  $F(u)=h$  near a $k$-singularity, the local behaviour of the map $w_{k,Z}$ must be known; to this end, it just suffices to study the Whitney map $w_{k}$. Setting $\textbf{t}:=(t,t_{1},\ldots,t_{k-1})$, we recall that the Whitney map is defined by $w_{k}(\textbf{t})=(t^{k+1}+t_{k-1}t^{k-1}+ \ldots + t_{2}t^{2}+t_{1}t,t_{1},\ldots,t_{k-1})$. Moreover, a point $\hat{\textbf{t}}$ is said to be a \textit{regular point} for $w_{k}$ if the linear map $w_{k}^{\,\prime}(\hat{\textbf{t}}):\mathbb{R}^{k} \rightarrow \mathbb{R}^{k}$  is an isomorphism. Consequently, for the map $w_{k}$  we have the following\\
\begin{proposition}\label{Pro113}\textit{ Let }$w_{k}:\mathbb{R}^{k} \rightarrow \mathbb{R}^{k}$  \textit{ as above}\textit{. Then the equation }$w_{k}(\textbf{t})=\textbf{s}$\textit{ has at most }$k+1$ solutions, $\forall \, \textbf{s} \in \mathbb{R}^{k}$.\textit{ Moreover, for any pair of neighbourhoods }$\hat{U}, \hat{V}$\textit{ of the origin, there exist points } $\hat{\textbf{s}}\in \hat{V}$\textit{ such that the equation }$w_{k}(\textbf{t})=\hat{\textbf{s}}$\textit{ has exactly } $k+1$\textit{ distinct solutions }$\hat{\textbf{t}}_{1},\ldots,\hat{\textbf{t}}_{k+1}$\textit{ in }$\hat{U}$\textit{ which are regular points for } $w_{k}$.
\end{proposition}
\indent In order to prove the above result we need the following
\begin{lemma}\label{Lem114}\textit{ Let } $k\geq 1$\textit{ be a fixed integer. Then, for each } $\varepsilon>0$\textit{, there exist }$\alpha_{0},\ldots,\alpha_{k-1}\in \mathbb{R}$\textit{ such that } $0<\lvert\alpha_{i}\lvert<\varepsilon,i=0,\ldots,k-1,$\textit{ and the polynomial } $x^{k+1}+\alpha_{k-1}x^{k-1}+\ldots+\alpha_{0}$\textit{ has exactly }$k+1$\textit{ distinct real roots }$x_{1},\ldots,x_{k+1}$\textit{ satisfying} 
\begin{equation*}
0<\lvert x_{i}\lvert<k \sqrt{\varepsilon},\;i=1,\ldots,k+1 .
\end{equation*}
\end{lemma}
\indent\textbf{ Proof. }We shall proceed by induction on the integer $k$. For $k=1$, let us take $\varepsilon>0$ and consider the polynomial $x^{2}+\alpha_{0}$. Then, for $-\varepsilon<\alpha_{0}<0$ this has two different real roots $x_{1}=\sqrt{-\alpha_{0}}$ and $x_{2}=\sqrt{-\alpha_{0}}$ satisfying  $0<\lvert x_{i}\lvert<\sqrt{-\alpha_{0}}<\sqrt{\varepsilon},i=1,2.$\\
Let us suppose the result to hold for the integer $k$ and prove it for $k+1$. Given $\varepsilon > 0$, by inductive hypothesis there exist $\alpha_{0},\ldots,\alpha_{k-1}\in \mathbb{R}$ such that $0<\lvert\alpha_{i}\lvert<\varepsilon,i=0,\ldots,k-1$, and the polynomial $x^{k+1}+\alpha_{k-1}x^{k-1}+\ldots+\alpha_{0}$ has  exactly $k+1$ distinct real roots $x_{1}, \ldots,x_{k+1}$ satisfying $0<\lvert x_{i}\lvert<k \sqrt{\varepsilon},i=1,\ldots,k+1$. Then it follows that the polynomial $Q(x):=x(x^{k+1}+\alpha_{k-1}x^{k-1}+\ldots+\alpha_{0})$ has $k+2$ distinct real roots $x_{1},\ldots,x_{k+1},x_{k+2}:=0$ satisfying $\lvert x_{i}\lvert<k\sqrt{\varepsilon},i=1,\ldots,k+2$. We now claim that:
\begin{align}\label{111}
\begin{split}
&\text{there exists }\delta_{\varepsilon}> 0\text{ such that, for }\lvert y \lvert<\delta_{\varepsilon},\text{ the equation }Q(x)=y\text{ has }k+2 \\
&\text{distinct real roots } z_{1}, \ldots. z_{k+2}\text{ satisfying }\lvert z_{i}-x_{i}\lvert < \sqrt{\varepsilon},i=1,\ldots,k+2.
\end{split}
\end{align}
Once (\ref{111}) is shown to be true we can easily conclude the proof of the Lemma. In fact, for $\beta_{0}$ such that $0<\lvert\beta_{0}\lvert < \text{min}(\varepsilon,\delta_{\varepsilon}),$ let us consider the polynomial
\begin{equation*}
Q(x)+\beta_{0}=x^{k+2}+\alpha_{k-1}x^{k}+\ldots+\alpha_{0}x + \beta_{0}=x^{k+2}+\beta_{k}x^{k}+\ldots+\beta_{1}x+\beta_{0} ,
\end{equation*}
with $\beta_{k}:=\alpha_{k-1},\ldots,\beta_{1}:=\alpha_{0}$. Then, by construction, one has $0<\lvert\beta_{i}\lvert < \varepsilon, i = 0,\ldots,k$. Furthermore, being $\lvert\beta_{0}\lvert < \delta_{\varepsilon}$, the above polynomial has $k+2$ distinct real roots $z_{1},\ldots, z_{k+1}, z_{k+2}$ which are not zero because $\beta_{0}\neq 0$, and which satisfy $0<\lvert z_{i}\lvert \leq \lvert z_{i}-x_{i}\lvert+\lvert x_{i}\lvert < \sqrt{\varepsilon} + k\sqrt{\varepsilon} =(k+1)\sqrt{\varepsilon},\,i=1,\ldots,k+2.$\\
\indent We still need to prove the validity of (1.1.1). By construction $Q(x)=\prod_{i=1}^{k+2}(x-x_{i})$, thus $Q'(x)=\sum_{i=1}^{k+2}\prod_{\!\!\!\!\scriptsize\begin{array}{c}
j=1 \\
j\neq 1
\end{array}}^{k+2}\!\!(x-x_{j})$ and therefore $Q'(x)=\prod_{\!\!\!\!\scriptsize\begin{array}{c}
j=1 \\
j\neq 1
\end{array}}^{k+2}\!\!(x_{i}-x_{j})\neq 0,i=1,\ldots,k+2,$ since the roots are distinct. By the Inverse Function Theorem (cf. e.g. \cite{A-M-R}, theorem 2.5.2),  $\forall \,i=1,\ldots,k+2$ there are intervals $A_{i},B_{i}$ in $\mathbb{R}$, with mutually disjoint $A_{i}$, such that $x_{i}\in A_{i},0\in B_{i}$ and $Q:A_{i}\rightarrow B_{i}$ is a diffeomorphism. Let us consider $(Q\lvert_{A_{i}})^{-1}:B_{i}\rightarrow A_{i}$. By continuity there exists $\delta_{\varepsilon}>0$ such that, for $\lvert y\lvert = \lvert y-0 \lvert < \delta_{\varepsilon}$, one has $\lvert(Q\lvert_{A_{i}})^{-1}(y)-(Q\lvert_{A_{i}})^{-1}(0)\lvert=\lvert(Q\lvert_{A_{i}})^{-1}(y)-x_{i}\lvert <\sqrt{\varepsilon},i=1,\ldots,k+2$. Then, by defining $z_{1}:=(Q\lvert_{A_{i}})^{-1}(y),i=1,\ldots, k+2$, one has that $Q(z_{i})=Q((Q\lvert_{A_{i}})^{-1}(y))=y$ and $\lvert z_{i}-x_{i}\lvert < \sqrt{\varepsilon}$, with the points $z_{i}$ being different from one another since $z_{i}\in A_{i}. \bracevert$\\
\par \textbf{Proof of Proposition \ref{Pro113}} The equation $w_{k}(\textbf{t})=\textbf{s}$, that is $w_{k}(t,t_{1},\ldots,t_{k-1})=(t^{k+1}+t_{k-1}t^{k-1}+ \ldots + t_{1}t,t_{1},\ldots,t_{k-1})=(s,s_{1},\ldots,s_{k-1})$, admits solutions \textit{iff}
\begin{align*}
t_{1}&=s_{1}\\
&\, :\\
t_{k-1}&=s_{k-1}\\
t^{k+1}+s_{k-1}t^{k-1}+\ldots +s_{1}t&=s \;.
\end{align*}
Since the t-polynomial $t^{k+1}+s_{k-1}t^{k-1}+\ldots +s_{1}t-s $ has at most $k+1$ real roots it follows that $w_{k}(\textbf{t})=\textbf{s}$ has at most $k+1$ solutions as well. Let $\hat{U}, \hat{V}$ be neighbourhoods of the origin and let us choose $\varepsilon>0$ such that the Cartesian products $(-M,M)^{k}$ and $(-\varepsilon,\varepsilon)^{k}$ are contained in $\hat{U}$ and $\hat{V}$, respectively, where $M:=(-\text{max}(k\sqrt{\varepsilon},\varepsilon)$. Consider, as in Lemma 1.1.4, $\alpha_{0},\ldots,\alpha_{k-1}\in \mathbb{R}$ such that $0<\lvert\alpha_{i}\lvert<\varepsilon,i=0,\ldots,k-1$ and that the polynomial $P(t):=t^{k+1}+\alpha_{k-1}t^{k-1}+\ldots+\alpha_{0}$ has $k+1$ distinct real roots $x_{1},\ldots,x_{k+1}$ with $0<\lvert x_{i}\lvert < k\sqrt{\varepsilon}, i=1,\ldots,k+1$. Set  $\hat{\textbf{s}}:=(-\alpha_{0},\alpha_{1},\ldots,\alpha_{k-1})\in \hat{V}$.  By construction the equation $w_{k}(\textbf{t})=\hat{\textbf{s}}$ has $k+1$ distinct roots $\hat{\textbf{t}}_{i}:=(x_{i},\alpha_{1},\ldots,\alpha_{k-1})\in \hat{U}, i=1,\ldots,k+1$. 
If we define $W_{k}(\textbf{t}):=t^{k+1}+t_{k-1}t^{k-1}+ \ldots + t_{1}t$ we have that $w_{k}(\textbf{t})=(W_{k}(\textbf{t}),t_{1},\ldots, t_{k-1})$, and then it follows that $w_{k}'(\textbf{t}):\mathbb{R}^{k} \rightarrow \mathbb{R}^{k}$ can be written in a matricial form as
\begin{equation*}
w_{k}'(\textbf{t})= \left[ \begin{array}{cc} 
\dfrac{\partial W_{k}}{\partial t}(\textbf{t}) \,\,* \,\,\cdots  & \qquad  \cdots\,\,*\\
0\\
\vdots\\
&\textbf{1}_{\mathbb{R}^{k-1}}\\
\vdots\\
0\end{array} \right],
\end{equation*}
where the symbols $*$ indicate entries that we do not need to specify. \medskip \\
Then one has that $w_{k}^{\prime}(\textbf{t})$ is invertible \textit{iff }$\dfrac{\partial W_{k}}{\partial t}(\textbf{t})\neq 0$. The last inequality just occurs for $\textbf{t}=\hat{\textbf{t}}_{i}=(x_{i},\alpha_{1},\ldots,\alpha_{k-1}),\, i=1,\ldots,k+1$. Indeed, being $P(t)=W_{k}(t, \alpha_{1},\ldots,\alpha_{k-1})+ \alpha_{0}$ one has that $\dfrac{\partial W_{k}}{\partial t}(\hat{\textbf{t}}_{i})=P'(x_{i})$. Since $P(t)$ has distinct roots $x_{i},i=1,\ldots,k+1$ then necessarily such roots are simple, that is $P'(x_{i})\neq 0$ (cf. proof of Lemma 1.1.4). $\bracevert$\\

Thanks to the previous proposition we are now able to prove the following result:
\begin{theorem}\label{Teo115}(Local Multiplicity Theorem).\textit{ Let }$U,V$\textit{ be open subsets of the B-spaces} $X,Y$\textit{ and }$F:U\subseteq X\rightarrow V\subseteq Y$\textit{ a }$C^{\infty}$ 0\textit{-Fredholm map. Let }$u_{o}\in U$\textit{ be a k}-singularity\textit{ for }$F$,\textit{ for some }$k\geq 1$.\textit{ Then there exist neighbourhoods }$U_{o}\subseteq U$\textit{ of }$u_{o}$ and $V_{o}\subseteq V$\textit{ of }$h_{o}:=F(u_{o})$\textit{ such that the equation }$F(u)=h$\textit{ has at most }$k+1$\textit{ solutions in }$U_{o}$\textit{ for }$h \in V_{o}$\textit{ and the set } 
\begin{center}
$\varOmega:=\{h\in V_{o}:F(u)=h \textit{ has }k+1 \textit{ distinct solutions in }U_{o}\}$
\end{center} \textit{is open in }$V_{o}$.\textit{ Moreover} $\varOmega$\textit{ is non-empty: in fact }$h_{o} \in \partial \varOmega$\textit{, where }$\partial$\textit{ means topological boundary. Precisely, for any pair of neighbourhoods }$\widetilde{U}_{o}\subseteq U$\textit{ of }$u_{o}$\textit{ and }$\widetilde{V}_{o}\subseteq V$\textit{ of }$h_{o}$\textit{ there exists }$\tilde{h}\in \widetilde{V}_{o}$\textit{ such that }$F(u)=\tilde{h}$\textit{ has exactly }$k+1$\textit{ distinct solutions }$\tilde{u}_{1},\ldots,\tilde{u}_{k+1}$\textit{ in }$\widetilde{U}_{o}$\textit{ which are regular points for }$F$.
\end{theorem}
\par \textbf{ Proof.} By using the Normal Form Theorem we can choose neighbourhoods $U_{o} \subseteq U$ of $u_{o}, V_{o}\subseteq V$ of $F(u_{o}),\hat{U}_{o}, \hat{V}_{o}\subseteq \mathbb{R}^{k}\times Z$ of the origin $\underline{0}$, along with a suitable $B$-space $Z$ and smooth surjective diffeomorphisms $\alpha: U_{o}\rightarrow \hat{U}_{o}, \beta : V_{o}\rightarrow \hat{V}_{o}$ such that $F(U_{o})\subseteq V_{o}, w_{k,Z}(\hat{U}_{o})\subseteq \hat{V}_{o},\alpha(u_{o})=\underline{0}$ and the following diagram is commutative:
\begin{align}\label{112}
\begin{tabular}{ccc}
&$F \quad $&\\
$u_{o} \in U_{o}\quad $&$ \rightarrow \quad $&$ V_{o}$ \bigskip  \\
$\alpha \downarrow $&&$ \downarrow \beta$ \bigskip \\ 
$\underline{0}\in \hat{U}_{o}\quad $&$ \rightarrow  \quad $&$ \hat{V}_{o}$. \\
&$w_{k,Z} \quad $&\\
\end{tabular}
\end{align}
Since $w_{k,Z}=w_{k}\times \textbf{1}_{\textbf{Z}}$, from Proposition 1.1.3 we obtain that the equation $w_{k,Z}(\textbf{t},z)=(w_{k}(\textbf{t}),z)=(\textbf{s},\zeta)$ has at most $k+1$ solutions in $\hat{U}_{o}$, for every $(\textbf{s},\zeta)\in \hat{V}_{o}$. Thus the equation $F(u)=h$ has at most $k+1$ solutions in $U_{o}$ for $h\in V_{o}$.\\
\indent Now suppose that we are given neighbourhoods $\widetilde{U}_{o}\subseteq U$ of $u_{o}$ and $\widetilde{V}_{o}\subseteq V$ of $h_{o}$. In order to prove the last part of the theorem it is not restrictive to assume that the neighbourhoods $\widetilde{U}_{o}$ and $\widetilde{V}_{o}$ are such that
\begin{equation*}
\widetilde{U}_{o}\subseteq U_{o},\widetilde{V}_{o}\subseteq V_{o}\, , \, \alpha(\widetilde{U}_{o})=\hat{U}\times W,\beta(\widetilde{V}_{o})=\hat{V} \times W,
\end{equation*}
where $\hat{U},\hat{V}$ are neighbourhoods of the origin of $\mathbb{R}^{k}$ and $W$ is a neighbourhood of the origin of $Z$.\\
\indent We contend that there exists a point $(\hat{\textbf{s}},\zeta)$ in $\hat{V}\times W$ such that $w_{k,Z}(\textbf{t},z)=(w_{k}(\textbf{t}),z)= (\hat{\textbf{s}},\zeta)$ has exactly $k+1$ distinct solutions $(\hat{\textbf{t}}_{i},\zeta)$ in $\hat{U}\times W,i=1, \ldots,k+1$ which are regular points for $w_{k,Z}$, i.e. $w'_{k,Z}(\hat{\textbf{t}}_{i},\zeta)= w_{k}'(\hat{\textbf{t}}_{i})\times \textbf{1}_{\textbf{Z}}$ is an isomorphism. Indeed, by using again Proposition 1.1.3, we can take $\hat{\textbf{s}}$ in $\hat{V}$ and $\hat{\textbf{t}}_{i}$ in $\hat{U}$ such that $w_{k}(\hat{\textbf{t}}_{i})=\hat{\textbf{s}}$ and $w'_{k}(\hat{\textbf{t}}_{i})$ is an isomorphism, $i=1,\ldots,k+1$. Then, for $\zeta$ arbitrarily chosen in $W$, let us take $z=\zeta$. It follows that $\tilde{h}:=\beta^{-1}(\hat{\textbf{s}},\zeta)$ is in $\widetilde{V}_{o}$ and the equation $F(u)=\tilde{h}$ has in $\widetilde{U}_{o}$ exactly  $k+1$ distinct solutions $\tilde{u}_{i}:= \alpha^{-1}(\hat{\textbf{t}}_{i},\zeta),i=1,\ldots, k+1$. Finally, the points $\tilde{u}_{i} $ are regular. In fact,  $F=\beta^{-1}w_{k,Z}\alpha$ and so  
\begin{center}
$F'(\tilde{u}_{i})=(\beta^{-1})'(w_{k,Z}(\hat{\textbf{t}}_{i},\zeta))w'_{k,Z}(\hat{\textbf{t}_{i}},\zeta)\alpha'(\tilde{u}_{i})$
\end{center}
is an isomorphism because $\alpha,\beta$ are diffeomorphisms and $w'_{k,Z}(\hat{\textbf{t}}_{i},\zeta)$ is an isomorphism. \\
\indent To show that $\varOmega$ is an open subset of $V_{o}$ we consider a point $\tilde{h}$ in $V_{o}$ with $k+1$ distinct solutions $\tilde{u}_{1},\ldots,\tilde{u}_{k+1}$ in $U_{o}$. Then $w_{k,Z}(\textbf{t},z)=\beta(\tilde{h})$ has $k+1$ distinct solutions in $\hat{U}_{o}$, namely $\alpha(\tilde{u}_{1}),\ldots,\alpha(\tilde{u}_{k+1})$. As seen in Proposition \ref{Pro113}, one can easily conclude that $\alpha(\tilde{u}_{1}),\ldots,\alpha(\tilde{u}_{k+1})$ are regular points for $w_{k,Z}$. Hence $\tilde{u}_{1},\ldots,\tilde{u}_{k+1}$ are regular points for $F$. By the Inverse Function Theorem, there exist neighbourhoods $U_{i}\subseteq U$ of $\tilde{u}_{i}$ and $V_{i}\subseteq V$ of $\tilde{h},i=1,\ldots,k+1$, such that $F:U_{i}\rightarrow V_{i}$ is a diffeomorphism. Up to shrinking, we can choose the subsets $U_{i}$ and $V_{i}$ such that $U_{i}\subseteq U_{o},V_{i}\subseteq V_{o}$ and $U_{i}\cap U_{j}= \emptyset,i\neq j$. Thus, if $h\in \cap^{k+1}_{i=1}V_{i}$, the equation $F(u)=h$ has a unique solution in $U_{i}$. Therefore $F(u) =h$ has $k+1$ distinct solutions in $U_{o}$.\\
\indent To conclude, since for any neighbourhoods $\widetilde{U}_{o}\subseteq U$ of $u_{o}$ and $\widetilde{V}_{o} \subseteq V$ of  $h_{o}$ there exists $\tilde{h}\in \widetilde{V}_{o} $ such that $F(u)=\tilde{h}$ has exactly $k+1$ distinct (regular) solutions in $\widetilde{U}_{o}$, it follows that $\varOmega$ is non-empty and more specifically that $h_{o}\in \partial \varOmega$: it is also clear that $h_{o}\notin \varOmega$ because, as a consequence of l.c.d. (\ref{112}), the equation $F(u)=h_{o}$ has $u_{o}$ as its unique solution. $\bracevert$  \\
\par  We recall that in \cite{B-D 1} we considered the first order periodic problem 
\begin{equation*}
\begin{cases}
u'+a(t)u^{2} + p(t)u^{4} = h  \qquad  \text{ in\:\:} (0,1)\\
u(0) = u(1) ,
\end{cases}
\end{equation*}
with $a(t)\in C^{0}([0,1])\setminus \{0\}$ and mean value zero, as well as $p(t) \in C^{0}([0,1])\setminus \{0\}$ with $\int_{0}^{1}\!p(t)\neq 0$. For this problem we proved that $u_{o}\equiv 0$ is a 3-singularity for the associated map $F\!:\!X:=C_{\#}^{1}([0,1])\rightarrow Y:= C^{0}([0,1])$, defined as $F(u)=u'+g(t,u)$, where $g(t,u)=a(t)u^{2} + p(t)u^{4}$. Thanks to the above Local Multiplicity Theorem we can say that there exist neighbourhoods $U_{o}\subseteq X$ of $u_{o}=0$ and $V_{o}\subseteq Y$ of $F(u_{o})=0$ such that the problem has at most four solutions in $U_{o}$ for $h\in V_{o}$. Moreover, there exist right-hand sides $h$ near $F(u_{o})=0$ with exactly four distinct solutions and, more precisely, the subset of such $h$ is a non-empty open subset. \\
Another example of swallow's tail singularity, this time in the class of differential problems of second order, is presented in Section \ref{ss25}. We explicitly note that in this case the method adopted to prove the result is completely different from the one used for the above first order problem.\\
\par 
Now we move onto $(k+1)$-transverse singularities; since this condition is weaker than the hypothesis of being $(k+1)$-singularity we obtain a somewhat less precise result which, however, could be equally useful in the applications.
\begin{theorem}\label{Teo116}
\textit{ Let }$U,V$\textit{ be open subsets of the B-spaces }$X, Y$\textit{ and }$F:U \subseteq X\rightarrow V\subseteq Y$\textit{ a }$C^{\infty}$ 0\textit{-Fredholm map. Let }$u_{o}\in U$\textit{ be a }$(k+1)$-transverse singularity\textit{ for }$F$\textit{, for some integer }$k\geq 1$.\textit{ Then, for any pair of neighbourhoods }$\widetilde{U}_{o} \subseteq U$\textit{ of }$u_{o}$\textit{ and }$\widetilde{V}_{o} \subseteq V$\textit{ of }$h_{o}:=F(u_{o})$,\textit{ there exists }$\tilde{h}\in \widetilde{V}_{o}$\textit{ such that the equation }$F(u)=\tilde{h}$\textit{ has at least }$k+1$\textit{ distinct solutions }$\tilde{u}_{1},\ldots,\tilde{u}_{k+1}$\textit{ in }$\widetilde{U}_{o}$ \textit{ which are regular points for }$F$.
\end{theorem}
\par \textbf{Proof.} We may take $\widetilde{U}_{o}, \widetilde{V}_{o}$ as open balls $B(u_{o};\varepsilon), B(h_{o};\varepsilon)$ for a given $\varepsilon> 0 $. By continuity there exists $\delta >0$ such that $\lVert u-u_{o}\lVert<\delta$ implies $\lVert F(u)-h_{o}\lVert<\varepsilon/2$. Since $u_{o}$ is $(k+1)$-transverse then, as we saw in Remark 2.4.6 in \cite{B-D 1}, near $u_{o}$ we have a one-codimensional stratification of manifolds
\begin{center}
$u_{o}\in S_{1_{k+1}}(F) \subseteq S_{1_{k}}(F) \subseteq \ldots \subseteq S_{1_{2}}(F)\subseteq S_{1_{1}}(F)=S_{1}(F) \subseteq X $.
\end{center}
Since $S_{1_{k+1}}(F)$ is, near $u_{o}$, a one-codimensional submanifold of $S_{1_{k}}(F)$ then $u_{o}$ is an accumulation point for the set $S_{1_{k}}(F)\setminus S_{1_{k+1}}(F)$. Hence, there exists $\tilde{u}\in B(u_{o};\text{min}(\delta,\varepsilon/2))$ such that $\tilde{u}\in S_{1_{k}}(F) \setminus S_{1_{k+1}}(F)$. From point b) of Theorem 2.5.4 in \cite{B-D 1} we have that for $\tilde{u}$ near $u_{o}$ the following equivalence holds:
\begin{center}
$\tilde{u}\in S_{1_{k}}(F)\setminus S_{1_{k+1}}(F) \Leftrightarrow
\tilde{u}\in S_{1}(F)$ and it is a $k$-singularity.\\
\end{center}
Thus $\tilde{u}$ is a $k$-singularity and  $\lVert F(\tilde{u})-h_{o}\lVert<\varepsilon/2$.\\
From Theorem \ref{Teo115} we deduce the existence of $\tilde{h}\in B(F(\tilde{u});\varepsilon/2)$ such that the equation $F(u)=\tilde{h}$ has exactly $k+1$ distinct solutions $\tilde{u}_{1},\ldots,\tilde{u}_{k+1}$ in $B(\tilde{u};\varepsilon/2)$, which are regular points for $F$. It follows that $\lVert \tilde{h}-h_{o}\lVert<\varepsilon$ and $\lVert \tilde{u_{1}}-u_{o}\lVert<\varepsilon$, for $i=1,\ldots,k+1. \bracevert$\\
\par  The next result is a straightforward consequence of the above theorem. It states that, near an  $\infty$-transverse singularity for $F$, the equation $F(u)=h$ has an arbitrarily large number of solutions. It is also useful to study the problems with solution curves seen in the Introduction, \cite{B-D 5}.\\
\begin{corollary}\label{Co117}
\textit{ Let }$U,V$\textit{ be open subsets of the B-spaces }$X,Y$\textit{ and }$F:U \subseteq X\rightarrow V \subseteq Y$\textit{ a }$C^{\infty}$ 0\textit{-Fredholm map. Let }$u_{o}\in U$\textit{ be an }$\infty$-transverse singularity\textit{ for }$F$.\textit{ Then, for any integer }$k$\textit{ and for any pair of neighbourhoods }$\widetilde{U}_{o}\subseteq U$\textit{ of }$u_{o}$\textit{ and }$\widetilde{V}_{o} \subseteq V$\textit{ of }$h_{o}:=F(u_{o})$\textit{ there exists }$\tilde{h}\in \widetilde{V}_{o}$\textit{ such that the equation }$F(u)=\tilde{h}$\textit{ has at least }$k+1$\textit{ distinct solutions }$\tilde{u}_{1}, \ldots,\tilde{u}_{k+1}$\textit{ in }$\widetilde{U}_{o}$\textit{ which are regular points for }$F$. 
\end{corollary}
\par \textbf{Proof. } By definition, $u_{o}$ is $(k+1)$-transverse for all $k\geq 1$.$\bracevert$ \\
\par In the previous results of existence and multiplicity we looked for \textit{regular} solutions $\tilde{u}_{i}$ of the equation $F(u)=h$. Near such solutions the map $F$ has nice properties which turn out to be quite useful in the applications. In fact, these properties guarantee that the problem $F(u)=h$ is well-posed in a neighbourhood of each $\tilde{u}_{i}$. This follows at once from the Inverse Function Theorem that enables us to say that the above map is a local diffeomorphism near a regular point. However, when we do not require the solutions of $F(u)=h$  to be regular points for the map $F$, we can generalize both Theorems \ref{Teo115} and \ref{Teo116} in the following way:
\begin{theorem}\label{Teo118}
\textit{ Let }$U,V$\textit{ be open subsets of the B-spaces }$X,Y$\textit{ and }$F:U \subseteq X \rightarrow V \subseteq Y$\textit{ a }$C^{\infty}$  0-\textit{Fredholm map. Let }$u_{o}\in U$\textit{ be a k}-transverse singularity\textit{ for }$F$,\textit{ for some integer }$k\geq 1$.\textit{ Then, for any pair of neighbourhoods }$\widetilde{U}_{o} \subseteq U$\textit{ of }$u_{o}$\textit{ and }$\widetilde{V}_{o}\subseteq V$\textit{ of }$h_{o}:=F(u_{o})$,\textit{ there exists }$\tilde{h} \in \widetilde{V}_{o}$\textit{ such that the equation }$F(u)=\tilde{h}$\textit{ has at least} $k+1$\textit{ distinct solutions in }$\widetilde{U}_{o}$.\textit{ More precisely, one of the following statements holds:}
\begin{itemize}
\item[i)]$F(u)=\tilde{h}$\textit{ has }$k+1$\textit{ distinct solutions }$\tilde{u}_{1},\ldots, \tilde{u}_{k+1}$ \textit{in }$\widetilde{U}_{o}$\textit{ which are regular points for }$F$; 
\item[ii)] $F(u)=\tilde{h}$\textit{ has a} regular curve \textit{of (singular) solutions in} $\widetilde{U}_{o}$, i.e.\textit{ there exist }$\varepsilon_{o}>0$ \textit{ and a }$C^{\infty}$ \textit{ map }$\gamma:(-\varepsilon_{o},\varepsilon_{o})\rightarrow \widetilde{U}_{o} \subseteq X$\textit{ such that }$\gamma^{\,'}(t)\neq 0$ \textit{ and }$F(\gamma(t))=\tilde{h}, \, \forall \, t\in (-\varepsilon_{o},\varepsilon_{o})$.
\end{itemize}
\end{theorem}
\par \textbf{Proof. }We can assume that $\widetilde{U}_{o}, \widetilde{V}_{o}$ are, respectively, open balls $B(u_{o};\varepsilon), B(h_{o};\varepsilon)$ for a given $\varepsilon >0$, and we can choose $\delta> 0$ such that if $\lVert u-u_{o}\lVert<\delta$ then $\lVert F(u)-h_{o}\lVert<\varepsilon/2$. Since $u_{o}$ is $k$-transverse then, by Remark 2.4.6 and Theorem 2.5.4 in \cite{B-D 1}, near $u_{o}$ we have:\\
- the inclusion of strata $S_{1_{k+1}}(F)\subseteq S_{1_{k}}(F)$, where $S_{1_{k}}(F)$ is a $k$-codimensional manifold of $X$ in a neighbourhood of $u_{o}$ and $u_{o}\in S_{1_{k}}(F)$; \\
- the disjoint set union 
\begin{center}
$S_{1_{k}}(F)=S_{1_{k+1}}(F)\cup(S_{1_{k}}(F)\setminus S_{1_{k+1}}(F)) ,$
\end{center}
where $S_{1_{k+1}}(F)=\{u\in S_{1_{k}}(F):N(F'(u))\subseteq T_{u}S_{1_{k}}(F) \}$ and $S_{1_{k}}(F)\setminus S_{1_{k+1}}(F)$ is the set of $k$-singularities for $F$ near $u_{o}$.\\
Set $\eta:=\text{min}(\delta,\varepsilon/2)$. Then the following two facts can alternatively occur:\\
j) $ B(u_{o};\eta)\cap (S_{1_{k}}(F)\setminus S_{1_{k+1}}(F))\neq \emptyset$;\\
in such a case, taken $\tilde{u}\in S_{1_{k}}(F)\setminus S_{1_{k+1}}(F)$, we proceed exactly as in the proof of  Theorem \ref{Teo116} to obtain part  i)  of the statement;\\
jj) $B(u_{o};\eta)\cap(S_{1_{k}}(F)\setminus S_{1_{k+1}}(F))= \emptyset$;\\
this means that $B(u_{o};\eta)\cap S_{1_{k}}(F)\subseteq S_{1_{k+1}}(F)$. Let $\varphi$ be a $C^{\infty}$ kernel fibering map for $F$ defined on $B(u_{o};\eta)$. Such a map always exists, up to shrinking $B(u_{o};\eta)$, thanks to Theorem 1.4.6 in \cite{B-D 1}. Given that  $S_{1_{k}}(F)\subseteq S_{1}(F)$ we have by definition that $\varphi(u)$ spans $N(F'(u))$, for any $u\in S_{1_{k}}(F)$, and by the above recalled characterization of $S_{1_{k+1}}(F)$ we get that $\varphi(u)\in T_{u}S_{1_{k}}(F), \, \forall \, u\in B(u_{o};\eta)\cap S_{1_{k}}(F)$. Thus $\varphi$ is a $C^{\infty}$ non-zero vector field on the manifold $B(u_{o};\eta)\cap S_{1_{k}}(F)$. We may integrate such a vector field on $S_{1_{k}}(F)$ near $u_{o}$. More precisely, there exist $\varepsilon_{o} > 0$, an open neighbourhood $U_{o}$ of $u_{o}$ on $B(u_{o};\eta)\cap S_{1_{k}}(F)$ and a $C^{\infty}$ map $\varGamma:(-\varepsilon_{o},\varepsilon_{o})\times U_{o}\rightarrow B(u_{o};\eta)\cap S_{1_{k}}(F)$ such that
\begin{equation*}
\varGamma(0,u)=u \; , \; \dfrac{\partial\varGamma}{\partial t}(t,u)=\varphi(\varGamma(t,u)) \; , \; \forall (t,u) \in (-\varepsilon_{o},\varepsilon_{o})\times U_{o}.
\end{equation*}
The map $\varGamma$ is the local flow of the vector field on the manifold and the map $\varGamma (\cdotp,u)$ is an integral curve of the vector field based at $u$ (cf. \cite{A-M-R} for a detailed proof of the local existence of such a flow). \\
Since $\varGamma(0,u_{o})=u_{o}$, by continuity we can choose $\varepsilon_{o}$ and $U_{o}$ such that $F(\varGamma((-\varepsilon_{o},\varepsilon_{o})\times U_{o}))\subseteq B(h_{o};\varepsilon)$. Let us take $\tilde{u}\in U_{o}$ and define $\tilde{h}:=F(\tilde{u})\in B(h_{o};\varepsilon)$. We contend that $F(\varGamma (t,\tilde{u}))=F(\tilde{u}), \, \forall \, t\in (-\varepsilon_{o}, \varepsilon_{o})$. It is enough to prove that $\dfrac{\partial F\circ\varGamma}{\partial t}(t,\tilde{u})=0, \, \forall \, t\in (-\varepsilon_{o}, \varepsilon_{o})$. Indeed, by applying the chain rule we get 
\begin{equation*}
\dfrac{\partial F\circ\varGamma}{\partial t}(t,\tilde{u})= F'(\varGamma(t,\tilde{u}))\dfrac{\partial \varGamma}{\partial t}(t,\tilde{u})=F'(\varGamma(t,\tilde{u}))\varphi(\varGamma(t,\tilde{u}))=0,
\end{equation*}
because $\varGamma(t,\tilde{u})\in B(u_{o};\eta)\cap S_{1_{k}}(F)\subseteq S_{1}(F)$. Finally, we note that the map $\varGamma(\cdotp,\tilde{u})$ is a \textit{regular curve} since, by construction,$\dfrac{\partial \varGamma}{\partial t}(t,\tilde{u})=\varphi(\varGamma(t,\tilde{u}))\neq 0$. Therefore part ii) of the statement is simply proved by taking $\gamma(t):= \varGamma(t,\tilde{u}), \, \forall \, t\in (-\varepsilon_{o},\varepsilon_{o}). \bracevert$ \\   \vspace{2pt}
\subsection{Some Remarks}\label{ss12}
\subsubsection{}\label{sss121} In Example 2.6.4 of \cite{B-D 1} we saw that, given an integer $k\geq1$ and a $B$-space $Z$, the points $(0,\ldots,0,z), z\in Z$, are $k$-singularities for the map
\begin{align*}
w_{k,Z}:\mathbb{R}^{k}\times Z &\rightarrow \mathbb{R}^{k}\times Z\\
(t,t_{1},\ldots,t_{k-1},z) &\mapsto (t^{k+1}+ t_{k-1}t^{k-1}+\ldots+t_{2}t^{2}+t_{1}t,t_{1},\ldots,t_{k-1}z) ,
\end{align*}
As we already said, this map is called the \textit{generalized Whitney map}.\\
\indent In fact, H. Whitney (cf. \cite{Wh}) studied maps $w:\mathbb{R}^{2}\rightarrow \mathbb{R}^{2}$ of class $C^{d}, d\geq 3$, and considered under which conditions the origin of $\mathbb{R}^{2}$ is, for the map $w$, what we now call a 1-singularity or 2-singularity. He proved that, for suitable changes of coordinates on the domain near 0 and on the target space near $w(0)$, the map $w$ looks like the map $w_{1,\mathbb{R}}(t,s):=(t^{2},s)$ or $w_{2,\{0\}}(t,s):=(t^{3}+st,s)$, respectively. We recall that the map $w_{1,\mathbb{R}}$ is called the \textit{Whitney fold}, while $w_{2,\{0\}}$ is the \textit{Whitney cusp}.\\
\indent A generalization of the Whitney approach, always in the finite-dimensional case, was given by B. Morin in \cite{Mo}. There he considered smooth maps $F:V\rightarrow W$, with $V$ and $W$ smooth manifolds of dimension $n$ and $p$ respectively, $n\leq p$, and studied a class of singular points for $F$ defined by suitable conditions on the \textit{Thom-Boardman strata} of singularities, cf. \cite{Bo}. When $n=p$ this class reduces to the so-called \textit{Morin singularities of order }$k,k\leq n$ (see also \cite{G-G}, chapter VII, for an algebraic approach using the \textit{local ring} of a singularity). In \cite{Mo} it is shown that, near a Morin singularity of order $k$, the map $F$ is locally equivalent to the map that we here call $w_{k,\mathbb{R}^{n-k}}$. More precisely, this result can be restated as the equivalence of the two following conditions:\\
$(\alpha)\; u_{o}\in V $ is a \textit{Morin singularity of order} $k$ for the map $F$;\\
$(\beta)$ there exists a $C^{\infty}$  l.c.d. given by\\
\begin{center}
\begin{tabular}{ccc}
&$F \quad $& \\
$u_{o} \in V\quad $&$ \rightarrow \quad $&$ W$ \bigskip  \\
$\sigma \downarrow $&&$ \downarrow \tau$ \bigskip \\ 
$\underline{0}\in \mathbb{R}^{k} \times \mathbb{R}^{n-k}\quad $&$ \rightarrow  \quad $&$\mathbb{R}^{k} \times \mathbb{R}^{n-k}$, \\
&$w_{k,\mathbb{R}^{n-k}} \quad \medskip $&\\
\end{tabular}
\end{center}
This generalizes Whitney's result, although the diffeomorphisms $\sigma,\tau$ are not found by \textit{ad hoc} procedures like in \cite{Wh} because their existence is deduced by the Malgrange Preparation Theorem (cf. \cite{G-G}, chapter IV).\\
A comparison of the above equivalence with the Normal Form Theorem (Theorem \ref{Teo111}) implies at once that the $k$-singularities in finite dimensions coincide with the Morin singularities of order $k$ in the Euclidean setting, i.e when $V=W=\mathbb{R}^{n}$.
\subsubsection{}\label{sss122} The existence of a normal form for $k$-singularities of $C^{\infty}$ maps was stated in Theorem \ref{Teo111}. Here we are going to show that an analogous result is false for $C^{\infty}$  maximal transverse and $\infty$-transverse singularities.\\
\indent Given the integers $k,n\in \mathbb{N},k\geq 1$, let us consider the map
\begin{align*}
F_{n}:\mathbb{R}^{k+1}&\rightarrow \mathbb{R}^{k+1}\\
(t,t_{1}\ldots,t_{k})&\mapsto((1-\delta_{0n})t^{n}+t_{k}t^{k}+t_{k-1}t^{k-1}+\ldots+t_{2}t^{2}+t_{1}t,t_{1},\ldots,t_{k-1},t_{k}),
\end{align*}
where $\delta_{0n} $ is the Kronecker delta. In Example 2.6.7 of \cite{B-D 1} we saw that the origin $\underline{0}$ is a $k$-maximal transverse singularity  for $F_{n}$, when $n=0$ or $n\geq k+3$. However, the maps $F_{0}$ and $F_{k+3}$  are not $C^{\infty}$ equivalent at the origin. In fact, the equation $F_{0}(t,t_{1},\ldots,t_{k})=\underline{0}$ admits the solution curve given by $\{(t,0,\ldots,0):t\in \mathbb{R}\}$, while $F_{k+3}(t,t_{1},\ldots,t_{k})=\underline{0}$ has the unique solution $\underline{0}$. Hence there cannot exist a normal form for $k$-maximal transverse singularities. Should it exist $F_{0}$ and $F_{k+3}$ would be equivalent at the origin to a same map and thus $F_{0}$ and $F_{k+3}$ would be locally equivalent. This is not possible since curves of solutions are invariant under changes of coordinates.
\indent Let us now study the case of  $\infty$-transverse singularities. For $\mathbb{R}_{1}:=\{t\in \mathbb{R}:\lvert t\lvert < 1\}$, let us consider the map
\begin{align*}
F_{\infty}:\mathbb{R}_{1}\times&\textit{l}^{\;2}(\mathbb{N})\subseteq \mathbb{R}\times\textit{l}^{\;2}(\mathbb{N})\mapsto \mathbb{R}\times\textit{l}^{\;2}(\mathbb{N})\\
&\underline{t}:=(t,t_{1},t_{2},t_{3},\ldots)\mapsto(\sum\limits_{h=1}^{\infty}t_{h}t^{h},t_{1},t_{2},t_{3},\ldots),
\end{align*}
and set $f(\underline{t})=\sum\limits_{h=1}^{\infty}t_{h}t^{h}$. In Example 2.6.9 of \cite{B-D 1} we saw that $\dfrac{\partial^{k}f}{\partial t^{k}}(\underline{0})=0,k\geq 0$ and $\dfrac{\partial^{k+1}f}{\partial t^{k}\partial t_{h}}(\underline{0})$ are l.i., $k,h\geq 1$, that is the origin $\underline{0}$ is an $\infty$-transverse singularity for $F_{\infty}$. Here the set of solutions to $F_{\infty}(\underline{t})=\underline{0}$ is the curve $\{(t,0,0,0,\ldots):t\in \mathbb{R}\}$. On the other hand, we can slightly modify $f$ in order to still have that $\underline{0}$ is an $\infty$-transverse singularity but at the same time to not have a solution curve.\\
For example, let $g:\mathbb{R}\rightarrow \mathbb{R}$ be defined as $g(t)=\text{exp}(-1/t^{2})$, for $t\neq 0,g(0)=0$. Then the $C^{\infty}$ function $g$ is \textit{flat} at $0$, i.e. $\dfrac{\partial^{k}g}{\partial t^{k}}(0)=0,k\geq 0$. If we introduce
\begin{center}
$\hat{f}(\underline{t})=f(\underline{t})+g(t)\, ,$
\end{center}
\begin{flushleft}

\end{flushleft}then $\dfrac{\partial^{k}\hat{f}}{\partial t^{k}}(\underline{0})=\dfrac{\partial^{k}f}{\partial t^{k}}(\underline{0})=0, \; \dfrac{\partial^{k+1}\hat{f}}{\partial t^{k}\partial t_{h}}(\underline{0})=\dfrac{\partial^{k+1}f}{\partial t^{k}\partial t_{h}}(\underline{0})$ are l.i., $k,h \geq 1.$\\
This implies that the origin $\underline{0}$ is an $\infty$-transverse singularity of the $C^{\infty}$ map $\hat{F}:\mathbb{R}_{1}\times\textit{l}^{\;2}(\mathbb{N})\subseteq \mathbb{R}\times\textit{l}^{\;2}(\mathbb{N})\rightarrow \mathbb{R}\times\textit{l}^{\;2}(\mathbb{N})$, defined as $\hat{F}:=\hat{f}\times \textbf{1}_{\textit{l}^{\;2}(\mathbb{N})}$. However, $\underline{0}$ is the unique solution of the equation $\hat{F}(t)=\underline{0}$ since $g(t)=0$ \textit{iff} $t=0$. 
Hence, as above, $F_{\infty}$ and $\hat{F}$ are not equivalent at $\underline{0}$ and thus we cannot have a normal form for such singularities.
\subsubsection{}. \label{sss123} In the above examples concerning $\infty$-transverse singularities we saw that for the analytic map $F_{\infty}$ the set $F_{{\infty}}^{-1}(\underline{0})$ is a curve passing through $\underline{0}$. It is possible to prove that such a result is true for all analytic maps having an $\infty$-transverse singularity, see Subsection \ref{sss131} below. On the other hand for the map $\hat{F}$, which is of class $C^{\infty}$ and is not analytic, the set $\hat{F}^{-1}(\underline{0})$ is a single point but it is not difficult to provide examples of maps with the same regularity and such that the pre-image of $\underline{0}$ is a sequence of points or a non-analytic curve, as we will see in Subsections \ref{sss132} and \ref{sss133} respectively.\\
\par
We conclude the section with two simple topological addenda to the Local Multiplicity Theorem, i.e. Theorem \ref{Teo115}, valid for integers $k$ odd or even, respectively. We start by giving a rather trivial analytical complement to Lemma \ref{Lem114}.
\subsubsection{}\label{sss124}(Analytical addendum to Lemma \ref{Lem114}).\textit{ Let } $k\geq 1$\textit{ be a fixed integer, }$\varepsilon >0$\textit{ a real number, }$\alpha_{0},\ldots,\alpha_{k-1}\in \mathbb{R}$ \textit{ and }$\lvert\alpha_{i}\lvert<\varepsilon,i=0,\ldots,k-1:$\textit{ there exists a real number }$\rho(\varepsilon)>0$\textit{ such that }$\rho(\varepsilon)\rightarrow 0$\textit{ if }$\varepsilon\rightarrow 0$\textit{ and such that whenever }$x_{0}$\textit{ is a real root of the polynomial }$x^{k+1}+\alpha_{k-1}x^{k-1}+\ldots+\alpha_{0}$\textit{ then }$ \lvert x_{0}\lvert < \rho(\varepsilon)$. \smallskip \\
\indent Proof. Let us consider $x_{0}^{k+1}+\alpha_{k-1}x_{0}^{k-1}+\ldots+\alpha_{0}=0$ with $\lvert \alpha_{i}\lvert <\varepsilon, i=0,\ldots, k-1$. Then $ \lvert x_{0}\lvert^{k+1} \leq\lvert \alpha_{k-1}\lvert \lvert x_{0}\lvert^{k-1}+\ldots+\lvert\alpha_{0}\lvert \leq \varepsilon(\lvert x_{0}\lvert^{k-1}+\ldots+1)$ and so $ g(\lvert x_{0}\lvert)\leq \varepsilon$, where $g(t):=\dfrac{t^{k+1}}{t^{k}+\ldots+1}$, for $t\geq 0$. It is easy to see that $g'(t)>0$, hence $g$ is increasing and invertible: this implies that $\lvert x_{0}\lvert \leq g^{-1}(\varepsilon)$. We can thus choose $\rho(\varepsilon):=g^{-1}(\varepsilon)$: since $ g(0)=0$, if $\varepsilon\rightarrow 0$ then $\rho(\varepsilon)\rightarrow 0.$ \\
\indent 
We notice that one can give a rough estimate of $\rho(\varepsilon)$ as max$\{(k\varepsilon)^{1/2},(k\varepsilon)^{1/k+1}\}$.
\subsubsection{}\label{sss125}(Topological addendum to Theorem \ref{Teo115} for $k$ odd).\textit{ Let }$U,V$ \textit{be open subsets of the B-spaces} $X,Y$\textit{ and }$F:U\subseteq X \rightarrow V \subseteq Y$\textit{ a }${C^{\infty}}$ 0-\textit{Fredholm map. Let }$u_{o}\in U$\textit{ be a }$k$-singularity\textit{ for }$F$,\textit{ for some }$k\geq 1,k$\textit{ odd. Then there exist neighbourhoods }$U_{o}\subseteq U$\textit{ of }$u_{o}$\textit{ and }$V_{o}\subseteq V$\textit{ of }$h_{o}=F(u_{o})$\textit{ and an open subset }$\varOmega_{o}\subseteq V_{o}$\textit{ such that:} \\
-$\;\,\forall \, h\in \varOmega_{o}$,\textit{ the equation }$F(u)=h$\textit{ has no solutions }$u\in U_{o}$;\\
-$\;\, h_{o}\in \partial \varOmega_{o} $, \textit{ and in particular }$\varOmega_{o}\neq \emptyset$.\smallskip \\
\indent Proof. By virtue of the Normal Form Theorem, see also l.c.d. (\ref{112}), it is sufficient to give the proof for a generalized Whitney map $w_{k,Z}$, with $Z$ a suitable $B$-space. Since such a map is the identity on $Z$ there is no loss of generality to give the proof of \ref{sss125} for the map $w_{k}:U_{o}\subseteq \mathbb{R}^{k}\rightarrow V_{o}\subseteq \mathbb{R}^{k}$, with $U_{o},V_{o}$ arbitrary open subsets such that $u_{o}=\textbf{0}\in U_{o},h_{o}=\textbf{0}\in V_{o}$ and $w_{k}(U_{o})\subseteq V_{o}$, where we recall that $w_{k}(t,t_{1},\ldots,t_{k-1})=(t^{k+1}+t_{k-1}t^{k-1}+\ldots+t_{2}t^{2}+t_{1}t,t_{1},\ldots,t_{k-1})$.\\
\indent Let us choose $\varepsilon> 0$ such that $B(\textbf{0};\varepsilon)\subseteq V_{o}$ and cl$(B(\textbf{0};\delta))\subseteq U_{o}$, where we define $\delta:=\sqrt{k} \text{ max} (\varepsilon,\rho(\varepsilon))$ with $\rho(\varepsilon)$ given by \ref{sss124}(cl$(\cdotp)$ means topological closure). Since $\rho(\varepsilon)\rightarrow 0$ if $\varepsilon \rightarrow 0$, we can always choose $\varepsilon$ and $\delta$ as such; moreover, such a choice implies that $w_{k}^{-1}(B(\textbf{0};\varepsilon))\subseteq B(\textbf{0};\delta)$. In fact, if $\textbf{t}:=(t,t_{1},\ldots,t_{k-1}),\textbf{s}:=(s,s_{1}, ...,s_{k-1})\in \mathbb{R}^{k}$ we  already saw that $w_{k}(\textbf{t})=\textbf{s}$ implies that $\textbf{t}=(t,s_{1}, \ldots,s_{k-1})$ with $t$ solving the equation $t^{k+1}+s_{k-1}t^{k-1}+\ldots+ s_{1}t-s=0$. If $\textbf{s}\in B(\textbf{0};\varepsilon)$ it follows that $ \lvert s\lvert,\lvert s_{1}\lvert,\ldots,\lvert s_{k-1}\lvert< \varepsilon$ and thus the last equality and \ref{sss124} imply that $\lvert t_{n}\lvert< \rho(\varepsilon)$: it is then easy to see that, by definition of $\delta,\lVert \textbf{t} \lVert < \delta$ i.e $\textbf{t}\in B(\textbf{0};\delta).$
\par Let us set $\varOmega_{o}:=\{\textbf{s}\in B(\textbf{0};\varepsilon):w_{k}(\textbf{t})=\textbf{s}\text{ has no solutions }\textbf{t}\in U_{o}\}$. We notice that $\textbf{0}\notin \varOmega_{o}$ because $w_{k}(\textbf{0})=\textbf{0}$; furthermore, by choosing $\textbf{s}_{n}:=(-1/n^{2},0,\ldots,0)$ for $n\geq 1$ the equation $w_{k}(\textbf{t})=\textbf{s}_{n}$ has globally no solutions because this is equivalent to solving $t^{k+1}=-1/n^{2}$ for $k+1$ an even integer. Both facts imply that $\textbf{0}\in \partial\varOmega_{o}$. \\
We now show that $\varOmega_{o}$ is open in $U_{o}$. Were this not true, there would exist $ \textbf{s}_{o}=(s,s_{1},\ldots,s_{k-1})\in B(\textbf{0};\varepsilon)$ and a sequence $\textbf{s}_{n}\in B(\textbf{0};\varepsilon), n\geq 1$, such that $\textbf{s}_{n}\rightarrow \textbf{s}_{o},w_{k}(\textbf{t})=\textbf{s}_{o}$ has no solutions in $U_{o}$ while $w_{k}(\textbf{t})=\textbf{s}_{n}$ has at least one solution $\textbf{t}_{n}\in U_{o}$ with  $\textbf{t}_{n}\in B(\textbf{0};\delta)$. Since $B(\textbf{0};\delta)$ is bounded, we can suppose without loss of generality that $\textbf{t}_{n}\rightarrow \textbf{t}\in$ cl$(B(\textbf{0};\delta))\subseteq U_{o}$. By continuity, $w_{k}(\textbf{t}_{n})=\textbf{s}_{n}\rightarrow \textbf{s}_{o}=w_{k}(\textbf{t})$ and this contradicts the choice of $\textbf{s}_{o}$. 
\subsubsection{}\label{sss126}(Topological addendum to Theorem \ref{Teo115} for $k$ even).\textit{ Let }$U,V$ \textit{be open subsets of the B-spaces} $X,Y$\textit{ and }$F:U\subseteq X \rightarrow V \subseteq Y$\textit{ a }${C^{\infty}}$ 0-\textit{Fredholm map. Let }$u_{o}\in U$\textit{ be a }$k$-singularity\textit{ for }$F$,\textit{ for some }$k\geq 1,k$\textit{ even. Then there exist neighbourhoods }$U_{o}\subseteq U$\textit{ of }$u_{o}, V_{o}\subseteq V$\textit{ of }$h_{o}=F(u_{o})$\textit{ and open subsets }$\varOmega_{1}\subseteq \widetilde\varOmega \subseteq V_{o}$\textit{ such that:} \\
-$\;\,\forall \, h\in \widetilde{\varOmega}$,\textit{ the equation }$F(u)=h$\textit{ has at least a solution }$u\in U_{o}$;\\
-$\;\,\forall \, h\in \varOmega_{1} $,\textit{ the equation }$F(u)=h$\textit{ has a unique solution }$u\in U_{o}$\textit{ and }$u$\textit{ is regular for }$F$;\\
-$\;\, h_{o}\in \partial \varOmega_{1}$,\textit{ and in particular }$\varOmega_{1}\neq \emptyset$.\textit{ Precisely, for any pair of neighbourhoods }$
\widetilde{U}_{o}\subseteq U$\textit{ of }$u_{o}$\textit{ and }$\widetilde{V}_{o}\subseteq V$\textit{ of }$h_{o}$\textit{ there exists }$\tilde{h}\in \widetilde{V}_{o}$\textit{ such that the equation }$F(u)=\tilde{h}$\textit{ has exactly one solution }$\tilde{u}_{1}$\textit{ in }$\widetilde{U}_{o}$\textit{ with }$\tilde{u}_{1}$\textit{ regular for }$F$. \smallskip \\
\indent Proof. Once again, it is sufficient to give the proof for the map $w_{k}:U_{o}\subseteq \mathbb{R}^{k}\rightarrow V_{o}\subseteq \mathbb{R}^{k}$, with $U_{o},V_{o}$, arbitrary open subsets such that $u_{o}=\textbf{0}\in U_{o}, h_{o}=\textbf{0}\in V_{o}$ and $w_{k}(U_{o})\subseteq V_{o}$. We start by choosing $\varepsilon$ and $\delta$ as in \ref{sss125}, i.e. such that $B(\textbf{0};\varepsilon)\subseteq V_{o}$, cl$(B(\textbf{0};\delta))\subseteq U_{o}$ and $w_{k}^{-1}(B(\textbf{0};\varepsilon))\subseteq B(\textbf{0};\delta)$. Let us define $\widetilde{\varOmega}:=B(\textbf{0};\varepsilon)$ and $\varOmega_{1}:= \{\textbf{s}\in B(\textbf{0};\varepsilon): w_{k}(\textbf{t})=\textbf{s}$\text{ has only one solution }$\textbf{t}\in U_{o}\text{ and it is regular}\}$. Since the polynomial $t^{k+1}+ t_{k-1}t^{k-1}+\ldots+t_{2}t^{2}+t_{1}t$ has odd degree $k+1$ then the equation $w_{k}(\textbf{t})=\textbf{s}$ has at least one solution $\textbf{t}$ for any given $\textbf{s}$: in particular, if $\textbf{s}\in B(\textbf{0};\varepsilon)$ then $\textbf{t}\in B(\textbf{0};\delta)\subseteq U_{o}$. Now, let us suppose by contradiction that $\varOmega_{1}$ is not open, and so there exists $\textbf{s}_{o}\in \varOmega_{1}$ and a sequence $\textbf{s}_{n}\in B(\textbf{0};\varepsilon), n\geq 1$, such that $\textbf{s}_{n}\rightarrow \textbf{s}_{o}$ and $\textbf{s}_{n}\notin \varOmega_{1}$. Thus $w_{k}(\textbf{t})=\textbf{s}_{o}$ has a unique solution $\textbf{t}_{o}\in U_{o}$, namely $\textbf{t}_{o}\in B(\textbf{0};\delta)$, and $\textbf{t}_{o}$ is regular: by the Inverse Function Theorem there exist open neighbourhoods $\hat{U}\subseteq B(\textbf{0};\delta)$ of $\textbf{t}_{o}$ and $\hat{V}\subseteq B(\textbf{0};\varepsilon)$ of $\textbf{s}_{o}$ such that $w_{k}(\hat{U})= \hat{V}$ and $w_{k}:\hat{U}\rightarrow \hat{V}$ is a diffeomorphism. We can assume that $\textbf{s}_{n}\in \hat{V}$ for all $n$, therefore  $w_{k}(\textbf{t})=\textbf{s}_{n}$ has a unique regular solution $\textbf{t}_{n}\in \hat{U}$. Since $\textbf{s}_{n}\notin \varOmega_{1}$, this means that $w_{k}(\textbf{t})=\textbf{s}_{n}$ has another solution $\textbf{z}_{n}\in U_{o}$: precisely, $\textbf{z}_{n}\in B(\textbf{0};\delta)\setminus \hat{U}$. Hence, up to subsequences, $\textbf{z}_{n}\rightarrow \textbf{z}_{o}\in$ cl$(B(\textbf{0};\delta))\setminus \hat{U}\subseteq U_{o}$. By continuity we get that $\textbf{z}_{o}$ solves $w_{k}(\textbf{z}_{o})=\textbf{s}_{o}$ with $\textbf{z}_{o}\neq \textbf{t}_{o}$, and this is absurd.\\
\indent To prove the last part let $\widetilde{U}_{o} \subseteq U_{o}$ and  $\widetilde{V}_{o} \subseteq U_{o}$ be neighbourhoods of $\textbf{0}\in \mathbb{R}^{k}$. We define $\textbf{s}_{n}= (1/n^{k+1},0,\ldots,0),n\geq 1$, so the equation $w_{k}(\textbf{t})=\textbf{s}_{n}$ has a unique solution $\textbf{t}_{n}=(1/n,0,\ldots,0)$ and it is not difficult to check that $\textbf{t}_{n}$ is regular for $w_{k}$. Since, for large $n,\,\textbf{s}_{n}\in \widetilde{V}_{o}$ and $\textbf{t}_{n}\in \widetilde{U}_{o}$ the result is proved. Eventually, we notice that $w_{k}(\textbf{t})=\textbf{0}$ has a unique non-regular solution given by $\textbf{t}=\textbf{0}$ which implies that $\textbf{0}\notin \varOmega_{1}$: this, combined with what we found above, proves that $\textbf{0}\in \partial \varOmega_{1}.$\\
\par We finally remark that for $k$ odd, say $k=2n+1$, we could locally exhibit, for all $h=0, 2,\ldots,2n,2n+2$, non-empty open regions $\varOmega_{h}$ with $h$ regular solutions. Analogously, for $k$ even, say $k=2m$, and for all $h=1,3,\ldots,2m-1,2m+1$, we could show that there exist non-empty open regions $\varOmega_{h}$ with $h$ regular solutions. Here we restricted our study to the cases $k$ odd, $h=0,2n+2$ (cf. \ref{sss125} above and Theorem \ref{Teo115}) and $k$ even, $h=1,2m+1$ (cf. \ref{sss126} above and Theorem \ref{Teo115}).\\ \vspace{2pt}
\subsection{Properties of $\infty$-transverse Singularities}\label{ss13}
\subsubsection{}\label{sss131}. Let $U,V$  be open subsets of the $B$-spaces $X,Y$ and let $F:U\subseteq X \rightarrow V \subseteq Y$ be a $C^{\omega}$ 0-Fredholm map with $u_{o}\in U$ an $\infty$-transverse singularity for $F$. We claim that there exists a solution curve passing through $u_{o}$, precisely there exist $\varepsilon_{o}>0$ and a $C^{\omega}$ map $\gamma:(-\varepsilon_{o},\varepsilon_{o})\rightarrow U \subseteq X$ such that $\gamma\,'(t)\neq 0$ and $F(\gamma(t))=F(u_{o}),\,\forall \, t \in (-\varepsilon_{o},\varepsilon_{o}).$ \\
In order to show such a property we note that analytic curves of solutions are invariant under analytic changes of coordinates. Then, by means of the Local  Representation Theorem (cf. Theorem 1.2.2 in \cite{B-D 1}) we can assume that the map we are studying  has the following form:
\begin{align*}
\varPhi:(-\varepsilon_{o},\varepsilon_{o}) \times W \subseteq \mathbb{R} \times \varXi &\rightarrow (-\varepsilon_{o},\varepsilon_{o})\times W \subseteq \mathbb{R} \times \varXi\\
\varPhi(t,\xi)&=(f(t,\xi),\xi), \, \forall \, (t,\xi)\in (-\varepsilon_{o},\varepsilon_{o})\times W,
\end{align*}
where $W$ is an open subset of a suitable $B$-spaces $\varXi$ and $f:U \subseteq \mathbb{R} \times \varXi \rightarrow \mathbb{R}$ is a $C^{\omega}$ function such that $f(0,0)=0$. In other words, we study an LS-map (see Definition 1.2.4 in \cite{B-D 1}) locally equivalent to the given map. Finally, without loss of generality, we can still assume that $(0,0)\in \mathbb{R} \times \varXi$ is an $\infty$-transverse singularity for $\varPhi$. This is because of the  Invariance Theorem (see Theorem 2.5.5 in \cite{B-D 1}), which assures us that the considered type of singularity is unaffected by changes of coordinates.\\
Since $(0,0)$ is an $\infty$-transverse singularity for $\varPhi$ and thus $k$-transverse for any integer $k$ we have, by Proposition 2.6.2(T) of \cite{B-D 1}, that $\dfrac{\partial^{h+1}f}{\partial t^{h+1}}(0,0)=0,\, \forall h\geq 0$. Let us now consider the analytic function $\phi:(-\varepsilon_{o},\varepsilon_{o}) \rightarrow \mathbb{R},\,\phi(t):=f(t,0)$. Given that $\dfrac{\partial^{h+1}\phi}{\partial t^{h+1}}(0)=0$ for any integer $h$ we obtain, from the analyticity of $\phi$, that $\phi(t)\equiv \phi(0)$ and thus that $f(t,0)\equiv f(0,0)=0$. Therefore  $\varPhi(t,0)\equiv (f(t,0),0)=(0,0)$, and $\gamma:(-\varepsilon_{o},\varepsilon_{o}) \rightarrow \mathbb{R} \times \varXi$, defined as $ \gamma(t):=(t,0),\,t \in (-\varepsilon_{o},\varepsilon_{o})$, is a $C^{\omega}$ curve of solutions to the equation $\varPhi(t,\xi)=(0,0)$.\\
\indent In conclusion, we have shown that the local existence of a solution curve is a geometrical property common to all analytic 0-Fredholm maps near an $\infty$-transverse singularity. We can give a simple example of such a behaviour by considering the $C^{\omega}$ map
\begin{align*}
F_{\infty,Z}:\mathbb{R}_{1} \times \textit{l}\,^{2}(\mathbb{N})\times Z \subseteq \mathbb{R}\times \textit{l}\,^{2}(\mathbb{N})\times Z&\rightarrow
\mathbb{R} \times \textit{l}\,^{2}(\mathbb{N})\times Z\\
(t,t_{1},t_{2},t_{3}, \ldots , z) &\mapsto (\sum_{h=1}^{\infty} t_{h}t^{h},t_{1},t_{2},t_{3}, \ldots , z).
\end{align*}
Indeed, referring to Subsection \ref{sss122}, one can easily show that the origin $\underline{0}$ is an $\infty$-transverse singularity for $F_{\infty,Z}$ and that the straight line $\{(t,0):t \in \mathbb{R}_{1}\}$ is a solution set for the equation $F_{\infty,Z}(\underline{0})=\underline{0}$.\\
\indent The ``polynomial'' form of the map $F_{\infty,Z}$, similar to that of Whitney maps, prompts the question as to whether the map $F_{\infty,Z}$ could be the Normal Form of $C^{\omega}$, 0-Fredholm maps $F:X\rightarrow X$ near an $\infty$-transverse singularity, at least in the simple case where X is a \textit{separable} Hilbert space. While this is still an open problem, we know that in the general case, with $F$ a $C^{\infty}$ map, the answer is negative, i.e. there is not a normal form for $\infty$-transverse singularities (see Subsection \ref{sss122}).
\subsubsection{}.\label{sss132} Let us consider again the map $F_{\infty}$ studied in Subsection \ref{sss122}, i.e.
\begin{align*}
F_{\infty}:\mathbb{R}_{1} \times \textit{l}\,^{2}(\mathbb{N}) \subseteq \mathbb{R}\times \textit{l}\,^{2}(\mathbb{N})&\rightarrow \mathbb{R} \times \textit{l}\,^{2}(\mathbb{N})\\
\underline{t}:=(t,t_{1},t_{2},t_{3}, \ldots) &\mapsto (f(\underline{t}),t_{1},t_{2},t_{3}, \ldots),
\end{align*}
where $\mathbb{R}_{1}:=\{t\in \mathbb{R}:\lvert t\lvert < 1\}$ and $f(\underline{t}):=\sum_{h=1}^{\infty} t_{h}t^{h}$. We recall that $\dfrac{\partial^{k}f}{\partial t^{k}}(\underline{0})=0,k\geq 0$, and  $\dfrac{\partial^{k+1}f}{\partial t^{k}\partial t_{h}}(\underline{0})$ are l.i., $k,h\geq 1$. Let us define the map $\bar{g}: \mathbb{R}\rightarrow \mathbb{R}$ as $\bar{g}(t):=\text{exp}(-1/t^{2})\text{sin}(1/t)$ for $t\neq 0$ and $\bar{g}(0):=0$. Then one gets $\dfrac{\partial^{k}\bar{g}}{\partial t^{k}}(0)=0,k\geq 0$, and thus, if we set $\bar{f}(\underline{t}):=f(\underline{t})+\bar{g}(t)$, it is easy to see that $\dfrac{\partial^{k}\bar{f}}{\partial t^{k}}(\underline{0})=\dfrac{\partial^{k}f}{\partial t^{k}}(\underline{0})$ and $\dfrac{\partial^{k+1}\bar{f}}{\partial t^{k}\partial t_{h}}(\underline{0})=\dfrac{\partial^{k+1}f}{\partial t^{k}\partial t_{h}}(\underline{0}),k,h\geq 1$. Hence one has that $\dfrac{\partial^{k}\bar{f}}{\partial t^{k}}(\underline{0})=0$ and  $\dfrac{\partial^{k+1}\bar{f}}{\partial t^{k}\partial t_{h}}(\underline{0})$ are l.i., for any $k,h\geq 1$; this implies that the origin $\underline{0}$ is an $\infty$-transverse singularity for
\begin{center}
$\bar{F}:\mathbb{R}_{1} \times \textit{l}\,^{2}(\mathbb{N}) \subseteq \mathbb{R}\times \textit{l}\,^{2}(\mathbb{N})\rightarrow \mathbb{R} \times \textit{l}\,^{2}(\mathbb{N}), \bar{F}(t,t_{1},t_{2},t_{3}, \ldots):= (\bar{f}(\underline{t}),t_{1},t_{2},t_{3}, \ldots).$
\end{center}
Moreover, we have that $(\bar{f}(\underline{t}),t_{1},t_{2},t_{3},\ldots)=(0,0,0,0,\ldots)$ \textit{iff} $t_{1}=t_{2}=t_{3}=\ldots=0$ and $\bar{g}(t)=0$, i.e. $t=0$ or $t=1/m\pi,m\in \mathbb{Z}\setminus \{0\}$. Therefore the solution set to $\bar{F}(t,t_{1},t_{2},t_{3},\ldots)=\underline{0}$ is given by the sequence $\{\underline{0}\}\cup\{(1/m\pi,0,0,0,\ldots),m\in \mathbb{Z}\setminus \{0\}\}$.
\subsubsection{}\label{sss133}. In order to find a non-analytic curve of solutions near an $\infty$-transverse singularity for a $C^{\infty}$ map which is not $C^{\omega}$ it is not convenient to consider an LS-map as we did in the previous examples. Hence we shall study a non-analytic perturbation on the first and the second coordinate of the map $F_{\infty}$ given in Subsection \ref{sss122}. For this it suffices to use again the function $g:\mathbb{R}\rightarrow \mathbb{R}$ defined as $g(t)=\text{exp}(-1/t^{2})$ for $t\neq 0, g(0)=0$. More precisely, we define 
\begin{align*}
\widetilde{F}:\mathbb{R}_{1} \times \textit{l}\,^{2}(\mathbb{N}) \subseteq \mathbb{R}\times \textit{l}\,^{2}(\mathbb{N})&\rightarrow \mathbb{R} \times \textit{l}\,^{2}(\mathbb{N})\\
\underline{t}:=(t,t_{1},t_{2},t_{3}, \ldots) &\mapsto (f(\underline{t})+g(t)t,t_{1}+g(t),t_{2},t_{3}, \ldots)
\end{align*}
with $f(\underline{t})=\sum_{h=1}^{\infty} t_{h}t^{h}$.\\
It is useful to write $\widetilde{F}(\underline{t})=(\tilde{f}_{0}(\underline{t}),\tilde{f}_{1}(\underline{t}),\tilde{f}_{2}(\underline{t}),\tilde{f}_{3}(\underline{t}),\ldots)$ where $\tilde{f}_{0}(\underline{t}):=(t_{1}+g(t))t + \sum_{h=2}^{\infty} t_{h}t^{h},\tilde{f}_{1}(\underline{t}):=t_{1}+g(t)$ and $\tilde{f}_{h}(\underline{t}):=t_{h}$, for $h\geq 2$.\\
One can easily see that the solutions to $\widetilde{F}(\underline{t})=\underline{0}$ are of the form $\{(t,-g(t),0,0,\ldots), t\in \mathbb{R}_{1}\}$, which is just a $C^{\infty}$ and non-analytic curve passing through $\underline{0}$.\\
This means that we only need to check that $\widetilde{F}$ is a 0-Fredholm map and $\underline{0}$ is an $\infty$-transverse singularity for $\widetilde{F}$. We shall sketch the direct proof below while at the end we will give a diagram-wise justification of this property.\\
It is clear that $\widetilde{F}^{'}(\underline{t}):\mathbb{R}\times \textit{l}\,^{2}(\mathbb{N})\rightarrow \mathbb{R} \times \textit{l}\,^{2}(\mathbb{N})$ can be represented by an infinite "matrix" with entries $\dfrac{\partial \tilde{f}_{h}}{\partial t_{\eta}}(\underline{t}), h,\eta\geq 0$ where $t_{0}\equiv t$. Hence
\begin{equation*}
\widetilde{F}'(\underline{t})= \left[ \begin{array}{cccc} 
\dfrac{\partial f}{\partial t}(\underline{t})+g(t)+g'(t)t &\quad t &\quad t^{2} &\quad \cdots  \bigskip\\
g'(t) &\quad 1 &\quad 0 &\quad\cdots \bigskip \\
0 &\quad 0 &\quad 1 &\bigskip \\
\vdots &\quad\vdots &\quad &\quad\ddots \medskip
\end{array} \right].
\end{equation*}
Then, by considering the second-order minor 
\begin{equation*}
M(\underline{t}):= \left[ \begin{array}{cc} 
\dfrac{\partial f}{\partial t}(\underline{t})+g(t)+g'(t)t &\quad t \bigskip\\
g'(t) &\quad 1  \medskip
\end{array} \right],
\end{equation*}
one gets that the point $\underline{t}$ is singular for $\widetilde{F}$ \textit{iff} $M$ is a singular matrix, i.e.
\begin{center}
det$\,M(\underline{t})=\dfrac{\partial f}{\partial t}(\underline{t})+g(t) = 0$.
\end{center}
In this case it is not difficult to verify that $(1,-g'(t))$ spans $N(M(\underline{t}))$. Therefore, when $\dfrac{\partial f}{\partial t}(\underline{t})+g(t)=0$ then dim$\,N(\widetilde{F}'(\underline{t}))=1$ and $\varphi(\underline{t}):=(1,-g'(t),0,0,\ldots)$ is a fibering map for the kernel of $\widetilde{F}'(t)$. In particular $\underline{0} \in S_{1}(\widetilde{F})=\{\underline{t} \in \mathbb{R}_{1} \times \textit{l}\,^{2}(\mathbb{N}):\text{det}\,M(\underline{t})=0\}$.\\
Additionally, with the usual identifications $\mathbb{R}^{\ast}\cong \mathbb{R}$ and $\textit{l}\,^{2}(\mathbb{N})^{\ast}\cong \textit{l}\,^{2}(\mathbb{N})$ the adjoint $\widetilde{F}'(\underline{t})^{\ast}:(\mathbb{R} \times \textit{l}\,^{2}(\mathbb{N}))^{\ast}\cong \mathbb{R} \times \textit{l}\,^{2}(\mathbb{N})\rightarrow (\mathbb{R} \times \textit{l}\,^{2}(\mathbb{N}))^{\ast}\cong \mathbb{R} \times \textit{l}\,^{2}(\mathbb{N})$ is given by the \textquotedblleft transpose\textquotedblright $\,$ of  $\widetilde{F}'(\underline{t})$, that is
\begin{equation*}
\widetilde{F}'(\underline{t})^{\ast}= \left[ \begin{array}{cccc} 
\dfrac{\partial f}{\partial t}(\underline{t})+g(t)+g'(t)t &\quad g'(t) &\quad 0 &\quad \cdots  \bigskip\\
t &\quad 1 &\quad 0 &\quad\cdots \bigskip \\
t^{2} &\quad 0 &\quad 1 &\bigskip \\
\vdots&\quad\vdots &\quad &\quad\ddots\medskip
\end{array} \right].
\end{equation*}
By considering
\begin{equation*}
M^{\ast}(\underline{t}):= \left[ \begin{array}{cc} 
\dfrac{\partial f}{\partial t}(\underline{t})+g(t)+g'(t)t &\quad g'(t) \bigskip\\
t &\quad 1  \medskip
\end{array} \right]
\end{equation*}
we get that $\widetilde{F}'(\underline{t})^{\ast}$ is singular \textit{iff} $\text{det}\,M^{\ast}(\underline{t}) =\text{det}\,M(\underline{t})=0$, that is $\underline{t} \in S_{1}(\widetilde{F})$. One can also check that the vector $\psi(\underline{t}):=(1,-t,-t^{2},-t^{3},\ldots)$ spans $N(\widetilde{F}'(\underline{t})^{\ast})$ for any $\underline{t}\in S_{1}(\widetilde{F})$. Thus dim$N(\widetilde{F}'(\underline{t})^{\ast})=1$ and $\psi(\underline t)$ is a cokernel fibering map for $\widetilde{F}$, see Definition 1.1.3 in \cite{B-D 1}. The above properties show that $\widetilde{F}$ is a 0-Fredholm map too.\\
By using Definition 1.1.5 in \cite{B-D 1}, and after some calculations, we obtain that\\
$J_{0}(\varphi,\psi)(\underline{t})=\psi(\underline{t})\widetilde{F}'(\underline{t})\varphi(\underline{t})=\dfrac{\partial f}{\partial t}(\underline{t})+g(t);$\\
$I_{1}(\varphi,\psi)(\underline{t})(r,v)=J'_{0}(\varphi,\psi)(\underline{t})(r,v)=(\dfrac{\partial f}{\partial t})'(\underline{t})(r,v)+rg'(t), r \in \mathbb{R}, v \in \textit{l}\,^{2}(\mathbb{N});$\\
$J_{1}(\varphi,\psi)(\underline{t})=I_{1}(\varphi,\psi)(\underline{t})\varphi(\underline{t})=I_{1}(\varphi,\psi)(\underline{t})(1,-g'(t),0,0,\ldots)=\\
=1\cdotp \dfrac{\partial^{2}f}{\partial t^{2}}(\underline{t})-g'(t)\cdotp \dfrac{\partial^{2}f}{\partial t \partial t_{1}}(\underline{t})+1\cdotp g'(t)=\dfrac{\partial^{2}f}{\partial t^{2}}(\underline{t})$ because $\dfrac{\partial^{2}f}{\partial t \partial t_{1}}(\underline{t})=1;$\\
$I_{2}(\varphi,\psi)(\underline{t})(r,v)=J'_{1}(\varphi,\psi)(\underline{t})(r,v)=(\dfrac{\partial^{2}f}{\partial t^{2}})'(\underline{t})(r,v);$\\
$J_{2}(\varphi,\psi)(\underline{t})=I_{2}(\varphi,\psi)(\underline{t})(1,-g'(t),0,0,\ldots)=1\cdotp \dfrac{\partial^{3}f}{\partial t^{3}}(\underline{t})-g'(t)\cdotp \dfrac{\partial^{3}f}{\partial t^{2} \partial t_{1}}(\underline{t})=\dfrac{\partial^{3}f}{\partial t^{3}}(\underline{t})$, since $\dfrac{\partial^{3}f}{\partial t^{2}\partial t_{1}}(\underline{t})=0$.\\
It is now clear how to continue in determining the functionals $I_{h},J_{h}$. In fact, it is easy to check that:
\begin{equation*}
J_{h}(\varphi,\psi)(\underline{0})=\dfrac{\partial^{h+1}f}{\partial t^{h+1}}(\underline{0})\, ,\, h\geq 0, I_{h}(\varphi,\psi)(\underline{0})=(\frac{\partial^{h}f}{\partial t^{h}})'(\underline{0}), h\geq 1.
\end{equation*}
These equalities imply that $\underline{0}$ is an $\infty$-transverse singularity for $\widetilde{F}$ because $\dfrac{\partial^{k} f}{\partial t^{k}}(\underline{0})=0,k\geq 0$, and  $\dfrac{\partial^{k+1} f}{\partial t^{k}\partial t_{j}}(\underline{0}),k,j \geq 1$, are l.i., as we recalled above.\\
\indent There is also an alternative proof of the fact that $\underline{0}$ is an $\infty$-transverse singularity for $\widetilde{F}$ which relies on the use of a suitable l.c.d.. Let us define the map $\alpha:\mathbb{R} \times \textit{l}\,^{2}(\mathbb{N})\rightarrow \mathbb{R} \times \textit{l}\,^{2}(\mathbb{N})$ as $\alpha(t,t_{1},t_{2},t_{3}, \ldots)= (t,t_{1}-g(t),t_{2},t_{3}, \ldots)$ and let us consider the diagram 
\qquad \begin{center}
\begin{tabular}{ccc}
&$F_{\infty} \quad $& \\
$\underline{0} \in \mathbb{R}_{1} \times \textit{l}\,^{2}(\mathbb{N}) \quad $&$ \rightarrow \quad $&$\mathbb{R} \times \textit{l}\,^{2}(\mathbb{N})$ \bigskip  \\
$\alpha \downarrow $&&$ \downarrow \beta$ \bigskip \\ 
$\underline{0}\in \mathbb{R}_{1} \times \textit{l}\,^{2}(\mathbb{N}) \quad $&$ \rightarrow  \quad $&$ \mathbb{R} \times \textit{l}\,^{2}(\mathbb{N})$ ,\\
&$\widetilde{F} \quad $&\\
\end{tabular}
\end{center}
where $\beta$ is the identity. It is easy to check that $\alpha$ is a diffeormorphism and the diagram is commutative. Since $\underline{0}$ is an $\infty$-transverse singularity for the map $F_{\infty}$ then, by the Invariance Theorem (cf \cite{B-D 1}, Theorem 2.5.5), $\alpha(\underline{0})= \underline{0}$ is an $\infty$-transverse singularity for $\widetilde{F}$.\\
As a matter of fact the straight line $\{(t,0):t \in \mathbb{R}_{1}\}$, which is an analytic solution curve for $F_{\infty}$, through the non-analytic diffeormorphism $\alpha$ is mapped into the non-analytic solution curve $\{(t,-g(t),0,0,\ldots), t\in \mathbb{R}_{1}\}$ for the map $\widetilde{F}$.
\subsubsection{}\label{sss134}. Let us consider again the maps $F_{\infty}$ and $g$ studied above and, for any $\varepsilon>0$, let us define the map 
\begin{align*}
F_{\infty,\varepsilon}:\mathbb{R}_{1}\times&\textit{l}^{\;2}(\mathbb{N})\subseteq \mathbb{R}\times\textit{l}^{\;2}(\mathbb{N})\mapsto \mathbb{R}\times\textit{l}^{\;2}(\mathbb{N})\\
&\underline{t}:=(t,t_{1},t_{2},t_{3},\ldots)\mapsto(\sum\limits_{h=1}^{\infty}t_{h}t^{h}+\varepsilon g(t),t_{1},t_{2},t_{3},\ldots).
\end{align*}
This is an arbitrarily small perturbation of the map $F_{\infty}$ and it can be shown that, because of non-analiticity, there are no solution curves to the equation $F_{\infty,\varepsilon}(\underline{t})=\underline{s}$. Hence a non-analytic modification of the map $F_{\infty}$ can annihilate the property of having solution curves; yet, as seen in \ref{sss122}, the origin $\underline{0}$ is still an $\infty$-transverse singularity. We wonder if we can produce a $C^{\infty}$ non-analytic perturbation of the map $F_{\infty}$ such that we do not have any $\infty$-transverse singularity.\\
\indent In general, it is an open problem to know if $C^{\omega}$ $\infty$-transverse singularity are stable under $C^{\infty}$ non-analytic perturbations or, even under a stricter request, $C^{\omega}$ perturbations of the map. One could raise a similar issue about the stability of $C^{\infty}$, non-analytic, $\infty$-transverse singularities, for example of the ones studied in Subsections \ref{sss122} and \ref{sss132}. \\ \vspace{8 pt}

\section{From Local to Pointwise Conditions}\label{s2}
\subsection{Pointwise Conditions for Low-Order Singularities}\label{ss21}
\quad  The classification of simple singularities introduced in Chapter 2 of \cite{B-D 1} for a smooth 0-Fredholm map $F$ is of a \textit{local nature}. Given a singularity $u_{o}$ for $F$ and a fixed fibering pair $(\varphi,\psi)$ near $u_{o}$ (cf. Definition 1.1.3 of \cite{B-D 1}), thanks to Definition 2.1.1 of \cite{B-D 1} we classify the singularity $u_{o}$ by means of the fibering functionals $J_{k},I_{k}$. By the very definition of these functionals (cf. Definition 1.1.6 of \cite{B-D 1}) this amounts to say that the iterated derivatives of both the map $F$ and the fibering maps $\varphi,\psi$ have to be computed on a suitable neighbourhood of $u_{o}$.\\
\indent  On the other hand, the conditions describing the type of low-order simple singularities found in the literature are of a \textit{pointwise nature}. This means that one computes the derivatives of  $F$ at $u_{o}$ and verifies that suitable conditions on these derivatives are satisfied. While this is certainly advantageous for the study of $k$-singularities with small values of $k$, it has the quite onerous drawback of requiring to check a fast growing number of conditions, for the derivatives of $F$ at $u_{o}$, as $k$ increases. However, this approach can be convenient when studying boundary value problems which exhibit low-order singularities, see e.g. Section \ref{ss25} below. Interestingly, the local approach seems to be more suitable when dealing with boundary value problems where $k$-singularities can occur for any integer $k$: it makes it possible to study all the $k$-singularities described in \cite{B-D 1} and to find more viable conditions in order to determine the kind of the studied singularity. This last assertion is clarified by the \textit{dual characterization} developed in \cite{B-D 3} and its use in the applications to boundary value problems, \cite{B-D 5}.\\
\indent The above considerations prompted us to show that, at least for low-order singularities, one can actually obtain a pointwise characterization of simple singularities from the previously introduced local conditions. We prove that the local knowledge of the iterated derivatives of $F$ and the fibering maps $\varphi,\psi$ on a neighbourhood of $u_{o}$ can be replaced by pointwise conditions only involving the derivatives of $F$ at $u_{o}$. \\
\indent Here, for historical reasons mainly, we restrict our study so as to provide a pointwise characterization of $k$-transverse singularities and $k$-singularities for $k=1,2,3,4$: this is accomplished in Theorems 2.1.1, ..., 2.1.4 below. We recall that, according to Thom's classification of the singularities for stable maps between four-dimensional manifolds, $k$-singularities for $k=1,2,3,4$ are usually named \textit{fold, cusp, swallow's tail} and \textit{butterfly} respectively. All the maps we consider in this chapter will be assumed to be of class $C^{\infty}$.
\begin{theorem}\label{Teo211}
 \textit{Let} $U,V$ \textit{be open subsets of the} $B$\textit{-spaces }$X,Y$\textit{ respectively and }$F:U\subseteq X \rightarrow V \subseteq Y$\textit{ a }$C^{\infty}$ 0-\textit{Fredholm map. Let }$u_{o}\in S_{1}(F)$\textit{ and let us consider }
\begin{equation}\label{211}
n_{o}\in N(F'(u_{o})) \setminus \{0\} \, , \; v_{o}\in X.
\end{equation}
\textit{Then we have the following}\\
\indent  \textbf{1-transversality condition:}\textit{ the point }$u_{o}$\textit{ is a }1-transverse singularity\textit{ if and only if there exist }$n_{o},v_{o}$\textit{ satisfying condition }(\ref{211})\textit{ with }$v_{o}$\textit{ such that}
\begin{equation}\label{212}
F''(u_{o})[n_{o}, v_{o}] \notin R(F'(u_{o})).
\end{equation}
\textit{In particular, we have the}\\
\indent \textbf{Fold condition: }\textit{the point }$u_{o}$\textit{ is a }fold (\textit{or} 1-singularity)\textit{ if and only if there exists }$n_{o}$\textit{ as in }(\ref{211})\textit{ such that}
\begin{equation}\label{213}
F''(u_{o})[n_{o}, n_{o}] \notin R(F'(u_{o})).
\end{equation}
\end{theorem}
\indent  It is straightforward to check that if the above fold condition is satisfied for a particular $n_{o} \in N(F'(u_{o})) \setminus \{0\}$ then it holds for all elements of $N(F'(u_{o})) \setminus \{0\}$.
\begin{theorem}\label{Teo212}
\textit{Let} $U,V$ \textit{be open subsets of the} $B$\textit{-spaces }$X,Y$\textit{ and }$F:U\subseteq X \rightarrow V \subseteq Y$\textit{ a }$C^{\infty}$ 0-\textit{Fredholm map. Let }$u_{o}\in S_{1}(F)$\textit{ be a } 1-transverse singularity \textit{which is not a} fold. \textit{Then there exist }$n_{o}, n_{1}, v_{o}, v_{1}\in X$\textit{ such that }
\begin{align}\label{214}
\begin{split}
&n_{o} \in N(F'(u_{o})) \setminus \{0\},\\
&F''(u_{o})[n_{o}, n_{o}] \in R(F'(u_{o})),\\
&F'(u_{o})n_{1}= -F''(u_{o})[n_{o},n_{o}],\\
&F''(u_{o})[n_{o}, v_{o}] \in R(F'(u_{o})),\\
&F'(u_{o})v_{1}= -F''(u_{o})[n_{o},v_{o}].\\
\end{split}
\end{align}
\textit{Then we have the following}\\
\indent \textbf{2-transversality condition:}\textit{ the point }$u_{o}$\textit{ is a }2-transverse singularity\textit{ if and only if there exist }$n_{o},n_{1},v_{o},v_{1}$\textit{  satisfying conditions }(\ref{214})\textit{ with }$v_{o}$\textit{  such that}
\begin{equation}\label{215}
F^{(3)}(u_{o})[n_{o},n_{o},v_{o}]+ 2F''(u_{o})[n_{o},v_{1}]+ F''(u_{o})[n_{1},v_{o}] \notin R(F'(u_{o})).
\end{equation}
\textit{In particular, we have the}\\
\indent \textbf{Cusp condition:}\textit{ the point }$u_{o}$\textit{ is a }cusp (\textit{or} 2-singularity)\textit{ if and only if there exist }$n_{o},n_{1}$\textit{ as in }(\ref{214})\textit{ such that}
\begin{equation}\label{216}
F^{(3)}(u_{o})[n_{o},n_{o},n_{o}]+ 3F''(u_{o})[n_{1},n_{o}] \notin R(F'(u_{o})).
\end{equation}
\end{theorem}
\begin{theorem} \label{Teo213}\textit{ Let} $U,V$ \textit{be open subsets of the} $B$\textit{-spaces} $X,Y$\textit{ and }$F:U\subseteq X \rightarrow V \subseteq Y$\textit{ a }$C^{\infty}$ 0-\textit{Fredholm map. Let }$u_{o}\in S_{1}(F)$\textit{ be a }2-transverse singularity \textit{which is not a} cusp. \textit{Then there exist }$n_{o}, n_{1}, n_{2}, v_{o}, v_{1}, v_{2}\in X$\textit{ such that }
\begin{align}\label{217}
\begin{split}
&n_{o} \in N(F'(u_{o})) \setminus \{0\},\\
&F''(u_{o})[n_{o}, n_{o}] \in R(F'(u_{o})),\\
&F'(u_{o})n_{1}= -F''(u_{o})[n_{o},n_{o}],\\
&F^{(3)}(u_{o})[n_{o},n_{o},n_{o}]+ 3F''(u_{o})[n_{1},n_{o}] \in R(F'(u_{o})),\\
&F'(u_{o})n_{2}=-\{F^{(3)}(u_{o})[n_{o},n_{o},n_{o}]+3F''(u_{o})[n_{1},n_{o}]\},\\
&F''(u_{o})[n_{o}, v_{o}] \in R(F'(u_{o})),\\
&F'(u_{o})v_{1}= -F''(u_{o})[n_{o},v_{o}],\\
&F^{(3)}(u_{o})[n_{o},n_{o},v_{o}]+ 2F''(u_{o})[n_{o},v_{1}]+F''(u_{o})[n_{1},v_{o}] \in R(F'(u_{o})),\\
&F'(u_{o})v_{2}=-\{F^{(3)}(u_{o})[n_{o},n_{o},v_{o}]+2F''(u_{o})[n_{o},v_{1}]+F''(u_{o})[n_{1},v_{o}]\}.
\end{split}
\end{align}
\textit{Then we have the following}\\
\indent \textbf{3-transversality condition:}\textit{ the point }$u_{o}$\textit{ is a }3-transverse singularity\textit{ if and only if there exist }$n_{o}, n_{1}, n_{2}, v_{o}, v_{1}, v_{2}$\textit{ satisfying conditions }(\ref{217})\textit{ with }$v_{o}$\textit{ such that }
\begin{equation}\label{218}
\begin{split}
&F^{(4)}(u_{o})[n_{o},n_{o},n_{o}, v_{o}]+ 3F^{(3)}(u_{o})[n_{o},n_{o},v_{1}]+3F^{(3)}(u_{o})[n_{1},n_{o},v_{o}]+\\
\hspace{30pt}+&3F''(u_{o})[n_{o},v_{2}]+3F''(u_{o})[n_{1},v_{1}]+F''(u_{o})[n_{2},v_{o}]\notin R(F'(u_{o})).
\end{split}
\end{equation} 
\textit{In particular, we have}\\
\indent \textbf{Swallow's tail condition:}\textit{ the point }$u_{o}$\textit{ is a }swallow's tail (\textit{or} 3-singularity)\textit{ if and only if there exist }$n_{o}, n_{1}, n_{2}$\textit{ as in }(\ref{217})\textit{ such that}
\begin{equation}\label{219}
\begin{split}
&F^{(4)}(u_{o})[n_{o},n_{o},n_{o}, n_{o}]+ 6F^{(3)}(u_{o})[n_{1},n_{o},n_{o}]+\\
\hspace{80pt} &+3F''(u_{o})[n_{1},n_{1}]+4 F''(u_{o})[n_{2},n_{o}]\notin R(F'(u_{o})).
\end{split}
\end{equation}
\end{theorem}
\begin{theorem}\label{Teo214} \textit{Let} $U,V$ \textit{be open subsets of the} $B$\textit{-spaces} $X,Y$\textit{ and }$F:U\subseteq X \rightarrow V \subseteq Y$\textit{ a }$C^{\infty}$ 0-\textit{Fredholm map. Let }$u_{o}\in S_{1}(F)$\textit{ be a } 3-transverse singularity \textit{which is not a} swallow's tail. \textit{Then there exist }$n_{o}, n_{1}, n_{2}, n_{3}\in X$\textit{ such that }
\begin{align}\label{2110}
\begin{split}
&n_{o} \in N(F'(u_{o})) \setminus \{0\},\\
&F''(u_{o})[n_{o}, n_{o}] \in R(F'(u_{o})),\\
&F'(u_{o})n_{1}= -F''(u_{o})[n_{o},n_{o}],\\
&F^{(3)}(u_{o})[n_{o},n_{o},n_{o}]+ 3F''(u_{o})[n_{1},n_{o}] \in R(F'(u_{o})),\\
&F'(u_{o})n_{2}=-\{F^{(3)}(u_{o})[n_{o},n_{o},n_{o}]+3F''(u_{o})[n_{1},n_{o}]\},\\
&F^{(4)}(u_{o})[n_{o},n_{o},n_{o}, n_{o}]+ 6F^{(3)}(u_{o})[n_{1},n_{o},n_{o}]+\\
& \qquad \qquad \quad+ 3F''(u_{o})[n_{1},n_{1}]+4 F''(u_{o})[n_{2},n_{o}]\in R(F'(u_{o})),\\
&F'(u_{o})n_{3}=-\{F^{(4)}(u_{o})[n_{o},n_{o},n_{o},n_{o}]+6F^{(3)}(u_{o})[n_{1},n_{o},n_{o}]+\\
& \qquad \qquad \quad+3F''(u_{o})[n_{1},n_{1}]+4F''(u_{o})[n_{2},n_{o}]\}.\\
\end{split}
\end{align}
\textit{Then, we have the following}\\
\indent \textbf{Butterfly condition:}\textit{ the point }$u_{o}$\textit{ is a }butterfly (\textit{or} 4-singularity)\textit{ if and only if there exist }$n_{o}, n_{1}, n_{2}, n_{3}$\textit{ as in }(\ref{2110})\textit{ such that}
\begin{align}\label{2111}
\begin{split}
&F^{(5)}(u_{o})[n_{o},n_{o},n_{o},n_{o},n_{o}]+10 F^{(4)}(u_{o})[n_{1},n_{o},n_{o},n_{o}]+10F^{(3)}(u_{o})[n_{2},n_{o},n_{o}]+\\
&+15F^{(3)}(u_{o})[n_{1},n_{1},n_{o}]+10F''(u_{o})[n_{2},n_{1}]+5F''(u_{o})[n_{3},n_{o}]\notin R(F'(u_{o})).\hspace{10pt}
\end{split}
\end{align}
\end{theorem}
\begin{remark}\label{Rem215}
The brief comment after Theorem \ref{Teo211} also applies to the conditions of $k$-singularity for $k=2,3,4$ expressed by formulas (\ref{216}), (\ref{219}), (\ref{2111}). For example, for $k=2$, if there exist elements $n_{o},n_{1}$ verifying formulas (\ref{214}) and satisfying the condition of 2-singularity (\ref{216}), then the latter is also satisfied by any vectors $m_{o},m_{1}$ verifying (\ref{214}). For the sake of brevity we omit to state the similar properties for $k=3,4$. 
The independence of the pointwise $k$-singularity conditions from $n_{o},\ldots,n_{k-1}$ is an easy consequence of Lemma \ref{Lem223}.\\
\indent We explicitly point out that, in the above conditions (\ref{213}), (\ref{216}), (\ref{219}) and  (\ref{2111}), only the iterated derivatives of the map $F$ are involved. Finally we anticipate that, in many examples associated with differential problems, the condition $w\notin R(F'(u_{o})), w\in Y$, can be easily checked by proving that a suitable integral is non-vanishing (cf. Section \ref{ss25}).
\end{remark}
\begin{remark}\label{Rem216}
Infinite-dimensional pointwise conditions for 1-transverse and fold singular points were first given in \cite{A-P}. Although no explicit connection with the two-dimensional Whitney fold was suggested, such conditions expressed quite natural geometrical assumptions. Indeed, it was shown there that the 1-transversality condition requires the singular set to be a hypersurface near the singularity and that, when the fold condition is also verified, the image of the singular set is, locally, a hypersurface in the target space. Since then, an extensive use of abstract pointwise conditions for 1-transverse singularities and 1-singularities has been adopted by a number of authors. In particular we can refer to \cite{B-C} and \cite{Ca}, where the fold nature of such singularities was first established. Later on, general pointwise conditions were also given to define 2-transverse and cusp singularities, cf. \cite{B-C-T} and \cite{L-M}. These definitions rely upon choosing suitable vectors in the domain space, but a proof of the independence from such choices is only found in \cite{B-C-T}. Less general definitions of cusp and butterfly singularities were also considered in \cite{C-D 1}, \cite{Ruf1} and \cite{Ruf2} respectively.\\
\indent 
We explicitly note that our pointwise conditions for $k$-transverse singularities with  $k=1,2$, that is folds and cusps, agree with those given in the above-mentioned papers. Hence one could assume the conditions shown in Theorems \ref{Teo211}, ..., \ref{Teo214} to be \textit{pointwise definitions} for the respective singularities.
\end{remark}

Theorems \ref{Teo211}, ..., \ref{Teo214} are proved in Section \ref{ss22}. The main analytical result involved in the demonstrations is a new characterization of singularities by the only means of the kernel fibering map $\varphi$: this is stated in Theorem \ref{Teo221}. Moreover, a more technical result is needed in order to show that the conditions given in Theorems \ref{Teo211}, ..., \ref{Teo214} for the vectors $n_{k},v_{k}$ are independent from the previously chosen elements $n_{i},v_{i},i=0,\ldots,k-1$. This is shown in Lemma \ref{Lem223} and it provides the fundamental pointwise relationship with the local conditions studied in Theorem \ref{Teo221}.\\
\indent The rest of this section focuses on the basic definitions which will then be employed throughout the remaining sections.
\begin{definition} \label{De217}
Let $F:U\subseteq X \rightarrow V \subseteq Y$ be a $C^{\infty}$ 0-Fredholm map between open subsets $U,V$ of the $B$-spaces $X,Y$, let $u_{o}\in S_{1}(F)$ and let $(\varphi,\psi)$ be a $C^{\infty}$ fibering pair for $F$ on $U_{o}\subseteq U,$ where $U_{o}$ is a suitable neighbourhood of $u_{o}$. We introduce the maps $\varphi_{k}(\varphi,\psi)\equiv \varphi_{k}:U_{o} \rightarrow X, k \geq 0$ as $\varphi_{o}:=\varphi$ and, for $k\geq 1$, as $\varphi_{k}(u):= \varphi^{\prime}_{k-1}(u)\varphi(u), u \in U_{o}$.
\end{definition}

We recall that the existence of a $C^{\infty}$ fibering pair for $F$ near a singular point is guaranteed by Theorem 1.4.6 of \cite{B-D 1}. The maps $\varphi_{k}$ are also of class $C^{\infty}$ and they are nothing but the iterated derivatives of $\varphi$ along $\varphi$ itself (see Definition 1.1.7 in \cite{B-D 1}). By means of these derivatives we introduce the following families of maps:
\begin{definition} \label{De218}
Let $F:U\subseteq X \rightarrow V \subseteq Y$ be a $C^{\infty}$ 0-Fredholm map
between open subsets $U,V$ of the $B$-spaces $X,Y$, let $u_{o}\in S_{1}(F)$ and let $(\varphi,\psi)$ be a $C^{\infty}$ fibering pair for $F$ on $U_{o}\subseteq U$, a suitable neighbourhood of $u_{o}$. We define
\begin{align*}
\tau_{1}(\varphi,\psi)(u) &\equiv \tau_{1}(u):X \rightarrow Y \, , \, \forall\, u\in U_{o},\quad &\text{ as } \quad \tau_{1}(u)v&:=F''(u)[\varphi(u),v], v \in X ,\\
\sigma_{1}(\varphi,\psi)&\equiv \sigma_{1}:U_{o} \rightarrow Y  &\text{ as } \quad \,\, \sigma_{1}(u)&:=\tau_{1}(u)\varphi(u), u\in U_{o} \;\; , 
\end{align*}
and, inductively for $k\geq 2$, 
\begin{align*}
\tau_{k}(\varphi,\psi)(u) \equiv \tau_{k}(u):X \rightarrow Y \, , \, \forall\, u\in U_{o},\quad &\text{ as } \quad \tau_{k}(u)v:=\sigma^{\prime}_{k-1}(u)v \,+ \\
& \qquad \qquad \quad \;\;\; +F''(u)[\varphi_{k-1}(u),v]\, , \, v \in X,\\
\sigma_{k}(\varphi,\psi)\equiv \sigma_{k}:U_{o} \rightarrow Y \hspace{80pt} &\text{ as } \quad \;\, \sigma_{k}(u):=\tau_{k}(u)\varphi(u), u\in U_{o} .\;\; \qquad \qquad  \qquad \qquad
\end{align*}
\end{definition}
\indent We explicitly note that $\tau_{k}(u)\in L(X,Y), \forall \,k\geq 1, \forall \, u\in U_{o}$, hence we can consider the map $\tau_{k}:U_{o}\rightarrow L(X,Y)$. Then one can verify that $\tau_{k}\in C^{\infty}(U_{o},L(X,Y))$ and $\sigma_{k}\in C^{\infty}(U_{o},Y)$. All these properties can be proved inductively.\\
\par
In order to obtain a pointwise classification for $k$-transverse singularities and $k$-singularities for $k=1,2,3,4$ it is useful to know the explicit form of the maps $\tau_{k},\sigma_{k}$ for $k=1,2,3,4$. Therefore we state
\begin{proposition} \label{Pro219} \textit{Let} $F:U\subseteq X \rightarrow V \subseteq Y$ \textit{be a }$C^{\infty}$ 0-\textit{Fredholm map between open subsets} $U,V$\textit{ of the }$B$\textit{-spaces }$X,Y$;\textit{ assume that }$u_{o}\in S_{1}(F)$ \textit{ and that }$(\varphi,\psi)$\textit{ is a }$C^{\infty}$ \textit{ fibering pair for }$F$\textit{ on }$U_{o}\subseteq U$,\textit{ with } $U_{o}$\textit{ a suitable neighbourhood of }$u_{o}$.\textit{ Then, for }$k=1,2,3,4$\textit{ and for }$u \in U_{o}, v \in X$, \textit{ the maps }$\tau_{k}, \sigma_{k}$\textit{ have the form} \medskip \\ 
$\tau_{1}(u)v=F''(u)[\varphi(u),v]\,; $ \vspace{10pt} \\
$\sigma_{1}(u)=F''(u)[\varphi(u),\varphi(u)]\,;$ \vspace{10pt} \\
$\tau_{2}(u)v=F^{(3)}(u)[\varphi(u),\varphi(u),v]+2F''(u)[\varphi'(u)v,\varphi(u)]+F''(u)[\varphi_{1}(u),v]\,;$\vspace{10pt}\\
$\sigma_{2}(u)=F^{(3)}(u)[\varphi(u),\varphi(u),\varphi(u)]+3F''(u)[\varphi_{1}(u),\varphi(u)]\,; \vspace{3pt}$ 
\begin{align*}
\tau_{3}(u)v&=F^{(4)}(u)[\varphi(u),\varphi(u),\varphi(u),v]+3F^{(3)}(u)[\varphi(u),\varphi(u),\varphi'(u)v]+\\
&+3F^{(3)}(u)[\varphi_{1}(u),\varphi(u),v]+3F''(u)[\varphi(u),\varphi^{\prime}_{1}(u)v]+ \qquad\qquad\qquad\qquad\qquad\qquad\qquad\qquad\qquad\qquad\\
&+3F''(u)[\varphi_{1}(u),\varphi^{\prime}(u)v]+F''(u)[\varphi_{2}(u),v]\, ;\qquad\qquad\qquad\qquad\qquad
\end{align*}
\begin{align*}
\sigma_{3}(u)&=F^{(4)}(u)[\varphi(u),\varphi(u),\varphi(u),\varphi(u)]+6F^{(3)}(u)[\varphi_{1}(u),\varphi(u),\varphi(u)]+ \qquad\qquad\qquad\qquad\qquad\qquad\qquad\qquad\qquad\qquad\\
&+4F''(u)[\varphi_{2}(u),\varphi(u)]+3F''(u)[\varphi_{1}(u),\varphi_{1}(u)]\, ;
\end{align*}
\begin{align*}
\tau_{4}(u)v&=F^{(5)}(u)[\varphi(u),\varphi(u),\varphi(u),\varphi(u),v]+4F^{(4)}(u)[\varphi(u),\varphi(u),\varphi(u),\varphi'(u)v] +\\
&+6F^{(4)}(u)[\varphi_{1}(u),\varphi(u),\varphi(u),v]+6F^{(3)}(u)[\varphi(u),\varphi(u),\varphi^{\prime}_{1}(u)v] +\\
&+12F^{(3)}(u)[\varphi_{1}(u),\varphi(u),\varphi^{\prime}(u)v]+4F^{(3)}(u)[\varphi_{2}(u),\varphi(u),v]+\\
&+3F^{(3)}(u)[\varphi_{1}(u),\varphi_{1}(u),v] +4F''(u)[\varphi(u),\varphi^{\prime}_{2}(u)v]+4F''(u)[\varphi_{2}(u),\varphi^{\prime}(u)v]+\\
&+6F''(u)[\varphi_{1}(u),\varphi^{\prime}_{1}(u)v]+F''(u)[\varphi_{3}(u),v]\, ;
\end{align*}
\begin{align*}
\sigma_{4}(u)&=F^{(5)}(u)[\varphi(u),\varphi(u),\varphi(u),\varphi(u),\varphi(u)]+10F^{(4)}(u)[\varphi_{1}(u),\varphi(u),\varphi(u),\varphi(u)] +\\
&+10F^{(3)}(u)[\varphi_{2}(u),\varphi(u),\varphi(u)]+15F^{(3)}(u)[\varphi_{1}(u),\varphi_{1}(u),\varphi(u)]+\\
&+5F''(u)[\varphi_{3}(u),\varphi(u)]+10F''(u)[\varphi_{2}(u),\varphi_{1}(u)].
\end{align*}
\end{proposition}
\textbf{Proof.} By Definition \ref{De218}, $\tau_{1}(u)v=F''(u)[\varphi(u),v], u \in U_{o}, v \in X$ and thus 
\begin{center}
$\sigma_{1}(u)=\tau_{1}(u)\varphi(u)=F''(u)[\varphi(u),\varphi(u)].$
\end{center}
A direct computation gives 
$\sigma^{\prime}_{1}(u)v=F^{(3)}(u)[\varphi(u),\varphi(u),v]+2F''(u)[\varphi'(u)v,\varphi(u)],\\
 u \in U_{o}, v \in X$. Hence, again by \ref{De218}, we get that 
\begin{align*}
\begin{split}
\tau_{2}(u)v &=\sigma^{\prime}_{1}(u)v + F''(u)[\varphi_{1}(u),v] =\\
&=F^{(3)}(u)[\varphi(u),\varphi(u),v]+ 2F''(u)[\varphi'(u)v,\varphi(u)]+ F''(u)[\varphi_{1}(u),v].
\end{split}
\end{align*}
Since $\varphi'(u)\varphi(u)=\varphi_{1}(u)$, cf. Definition \ref{De217}, we conclude that
\begin{center}
$\sigma_{2}(u)=\tau_{2}(u)\varphi(u)=F^{(3)}(u)[\varphi(u),\varphi(u),\varphi(u)] +3F''(u)[\varphi_{1}(u),\varphi(u)].$
\end{center}
The formulas for $\tau_{3}(u),\sigma_{3}(u),\tau_{4}(u), \sigma_{4}(u)$ can be deduced in the same way. $\bracevert$ \\
\par In order to avoid cumbersome notations, we just assign suitable labels to the maps involved in Theorems \ref{Teo211}, ..., \ref{Teo214}.
\begin{definition} \label{De2110}
Let $F:U\subseteq X \rightarrow V \subseteq Y$ be a $C^{\infty}$ map
between open subsets $U,V$ of the $B$-spaces $X,Y$ and let $u_{o}\in U$ be a given point. For any $n_{o},n_{1},n_{2},n_{3},v_{o},v_{1},v_{2}, v_{3} \in X$ we define the following maps:
\begin{align*}
T_{1}(u_{o}) \equiv\;&\;T_{1}:X^{2} \rightarrow Y \, \, \,\quad \text{ as } \quad
T_{1}[n_{o},v_{o}]:=F''(u_{o})[n_{o},v_{o}]\, ; \\
\\
\varSigma_{1}(u_{o})\equiv & \; \varSigma_{1}:X \rightarrow Y \, \, \, \;\quad \text{ as } \quad  \varSigma_{1}[n_{o}]:=F''(u_{o})[n_{o},n_{o}] \; ; \quad  \\
\\
T_{2}(u_{o}) \equiv & \;T_{2}:X^{4} \rightarrow Y \, \, \,\quad \text{ as } \quad \\
&\;T_{2}[n_{o},n_{1},v_{o},v_{1}]:=F^{(3)}(u_{o})[n_{o},n_{o,}v_{o}]+2F''(u_{o})[n_{o},v_{1}]+F''(u_{o})[n_{1},v_{o}]\, ; \\
\\
\varSigma_{2}(u_{o})\equiv & \;\varSigma_{2}:X^{2} \rightarrow Y \, \, \quad \text{ as } \quad  \varSigma_{2}[n_{o},n_{1}]:=F^{(3)}(u_{o})[n_{o},n_{o},n_{o}]+3F''(u_{o})[n_{1},n_{o}] \, ; \\
\\
T_{3}(u_{o}) \equiv & \;T_{3}:X^{6} \rightarrow Y \, \, \,\quad \text{ as } \quad \\
&\;T_{3}[n_{o},n_{1},n_{2},v_{o},v_{1},v_{2}]:=F^{(4)}(u_{o})[n_{o},n_{o},n_{o},v_{o}]+3F^{(3)}(u_{o})[n_{o},n_{o},v_{1}]+\\
&\;+3F^{(3)}(u_{o})[n_{1},n_{o},v_{o}]+3F''(u_{o})[n_{o},v_{2}]+3F''(u_{o})[n_{1},v_{1}]+F''(u_{o})[n_{2},v_{o}]\, ;\\
\\
\varSigma_{3}(u_{o})\equiv & \;\varSigma_{3}:X^{3} \rightarrow Y \, \, \,\quad \text{ as } \quad \\
&\;\varSigma_{3}[n_{o},n_{1},n_{2}]:=F^{(4)}(u_{o})[n_{o},n_{o},n_{o},n_{o}]+6F^{(3)}(u_{o})[n_{1},n_{o},n_{o}]+\\
&\;+3F''(u_{o})[n_{1},n_{1}]+4F''(u_{o})[n_{2},n_{o}] \, ; \\
\\
T_{4}(u_{o}) \equiv &\;T_{4}:X^{8} \rightarrow Y \, \, \,\quad \text{ as } \quad \\
&\;T_{4}[n_{o},n_{1},n_{2},n_{3},v_{o},v_{1},v_{2},v_{3}]:=F^{(5)}(u_{o})[n_{o},n_{o},n_{o},n_{o},v_{o}]+\\
&\;+4F^{(4)}(u_{o})[n_{o},n_{o},n_{o},v_{1}]+6F^{(4)}(u_{o})[n_{1},n_{o},n_{o},v_{o}]+6F^{(3)}(u_{o})[n_{o},n_{o},v_{2}]+\\
&\;+12F^{(3)}(u_{o})[n_{1},n_{o},v_{1}]+4F^{(3)}(u_{o})[n_{2},n_{o},v_{o}]+3F^{(3)}(u_{o})[n_{1},n_{1},v_{o}]+\\
&\;+4F''(u_{o})[n_{o},v_{3}]+4F''(u_{o})[n_{2},v_{1}]+6F''(u_{o})[n_{1},v_{2}]+F''(u_{o})[n_{3},v_{o}]\, ;\\
\\
\varSigma_{4}(u_{o})\equiv &\;\varSigma_{4}:X^{4} \rightarrow Y \, \, \,\quad \text{ as } \quad \\
&\;\varSigma_{4}[n_{o},n_{1},n_{2},n_{3}]:=F^{(5)}(u_{o})[n_{o},n_{o},n_{o},n_{o},n_{o}]+10F^{(4)}(u_{o})[n_{1},n_{o},n_{o},n_{o}]+\\
&\;+10F^{(3)}(u_{o})[n_{2},n_{o},n_{o}]+15F^{(3)}(u_{o})[n_{1},n_{1},n_{o}]+10F''(u_{o})[n_{2},n_{1}]+\\
&\;+5F''(u_{o})[n_{3},n_{o}] \,.
\end{align*}
As a matter of notation, when there is no risk of confusion about the point $u_{o}$ we will always write $T_{k} ,\varSigma_{k}$ instead of $T_{k}(u_{o}) ,\varSigma_{k}(u_{o})$.
\end{definition}
\begin{remark} \label{Rem2111}
We note that the pointwise conditions of Theorems \ref{Teo211}, ..., \ref{Teo214} can be restated in terms of the maps $T_{k}$ and $\varSigma_{k}, k=1,2,3,4$. In a way, we could consider the maps $T_{k} ,\varSigma_{k}$ as a sort of ``polynomials'' in the ``indeterminates'' $n_{o},n_{1},n_{2},n_{3},v_{o},v_{1},v_{2},v_{3}$. In this view we then remark that, by evaluating the ``polynomials'' $T_{k}$ and $\varSigma_{k}$ at $\varphi_{h}(u_{o}), \varphi'_{h}(u_{o})v$, we get $\tau_{k}(u_{o})v,\sigma_{k}(u_{o})$. More precisely we obtain that, $\forall v\in X$, 
\begin{align*}
\begin{split}
\tau_{1}(u_{o})v &= T_{1}[\varphi(u_{o}),v] \, ;\\
\sigma_{1}(u_{o}) &= \varSigma_{1}[\varphi(u_{o})] \, ;\\
\tau_{2}(u_{o})v &= T_{2}[\varphi(u_{o}),\varphi_{1}(u_{o}),v,\varphi'(u_{o})v] \, ;\\
\sigma_{2}(u_{o}) &= \varSigma_{2}[\varphi(u_{o}),\varphi_{1}(u_{o})] \, ;\\
\tau_{3}(u_{o})v &= T_{3}[\varphi(u_{o}),\varphi_{1}(u_{o}),\varphi_{2}(u_{o}),v,\varphi'(u_{o})v,\varphi_{1}'(u_{o})v] \, ;\\
\sigma_{3}(u_{o}) &= \varSigma_{3}[\varphi(u_{o}),\varphi_{1}(u_{o}),\varphi_{2}(u_{o})] \, ;\\
\qquad\tau_{4}(u_{o})v &=T_{4}[\varphi(u_{o}),\varphi_{1}(u_{o}),\varphi_{2}(u_{o}),\varphi_{3}(u_{o}),v,\varphi'(u_{o})v,\varphi_{1}'(u_{o})v,\varphi_{2}'(u_{o})v] \, ;\\
\sigma_{4}(u_{o}) &= \varSigma_{4}[\varphi(u_{o}),\varphi_{1}(u_{o}),\varphi_{2}(u_{o}),\varphi_{3}(u_{o})] \,.
\end{split}
\end{align*}
From these identities we see that the equalities $\sigma_{k}(u_{o})=\tau_{k}(u_{o})\varphi(u_{o})$  could be reformulated in terms of the maps $T_{k}$ and $\varSigma_{k}$ as 
\begin{align}\label{2112}
\begin{split}
T_{1}[n_{o},n_{o}] &= \varSigma_{1}[n_{o}] \, ,\\
T_{2}[n_{o},n_{1},n_{o},n_{1}] &= \varSigma_{2}[n_{o},n_{1}] \, ,\\
T_{3}[n_{o},n_{1},n_{2},n_{o},n_{1},n_{2}] &= \varSigma_{3}[n_{o},n_{1},n_{2}] \, ,\\
T_{4}[n_{o},n_{1},n_{2},n_{3},n_{o},n_{1},n_{2},n_{3}] &= \varSigma_{4}[n_{o},n_{1},n_{2},n_{3}]\, . \qquad\quad \\
\end{split}
\end{align}
\end{remark}
\vspace{2pt}
\subsection{Proof of the Pointwise Conditions}\label{ss22}
\quad The maps $\tau_{k}$ and $\sigma_{k}$ allow us to give a useful characterization of the singularities and to establish some important equations on the singular strata. In this characterization only the iterated derivatives $\varphi_{k}$ of a given kernel fibering map are involved, as shown in Theorem \ref{Teo221} below. We note that for the proofs of Theorems \ref{Teo211}, ..., \ref{Teo214} we only need Theorem \ref{Teo221} for $k = 1, 2, 3, 4$. However, we state the result for any integer $k$ as both the statement and the proof gain in elegance and brevity. For the proof of Theorem \ref{Teo221} we refer to Section \ref{ss23}. We recall that the symbol $\mathcal{P}(F,u_{o})$ denotes the set of all smooth (germs of) fibering pairs for $F$ at $u_{o}$, cf. Subsection 1.4.1 of \cite{B-D 1}.
\begin{theorem}\label{Teo221}
\textit{Let} $U,V$ \textit{be open subsets in the} $B$\textit{-spaces} $X,Y$\textit{ and let }$F:U\subseteq X \rightarrow V \subseteq Y$\textit{ be a }$C^{\infty}$ 0-\textit{Fredholm map. Given }$u_{o}\in S_{1}(F)$\textit{ and }$(\varphi,\psi) \in \mathcal{P}(F,u_{o})$ \textit{ we assume that, for some integer }$k\geq 0, u_{o}$\textit{ is }$k$-transverse \textit{and not a }$k$-singularity \textit{(for $k=0$ this is an empty condition). Then the following equivalences hold:}
\begin{equation*}
u_{o}\textit{ is a }(k+1)\text{-transverse singularity}  \Leftrightarrow (LT_{k+1}):
\begin{cases}
\exists \; v\in X:\tau_{1}(u_{o})v \in R(F'(u_{o})),\\
\quad \quad , \ldots,\tau_{k}(u_{o})v \in R(F'(u_{o})),\\
\qquad \quad \; \tau_{k+1}(u_{o})v \notin R(F'(u_{o})).
\end{cases}
\end{equation*}
\textit{In particular,}
\begin{equation*}
u_{o}\textit{ is a }(k+1)\text{-singularity}\; \Leftrightarrow \;(LS_{k+1}): \quad \sigma_{k+1}(u_{o}) \notin R(F'(u_{o})) .
\end{equation*}
\textit{Moreover, if }$k\geq 1$\textit{ we have that, }$\forall \, h=1,\ldots,k$\textit{, the vectors }$\varphi_{h}(u_{o}) \in X$\textit{ satisfy the identities} \medskip\\
$(i_{h}) \qquad\qquad \qquad\qquad\qquad F'(u_{o})\varphi_{h}(u_{o}) = -\sigma_{h}(u_{o})$.
\end{theorem}
\begin{remark}\label{Rem222}
The above theorem is an intermediate step in our transition from local to pointwise conditions characterizing simple singularities. In fact, to state this theorem we only need to know the kernel map $\varphi$ locally, near the singularity $u_{o}$, as only the iterated derivatives of $\varphi$ are used to evaluate the maps $\tau_{h}$ and $\sigma_{h}$ at $u_{o}$. Note, however, that the cokernel $\psi(u_{o})$ is also used in the above statement (even if in an implicit way): this is because one has to decide whether or not $\tau_{h}(u_{o})v$ belongs to $R(F'(u_{o}))$ for $h=1, \ldots,k+1$. Indeed, by Definition 1.1.3 in \cite{B-D 1} we have that $\psi(u_{o}) \in Y^{\ast}\setminus \{0\}$ is such that $\psi(u_{o})w=0$\textit{ iff }$w\in R(F'(u_{o}))$. Finally, in the following we will often use the fact that there always exist vectors $v\in X$ such that $\tau_{h}(u_{o})v \in R(F'(u_{o})), h=1,\ldots,k$. In fact, $\tau_{h}(u_{o})v \in R(F'(u_{o}))$\textit{ iff }$\psi(u_{o})\tau_{h}(u_{o})v = 0$. By recalling that $\psi(u_{o})\tau_{h}(u_{o})v \in X^{\ast}$ it follows that the subset of such vectors $v$ is a subspace of $X$ which has at most codimension $k$.
\end{remark}
In order to apply the previous theorem to an arbitrary $(k+1)$-uple of vectors $(n_{o},\ldots,n_{k})$ satisfying formulas (\ref{211}), (\ref{214}), (\ref{217}) and (\ref{2110}) of Theorems \ref{Teo211}, ..., \ref{Teo214}, and thus prove the pointwise conditions for $k$-singularities given by formulas (\ref{213}), (\ref{216}), (\ref{219}) and (\ref{2111}), we need to know how the $(k+1)$-uples $(n_{o},\ldots,n_{k})$ and $(\varphi(u_{o}),\ldots,\varphi_{k}(u_{o}))$ behave under the action of the map $\varSigma_{k+1}$ for a given $(\varphi,\psi)\in \mathcal{P}(F,u_{o})$. Although we cannot in general expect that $\varSigma_{k+1}[n_{o},\ldots,n_{k}]=\varSigma_{k+1}[\varphi(u_{o}),$\\
$\ldots,\varphi_{k}(u_{o})]$, it is nonetheless possible to find suitable relations among the vectors $n_{h}$ and $\varphi_{h}(u_{o})$ that guarantee a weaker result which suffices to prove the pointwise conditions. These relations, and analogous ones for the map $T_{k+1},k=0,1,2,3$ are given in the following 
\begin{lemma}\label{Lem223}(Combinatorial Lemma). \textit{Let} $U,V$ \textit{be open subsets in the} $B$-spaces $X,Y$\textit{ and }$F:U\subseteq X \rightarrow V \subseteq Y$\textit{ a }$C^{\infty}$ 0-\textit{Fredholm map. Let }$u_{o}\in S_{1}(F)$\textit{ and }$(\varphi,\psi) \in \mathcal{P}(F,u_{o})$.\textit{ For a given integer }$k=0,\ldots,3$\textit{ let us assume that }$u_{o}$\textit{ is a }$k$-transverse singularity\textit{ which is not a }$k$-singularity\textit{ (empty condition for }$k=0$\textit{) and that there exist vectors }$n_{o},\ldots,n_{k}, v_{o},\ldots,v_{k}$\textit{ satisfying conditions  }$(\ref{211})$\textit{ for k=0,  } $(\ref{214})$\textit{ for k=1,  }$ (\ref{217})$\textit{ for k=2,  }$(\ref{2110})$\textit{ for k=3}.\textit{ Then, for any }$h=0,\ldots,k$,\textit{ the following equivalences hold:} 
\begin{align*}
(s_{h}) \quad &\varSigma_{h+1}[\varphi(u_{o}),\ldots,\varphi_{h}(u_{o})]\in R(F'(u_{o})) \Leftrightarrow \varSigma_{h+1}[n_{o},\ldots,n_{h}]\in R(F'(u_{o}));\\
(t_{h}) \quad &T_{h+1}[\varphi(u_{o}),\ldots,\varphi_{h}(u_{o}),v_{o},\ldots,\varphi'_{h-1}(u_{o})v_{o}]\in R(F'(u_{o}))\\
& \hspace{175pt}\Leftrightarrow T_{h+1}[n_{o},\ldots,n_{h},v_{o},\ldots,v_{h}]\in R(F'(u_{o})).
\end{align*}
\end{lemma}
\par By means of Remark \ref{Rem2111}, Theorem \ref{Teo221} and Lemma \ref{Lem223} (whose proof is postponed to Section \ref{ss24}) we are now able to prove Theorems \ref{Teo211}, ..., \ref{Teo214}. \\
\par
\textbf{Proof of Theorem \ref{Teo211}} Let $(\varphi,\psi) \in \mathcal{P}(F,u_{o})$ be a given f-pair. Since, by construction, $\varphi(u_{o})\in N(F'(u_{o}))$ and $\varphi(u_{o})\neq 0$, we have that $N(F'(u_{o}))\setminus\{0\} \neq\emptyset$. Now we prove in a separate way the fold and the 1-transversality conditions. \medskip\\
\indent\textit{Fold condition.} By Theorem \ref{Teo221} we know that $u_{o}$ is a fold \textit{iff} condition $(LS_{1})$ holds, i.e. \textit{iff} $\sigma_{1}(u_{o})\notin R(F'(u_{o}))$. Moreover, from Remark \ref{Rem2111} we have that $\sigma_{1}(u_{o})=\varSigma_{1}[\varphi(u_{o})]$. By taking $n_{o} \in N(F'(u_{o}))\setminus\{0\}$ we get, thanks to $(s_{o})$ of Lemma \ref{Lem223}, that $\varSigma_{1}[\varphi(u_{o})] \in R(F'(u_{o}))\textit{ iff } \varSigma_{1}[n_{o}]\in R(F'(u_{o}))$ and hence we can conclude that $u_{o}$ is a fold \textit{iff}  $\varSigma_{1}[n_{o}]\notin R(F'(u_{o}))$.\medskip\\
\indent \textit{1-transversality condition - Proof of}$\Leftarrow$. We choose $(n_{o},v_{o})\in (N(F'(u_{o}))\setminus \{0\})\times X$ such that $T_{1}[n_{o},v_{o}]\notin R(F'(u_{o}))$. By the equivalence $(t_{o})$ of Lemma \ref{Lem223} we have that $T_{1}[\varphi(u_{o}),v_{o}]\notin R(F'(u_{o}))$. On the other hand, from Remark \ref{Rem2111} we get $T_{1}[\varphi(u_{o}),v_{o}]=\tau_{1}(u_{o})v_{o}$. Hence $\tau_{1}(u_{o})v_{o} \notin R(F'(u_{o}))$. Thus $(LT_{1})$ of Theorem \ref{Teo221} is satisfied with $v=v_{o}$ and we conclude that $u_{o}$ is 1-transverse.\medskip\\
\indent \textit{1-transversality condition - Proof of}$\Rightarrow$. Let $u_{o}$ be a 1-transverse singularity. Then, from $(LT_{1})$ of Theorem \ref{Teo221}, there exists $v_{o}\in X$ such that $\tau_{1}(u_{o})v_{o} \notin R(F'(u_{o}))$. Hence from Remark \ref{Rem2111} one has that $T_{1}[\varphi(u_{o}),v_{o}]\notin R(F'(u_{o}))$ and, for $n_{o}\in N(F'(u_{o}))\setminus \{0\}$, thanks to $(t_{o})$ of Lemma \ref{Lem223}  we obtain that $T_{1}[n_{o},v_{o})]\notin R(F'(u_{o})). \bracevert$\\
\par
\textbf{Proof of Theorem \ref{Teo212}} Let $(\varphi,\psi) \in \mathcal{P}(F,u_{o})$ be a given f-pair. By hypothesis $u_{o}$ is a 1-transverse singularity but not a fold, hence the equality $(i_{1})$ of Theorem \ref{Teo221} holds, i.e. $F'(u_{o})\varphi_{1}(u_{o})= -\sigma_{1}(u_{o})$. On the other hand $\sigma_{1}(u_{o})=\varSigma_{1}[\varphi(u_{o})]$. Since $\varphi(u_{o})\in N(F'(u_{o}))\setminus \{0\}$ and $F'(u_{o})\varphi_{1}(u_{o})= -\varSigma_{1}[\varphi(u_{o})]$ then, by choosing $n_{o}:=\varphi(u_{o})$ and $n_{1}:=\varphi_{1}(u_{o})$, the first three conditions of (\ref{214}) are satisfied. This also implies that there exist vectors $n_{o},n_{1},v_{o},v_{1}$ satisfying all the conditions of (\ref{214}): it suffices to choose $v_{o} :=n_{o},v_{1}:=n_{1}$ with $n_{o},n_{1}$ as above. Now we separately prove the cusp and the 2-transversality conditions. \medskip\\
\indent \textit{ Cusp condition.} Since $u_{o}$ is a 1-transverse singularity but not a fold, Theorem \ref{Teo221} states that $u_{o}$ is a cusp \textit{iff} condition $(LS_{2})$ holds, i.e. \textit{iff} $\sigma_{2}(u_{o})\notin R(F'(u_{o}))$. We recall that, by Remark \ref{Rem2111}, $\sigma_{2}(u_{o})=\varSigma_{2}[\varphi(u_{o}),\varphi_{1}(u_{o})]$. As shown above, there exist $n_{o},n_{1},v_{o},v_{1}$ satisfying conditions \ref{214}, hence the assumptions of Lemma \ref{Lem223} are verified. From the equivalence $(s_{1})$ of Lemma \ref{Lem223} we have that $\varSigma_{2}[\varphi(u_{o}),\varphi_{1}(u_{o})] \notin R(F'(u_{o}))$ \textit{iff }$\varSigma_{2}[n_{o},n_{1}]\notin R(F'(u_{o}))$. Therefore $u_{o}$ is a cusp \textit{iff} $\varSigma_{2}[n_{o},n_{1}]\notin R(F'(u_{o}))$ and this is equivalent to saying that (\ref{216}) is satisfied.\medskip\\
\indent \textit{2-transversality condition - Proof of}$\Leftarrow$. Let $n_{o},n_{1},v_{o},v_{1}$  satisfy conditions (\ref{214}) and such that (\ref{215}) is verified, i.e. $T_{2}[n_{o},n_{1},v_{o},v_{1}]\notin R(F'(u_{o}))$.
From the fourth condition in (\ref{214}) we have that $T_{1}[n_{o},v_{o}]\in R(F'(u_{o}))$. By $(t_{1})$ and $(t_{o})$ of Lemma \ref{Lem223} we obtain that $T_{2}[\varphi(u_{o}),\varphi_{1}(u_{o}),v_{o},\varphi'(u_{o})v_{o}]\notin R(F'(u_{o}))$\ and $T_{1}[\varphi (u_{o}),v_{o}]\in R(F'(u_{o}))$. Thanks to Remark \ref{Rem2111} we know that $T_{2}[\varphi(u_{o}),\varphi_{1}(u_{o}),v_{o},\varphi'(u_{o})v_{o}]=\tau_{2}(u_{o})v_{o} \notin  R(F'(u_{o}))$ and $T_{1}[\varphi (u_{o}),v_{o}]= \tau_{1}(u_{o})v_{o} \in R(F'(u_{o}))$. By taking $v=v_{o}$ we have that $(LT_{2})$ of Theorem \ref{Teo221} is satisfied, hence we can conclude that $u_{o}$ is a 2-transverse singularity.\medskip\\
\indent \textit{2-transversality condition - Proof of}$\Rightarrow$. Let $u_{o}$ be a 2-transverse singularity. Then $(LT_{2})$ of Theorem \ref{Teo221} holds; let us call $v_{o}$ the element of $X$ such that $\tau_{1}(u_{o})v_{o} \in R(F'(u_{o}))$ and $\tau_{2}(u_{o})v_{o} \notin R(F'(u_{o}))$. Then, by Remark \ref{Rem2111} we have that $T_{1}[\varphi (u_{o}),v_{o}] \in R(F'(u_{o}))$ and $T_{2}[\varphi(u_{o}),\varphi_{1}(u_{o}),v_{o},\varphi'(u_{o})v_{o}] \notin R(F'(u_{o}))$. Now let $n_{o},n_{1}$ be such that the first three conditions of (\ref{214}) are satisfied. Since $u_{o} \in S_{1}(F)$ we can use Lemma \ref{Lem223} for $k=0$. Equivalence $(t_{o})$ implies that $T_{1}[n_{o},v_{o}]\in R(F'(u_{o}))$.  Then there exists $v_{1} \in X$ such that $F'(u_{o})v_{1}= -T_{1}[n_{o},v_{o}]$, hence $n_{o},n_{1},v_{o},v_{1}$ satisfy all the conditions \ref{214} and, from $(t_{1})$ of Lemma \ref{Lem223}, we obtain that $T_{2}[n_{o},n_{1},v_{o},v_{1}]\notin R(F'(u_{o}))$ which is equivalent to condition \ref{215}. $\bracevert$\\
\par
\textbf{Proof of Theorem \ref{Teo213}}. By hypothesis we know that $u_{o}$ is a 2-transverse singularity but not a cusp and so the equalities $(i_{1}),(i_{2})$ of Theorem \ref{Teo221} hold for a given f-pair $(\varphi,\psi) \in \mathcal{P}(F,u_{o}).$ i.e. $F'(u_{o})\varphi_{1}(u_{o})= -\sigma_{1}(u_{o})$ and $F'(u_{o})\varphi_{2}(u_{o})= -\sigma_{2}(u_{o})$. Since, by Remark \ref{Rem2111}, $\sigma_{1}(u_{o})=\varSigma_{1}[\varphi (u_{o})]$ and $\sigma_{2}(u_{o})=\varSigma_{2}[\varphi (u_{o}),\varphi_{1} (u_{o})]$ we easily get that $n_{o}:=\varphi (u_{o}),n_{1}:=\varphi_{1}(u_{o}),n_{2}:=\varphi_{ 2}(u_{o})$ satisfy the first five conditions of (\ref{217}). This implies that there exist elements $n_{o},n_{1},n_{2},v_{o},v_{1},v_{2}$ which satisfy all the conditions in (\ref{217}). It is then sufficient to take $v_{o}:=n_{o},v_{1}:=n_{1},v_{2}:=n_{2}$ with $n_{o},n_{1},n_{2}$ as above and apply formula (\ref{2112}) of Remark \ref{2111}. Once again, we treat the proof of the swallow's tail and 3-transversality conditions separately. \medskip \\
\indent\textit{ Swallow's tail condition.} We know that $u_{o}$ is a 2-transverse singularity but it is not a cusp. Hence, from Theorem \ref{Teo221}, $u_{o}$ is a swallow's tail\textit{ iff } condition  $(LS_{3})$ is true, i.e. \textit{iff} $\sigma_{3}(u_{o}) \notin R(F'(u_{o}))$. Moreover, from Remark \ref{Rem2111} we know that $\sigma_{3}(u_{o})=\varSigma_{3}[\varphi(u_{o}),\varphi_{1}(u_{o}),\varphi_{2}(u_{o})]$. Now let us consider $n_{o},n_{1},n_{2}$ such that they satisfy conditions (\ref{217}). By using $(s_{2})$ of Lemma \ref{Lem223} we get that 
$\varSigma_{3}[\varphi(u_{o}),\varphi_{1}(u_{o}),\varphi_{2}(u_{o})] \notin R(F'(u_{o})) $\textit{ iff }$\varSigma_{3}[n_{o},n_{1},n_{2}]\notin R(F'(u_{o}))$. Thus $u_{o}$ is a swallow's tail\textit{ iff } $\varSigma_{3}[n_{o},n_{1},n_{2}]\notin R(F'(u_{o}))$ as we wished to prove. \medskip\\
\indent \textit{3-transversality condition - Proof of} $\Leftarrow$. Let us consider $n_{o},n_{1},n_{2},v_{o},v_{1},v_{2}$ that satisfy conditions (\ref{217}) and such that (\ref{218}) is true, i.e. $T_{3}[n_{o},n_{1},n_{2},v_{o},v_{1},v_{2}]\notin R(F'(u_{o}))$. From the sixth and eighth conditions of (\ref{217}) we get that $T_{1}[n_{o},v_{o}]$ and $T_{2}[n_{o},n_{1},v_{o},v_{1}]$ belong to $ R(F'(u_{o}))$. These facts and the equivalences $(t_{o}), (t_{1})$ and $(t_{2})$ of Lemma \ref{Lem223} yield that 
\begin{align*}
T_{1}[\varphi(u_{o}),v_{o}]&\in R(F'(u_{o})), \\
T_{2}[\varphi(u_{o}),\varphi_{1}(u_{o}),v_{o},\varphi'(u_{o})v_{o}]&\in R(F'(u_{o})), \\
T_{3}[\varphi(u_{o}),\varphi_{1}(u_{o}),\varphi_{2}(u_{o}),v_{o},\varphi'(u_{o})v_{o},\varphi'_{1}(u_{o})v_{o}]&\notin R(F'(u_{o})). \qquad \qquad
\end{align*}
Then Remark \ref{Rem2111} implies that $\tau_{1}(u_{o})v_{o}\in R(F'(u_{o})), \,\tau_{2}(u_{o})v_{o}\in R(F'(u_{o}))$ and $\tau_{3}(u_{o})v_{o}\notin R(F'(u_{o}))$. By choosing $v=v_{o}$ we get that $(LT_{3})$ of Theorem \ref{Teo221} is verified and hence we can deduce that $u_{o}$ is a 3-transverse singularity. \medskip \\
\indent \textit{3-transversality condition - Proof of} $\Rightarrow$. Here we assume that $u_{o}$ is a 3-transverse singularity. Then, by $(LT_{3})$ of Theorem \ref{Teo221}, there exists $v_{o}\in X$ such that $\tau_{1}(u_{o})v_{o}\in R(F'(u_{o})), \,\tau_{2}(u_{o})v_{o}\in R(F'(u_{o}))$ and $\tau_{3}(u_{o})v_{o}\notin R(F'(u_{o}))$. From Remark \ref{2111} one gets that $T_{1}[\varphi(u_{o}),v_{o}]\in R(F'(u_{o})), \, T_{2}[\varphi(u_{o}),\varphi_{1}(u_{o}),v_{o},\varphi'(u_{o})v_{o}]\in R(F'(u_{o}))$ and $T_{3}[\varphi(u_{o}),\varphi_{1}(u_{o}),\varphi_{2}(u_{o}),v_{o},\varphi'(u_{o})v_{o},\varphi'_{1}(u_{o})v_{o}]\notin R(F'(u_{o}))$. Let us consider $n_{o},n_{1},n_{2}$ which satisfy the first five conditions of (\ref{217}). Since $u_{o}\in S_{1}(F)$ we can apply Lemma \ref{Lem223} for $k=0$ and, by  equivalence $(t_{o})$, it follows that  
$T_{1}[n_{o},v_{o}]\in R(F'(u_{o}))$. Then there exists $v_{1}\in X$ such that $F'(u_{o})v_{1} = -T_{1}[n_{o},v_{o}]$ and we obtain that $n_{o},n_{1},v_{o},v_{1}$ satisfy all the conditions of (\ref{214}). On the other hand, $u_{o}$ is a 1-transverse singularity which is not a 1-singularity and thus, from $(t_{1})$ of Lemma \ref{Lem223}, we get that $T_{2}[n_{o},n_{1},v_{o},v_{1}]\in R(F'(u_{o}))$. Hence there exists $v_{2}\in X$ such that $F'(u_{o})v_{2} = -T_{2}[n_{o},n_{1},v_{o},v_{1}]$ and $n_{o},n_{1},n_{2},v_{o},v_{1},v_{2}$ satisfy all the conditions of (\ref{217}). Finally, since $u_{o}$ is also a 2-transverse singularity which is not a 2-singularity we can use Lemma \ref{Lem223} for $k=2$ and, by equivalence $(t_{2})$, we can conclude that $T_{3}[n_{o},n_{1},n_{2},v_{o},v_{1},v_{2}]\notin R(F'(u_{o}))$  which is equivalent to stating that condition (\ref{218}) is satisfied. $\bracevert$\\
\par
\textbf{ Proof of Theorem \ref{Teo214}.} Let $(\varphi,\psi) \in \mathcal{P}(F,u_{o})$ be a given f-pair. Since $u_{o}$ is a 3-transverse singularity and it is not a swallow's tail, identities $(i_{1}), (i_{2})$ and $(i_{3})$ of Theorem \ref{Teo221} are true and thus $F'(u_{o})\varphi_{1}(u_{o})= -\sigma_{1}(u_{o}), F'(u_{o})\varphi_{2}(u_{o})= -\sigma_{2}(u_{o})$ and $F'(u_{o})\varphi_{3}(u_{o})= -\sigma_{3}(u_{o})$.  On the other hand, from Remark \ref{Rem2111} we have that $\sigma_{1}(u_{o})=\varSigma_{1}[\varphi(u_{o})], \sigma_{2}(u_{o})=\varSigma_{2}[\varphi (u_{o}),\varphi_{1} (u_{o})]$ and $\sigma_{3}(u_{o}) =\varSigma_{3}[\varphi(u_{o}),\varphi_{1}(u_{o}),\varphi_{2}(u_{o})]$. Hence there exist vectors $n_{o}:=\varphi(u_{o}), n_{1}:=\varphi_{1}(u_{o}), n_{2}:=\varphi_{2}(u_{o}), n_{3}:=\varphi_{3}(u_{o})$  which satisfy the conditions of (\ref{2110}). We can now prove the butterfly condition. \medskip\\
\indent\textit{ Butterfly condition.} By hypothesis $u_{o}$ is a 3-transverse singularity which is not a swallow's tail. By Theorem \ref{Teo221} we know that $u_{o}$ is a butterfly\textit{ iff } equivalence $(LS_{4})$ holds, i.e.\textit{ iff }$\sigma_{4}(u_{o})\notin R(F'(u_{o}))$.  Moreover, Remark \ref{Rem2111} implies that $\sigma_{4}(u_{o})=\varSigma_{4}[\varphi(u_{o}),\varphi_{1}(u_{o}),\varphi_{2}(u_{o}),\varphi_{3}(u_{o})]$ and for $n_{o},n_{1},n_{2},n_{3}$ satisfying the conditions (\ref{2110}) we obtain that $\varSigma_{4}[\varphi(u_{o}),\varphi_{1}(u_{o}),\varphi_{2}(u_{o}),\varphi_{3}(u_{o})]\notin R(F'(u_{o}))$\textit{ iff }$\varSigma_{4}[n_{o},n_{1},n_{2},n_{3}]\notin R(F'(u_{o}))$ from  equivalence $(s_{3})$ of Lemma \ref{Lem223}. Thus $u_{o}$ is a butterfly\textit{ iff }$\varSigma_{4}[n_{o},n_{1},n_{2},n_{3}]\notin R(F'(u_{o}))$.$\bracevert$\\
\vspace{2pt}
\subsection{Local Equations on the Strata of Singularities} \label{ss23}
\quad  The aim of this section is to prove Theorem \ref{Teo221}. To this end we first introduce an important tool given by Theorem \ref{Teo231}. This result, under suitable hypotheses of transversality, provides some peculiar identities on each singular stratum $S_{1_{k}}(F)$ (for the definition and properties of the singular strata cf. Section 2.4 in \cite{B-D 1}). These equalities are formulated in terms of the maps $\varphi_{k},\sigma_{k}$ and $\tau_{k}$ introduced in Definitions \ref{De217} and \ref{De218}. We note that in the proof of Theorems \ref{Teo221} and \ref{Teo231} an extensive use is made of some known facts concerning the characterization of singular strata; these facts were summarized in Remark 2.4.6 and Theorem 2.5.4 of \cite{B-D 1}.
\begin{theorem}\label{Teo231}
\textit{Let} $U,V$ \textit{be open subsets in the} $B$\textit{-spaces }$X,Y$\textit{ and }$F:U\subseteq X \rightarrow V \subseteq Y$\textit{ a }$C^{\infty}$ 0-\textit{Fredholm map. Let }$u_{o}\in S_{1}(F)$\textit{ and }$(\varphi,\psi) \in \mathcal{P}(F,u_{o})$.\textit{ Let us suppose that for some integer }$k\geq 1$\textit{ the point }$u_{o}$\textit{ is a }$k$-transverse singularity\textit{ for }$F$.\textit{ Then, on a suitable neighbourhood of }$u_{o}$,\textit{ for }$h=1,\ldots,k$,\textit{ the following equalities hold: }
\begin{align*}
&(j_{h})\quad &F'(u)\varphi^{\prime}_{h-1}(u)v = -\tau_{h}(u)v \; ,\quad\, \forall \, u \in S_{1_{h}}(F)\, , \, \forall \, v\in T_{u}S_{1_{h}}(F);\\
\\
&(jj_{h})\quad &I_{h+1}(u)v = \psi(u)\tau_{h+1}(u)v\;,\quad \forall \, u \in S_{1_{h+1}}(F)\, , \, \forall \, v\in T_{u}S_{1_{h}}(F).\\
\end{align*}
\end{theorem}
\textbf{Proof. }The proof is by induction and it will be divided in three steps. In the first one we deduce some useful identities, valid in a neighbourhood of $u_{o}$ and without assumptions of transversality. In the second step we analyse the case where $u_{o}$ is a 1-transverse singularity, while in the third one we study what happens when $u_{o}$ is $k$-transverse for $k \geq 2$.\\
\par \textit{ Step 1.} Let $u_{o}$ be a simple singularity for $F$ and let us suppose that $(\varphi,\psi)$ is an f-pair for $F$ on a suitable neighbourhood $U_{o}$ of $u_{o}$ which we can always assume, up to shrinking, to be the whole of $U$. Then, according to the definition of the fibering functionals $J_{k}$ and $I_{k}$ (cf. Definition 1.1.5 in \cite{B-D 1}), we consider the functional
\begin{equation*}
J_{o}(u)=\psi(u)F'(u)\varphi(u), \, u\in U;
\end{equation*}
by differentiating this identity we obtain for $u\in U$ and $v\in X$ that
\begin{equation*}
I_{1}(u)v=J^{\prime}_{o}(u)v = (\psi^{\prime}(u)v)F'(u)\varphi(u)+ \psi(u)F''(u)[\varphi(u),v]+\psi(u)F'(u)(\varphi^{\prime}(u)v).
\end{equation*}
Moreover, by the very definition of the fibering maps $\varphi$ and $\psi$ we have that
\begin{equation}\label{231}
F'(u)\varphi(u)= 0 \in Y,\;\;  u\in S_{1}(F),
\end{equation}
\begin{equation}\label{232}
\psi(u)F'(u) = 0 \in X^{\ast},\;\;  u\in S_{1}(F);
\end{equation}
hence, by the above formulas and  Definition \ref{De218}, we get the equalities
\begin{equation}\label{233}
I_{1}(u)v=\psi(u)F''(u)[\varphi(u),v]=\psi(u)\tau_{1}(u)v, \;\; u \in S_{1}(F), v \in X,
\end{equation}
\begin{equation}\label{234}
J_{1}(u)=I_{1}(u)\varphi(u)=\psi(u)\tau_{1}(u)\varphi(u)=\psi(u)\sigma_{1}(u)\,, \;\, u \in S_{1}(F).
\end{equation}
\par \textit{ Step 2.} Now let us suppose that $u_{o}$ is a 1-transverse singularity for $F$. Then, near $u_{o}$, we have that $S_{1_{1}}(F)=S_{1}(F)$ is a one-codimensional manifold of $X$, cf. Remark 2.4.6, d), of \cite{B-D 1}. By differentiating the equalities (\ref{231}), (\ref{232}), (\ref{234}) on the manifold $S_{1}(F)$ near $u_{o}$ we obtain, for $u \in S_{1}(F), v \in T_{u}S_{1}(F)$:
\begin{equation}\label{235}
\qquad\qquad(F'(u)\varphi(u))'v = F''(u)[\varphi(u),v]+F'(u)(\varphi^{\prime}(u)v)= 0 \in Y,
\end{equation}
\begin{equation}\label{236}
\qquad\qquad(\psi(u)F'(u))'v=(\psi^{\prime}(u)v)F'(u)+\psi(u)F''(u)[v,\cdotp] = 0 \in X^{\ast},
\end{equation}
\begin{equation}\label{237}
\qquad \;I_{2}(u)v=J^{\prime}_{1}(u)v = (\psi(u)\sigma_{1}(u))'v = (\psi'(u)v)\sigma_{1}(u)+\psi(u)(\sigma^{\prime}_{1}(u)v).
\end{equation}
If we write (\ref{235}) in the form
\begin{equation}\label{238}
F'(u)\varphi'(u)v = -F''(u)[\varphi(u),v] = -\tau_{1}(u)v \; , \;\; u \in S_{1}(F), v \in T_{u}S_{1}(F)
\end{equation}
and recall that $S_{1_{1}}(F)=S_{1}(F)$ and $\varphi_{o} = \varphi$ we obtain the identity $(j_{1})$ of the statement. Furthermore, (\ref{236}) can be rewritten as the following equality in $X^{\ast}$:
 \begin{equation*}
(\psi'(u)v)F'(u)= -\psi(u)F''(u)[v,\cdotp].
\end{equation*}
Applying this identity to the elements $\varphi_{m}(u)\in X$ (cf. Definition \ref{De217}) for all $m\geq 0$ we get another identity 
\begin{equation}\label{239}
(\psi'(u)v)F'(u)\varphi_{m}(u)= -\psi(u)F''(u)[v,\varphi_{m}(u)]\, , \; u \in S_{1}(F), v \in T_{u}S_{1}(F).
\end{equation}
\indent If $S_{1_{2}}(F)= \emptyset$ there is nothing to prove. When $S_{1_{2}}(F)\neq \emptyset$ let us consider $u$ near $u_{o}$ such that $u \in S_{1_{2}}(F)\subseteq S_{1_{1}}(F) = S_{1}(F)$. This amounts to saying that $\varphi(u) \in T_{u}S_{1_{1}}(F)$ as shown in Theorem 2.5.4, c), of \cite{B-D 1}. Then by applying (\ref{238}) with $v=\varphi(u)$ we can write that 
\begin{equation*}
F'(u)\varphi_{1}(u)= F'(u)\varphi'(u)\varphi(u) = -\tau_{1}(u)\varphi(u) = -\sigma_{1}(u) \,,\;u \in S_{1_{2}}(F),  
\end{equation*}
where the first and the last equality follow from Definitions \ref{De217} and \ref{De218}.
This last equality and (\ref{239}) give  
 \begin{align}\label{2310}
\begin{split}
(\psi'(u)v)\sigma_{1}(u) &=-(\psi'(u)v)F'(u)\varphi_{1}(u)=\\
&=\psi(u)F''(u) [v,\varphi_{1}(u)], u\in S_{1_{2}}(F), v \in T_{u}S_{1_{1}}(F).
\end{split}
\end{align} 
If we now use  (\ref{237}), (\ref{2310}) and recall Definition \ref{De218} we get
\begin{align*}
\begin{split}
I_{2}(u)v &= \psi(u)F''(u)[v,\varphi_{1}(u)]+ \psi(u)(\sigma^{\prime}_{1}(u)v)=\\
&= \psi(u)(F''(u)[v,\varphi_{1}(u)]+ \sigma^{\prime}_{1}(u)v)= \psi(u)\tau_{2}(u)v \,,\; u\in S_{1_{2}}(F), v \in T_{u}S_{1_{1}}(F),
\end{split}
\end{align*} 
that is equality $(jj_{1})$. Therefore the theorem is proved when $k=1$.\\
\par \textit{ Step 3.} Let us assume that $u_{o}$ is a $k$-transverse singularity for $k\geq 2$. We shall prove formulas $(j_{h}), (jj_{h})$ by induction on $h$. Since $u_{o}$ is also 1-transverse we know that such equalities hold for $h=1$, as was proved in the previous step. Now we suppose that $(j_{h}),(jj_{h})$ hold for an integer $h$, with $1\leq h\leq k-1$ and we will prove the corresponding identities for $h+1$.\\
\indent By hypothesis the equality $(j_{h})$ is true, i.e. \\
\begin{equation*}
F'(u)\varphi^{\prime}_{h-1}(u)v=-\tau_{h}(u)v \,,\; u\in S_{1_{h}}(F), \, v \in T_{u}S_{1_{h}}(F).
\end{equation*}
Since $u\in S_{1_{h+1}}(F)$ is equivalent to saying that $u\in S_{1_{h}}(F)$ and $\varphi(u)\in T_{u}S_{1_{h}}(F)$, cf. Theorem 2.5.4, c), of \cite{B-D 1}, we have the following identity:
\begin{equation}\label{2311}
F'(u)\varphi_{h}(u)=F'(u)\varphi^{\prime}_{h-1}(u)\varphi(u)=-\tau_{h}(u)\varphi(u)= -\sigma_{h}(u), \; u\in S_{1_{h+1}}(F)
\end{equation}
where, again, we used Definitions \ref{De217} and \ref{De218}.
Now we recall that $u_{o}$ is $k$-transverse and $1\leq h\leq k-1$, hence each stratum $S_{1_{h+1}}(F)$ is a manifold near $u_{o}$ (cf. \cite{B-D 1}, Remark 2.4.6). Therefore we can differentiate equality (\ref{2311}) on the manifold $S_{1_{h+1}}(F)$ and, in this way, for all $u\in S_{1_{h+1}}(F)$ and $v \in T_{u}S_{1_{h+1}}(F)$ we obtain that
\begin{equation*}
(F'(u)\varphi_{h}(u))'v = F''(u)[v,\varphi_{h}(u)]+ F'(u)(\varphi^{\prime}_{h}(u)v) = -\sigma^{\prime}_{h}(u)v.
\end{equation*}
The last equality can be rewritten, by means of  Definition \ref{De218}, as
\begin{equation}\label{2312}
F'(u)\varphi^{\prime}_{h}(u)v = -\sigma^{\prime}_{h}(u)v -F''(u)[v,\varphi_{h}(u)]= -\tau_{h+1}(u)v
\end{equation}
for $u\in S_{1_{h+1}}(F)$, $v \in T_{u}S_{1_{h+1}}(F)$. This proves identity $(j_{h+1})$.\\
\indent  Once again, by Theorem 2.5.4, c), of \cite{B-D 1} we know that for $1\leq h\leq k-1$ the stratum $S_{1_{h+2}}(F)$ is given by all $u\in S_{1_{h+1}}(F)$ such that $\varphi(u)\in T_{u}S_{1_{h+1}}(F)$ (where we are assuming that $S_{1_{k+1}}(F)$ is not empty because otherwise there is nothing to prove for $h=k-1$). In particular, thanks to Definitions \ref{De217} and \ref{De218}, equality (\ref{2312}) with $v=\varphi(u)$ yields 
\begin{equation*}
F'(u)\varphi_{h+1}(u) = F'(u)\varphi^{\prime}_{h}(u)\varphi(u)=-\tau_{h+1}(u)\varphi(u)=-\sigma_{h+1}(u), u\in S_{1_{h+2}}(F).
\end{equation*}
If we set $m=h+1$ this last identity and (\ref{239}) allow us to write
\begin{align}\label{2313}
\begin{split}
(\psi'(u)v)\sigma_{h+1}(u) &=-(\psi'(u)v)F'(u)\varphi_{h+1}(u)=\\
&=\psi(u)F''(u) [\varphi_{h+1}(u),v], u\in S_{1_{h+2}}(F), v \in T_{u}S_{1}(F)
\end{split}
\end{align}
with $S_{1_{h+2}}(F)\subseteq S_{1_{1}}(F) = S_{1}(F)$.\\
By hypothesis equality $(jj_{h})$ is true, i.e. for $1\leq h\leq k-1$ we have that
\begin{equation}\label{2314}
I_{h+1}(u)v=\psi(u)\tau_{h+1}(u)v,\, u\in S_{1_{h+1}}(F), \, v\in T_{u}S_{1_{h}}(F).
\end{equation}
Then, according to the definition of the fibering functionals $J_{h}$ and $I_{h}$, we can write
\begin{equation*}
J_{h+1}(u) = I_{h+1}(u)\varphi(u)=\psi(u)\tau_{h+1}(u)\varphi(u)=\psi(u)\sigma_{h+1}(u)\,, \;\, u \in S_{1_{h+1}}(F),
\end{equation*}
where the second equality is just (\ref{2314}) expressed for $v=\varphi(u)$, since $u\in S_{1_{h+1}}(F)$\textit{ iff }$u \in S_{1_{h}}(F)$ with $\varphi(u)\in T_{u}S_{1_{h}}(F)$, and the last equality follows from Definition \ref{De218}.\\
By differentiating this equality on $S_{1_{h+1}}(F)$ and thanks to Definition 1.1.5 in \cite{B-D 1} we obtain that, for $u\in S_{1_{h+1}}(F)$ and $v\in T_{u}S_{1_{h+1}}(F)$,
\begin{equation}\label{2315}
I_{h+2}(u)v = J^{\prime}_{h+1}(u)v = (\psi'(u)v)\sigma_{h+1}+\psi(u)\sigma^{\prime}_{h+1}(u)v.
\end{equation}
Recalling that $S_{1_{h+2}}(F)\subseteq S_{1_{h+1}}(F)$ and $T_{u}S_{1_{h+1}}(F)\subseteq T_{u}S_{1}(F)$, identity (\ref{2313}) can be used in (\ref{2315}) to deduce that, for $u \in S_{1_{h+2}}(F)$ and $v\in T_{u}S_{1_{h+1}}(F)$, 
\begin{align*}
\begin{split}
I_{h+2}(u)v &= \psi(u)F''(u)[\varphi_{h+1}(u),v]+\psi(u)\sigma^{\prime}_{h+1}(u)v=\\
&=\psi(u)(F''(u)[\varphi_{h+1}(u),v]+\sigma^{\prime}_{h+1}(u)v).
\end{split}
\end{align*}
By virtue of  Definition \ref{De218} we can conclude that 
\begin{equation*}
I_{h+2}(u)v= \psi(u)\tau_{h+2}(u)v \,,\, u \in S_{1_{h+2}}(F), \, v\in T_{u}S_{1_{h+1}}(F),
\end{equation*}
i.e. we proved formula $(jj_{h+1}).  \bracevert $\\
\par  In order to demonstrate Theorem \ref{Teo221} we still need a result, given by Proposition \ref{Pro233}, which characterizes the tangent space at each point of a singular stratum by means of the maps $\tau_{h}$. To this end it is convenient to rewrite the conditions of $(k+1)$-singularity and $(k+1)$-transversality (see Definition 2.1.1 in \cite{B-D 1}) as follows:
\begin{lemma}\label{Lem232} \textit{ Let} $U,V$ \textit{be open subsets in the} $B$\textit{-spaces }$X,Y$\textit{ and }$F:U\subseteq X \rightarrow V \subseteq Y$\textit{ a }$C^{\infty}$ 0-\textit{Fredholm map. Let }$u_{o}\in S_{1}(F)$\textit{ with }$(\varphi,\psi) \in \mathcal{P}(F,u_{o})$\textit{ and let us assume that, for some integer }$k\geq 1, u_{o}$\textit{ is }$k$-transverse\textit{ for }$F$. \textit{Then we have that }$u_{o}$\textit{ is }$h$-transverse\textit{ but not an }$h$-singularity\textit{ for all }$h$\textit{ such that }$1\leq h\leq k-1$. \textit{Moreover, if }$u_{o}$ is $k$-transverse\textit{ and not a }$k$-singularity\textit{ we have that}
 \begin{equation*}
u_{o}\textit{ is }(k+1)\text{-transverse} \Leftrightarrow \exists \, v\in \cap^{k}_{j=1}N(I_{j}(u_{o})): I_{k+1}(u_{o})v\neq 0,
\end{equation*}
\textit{and, in particular,}
\begin{center}
$ u_{o}\textit{ is a }(k+1)\text{-singularity} \Leftrightarrow J_{k+1}(u_{o})\neq 0.\hspace{130pt}$
\end{center}
\end{lemma}
\indent \textbf{Proof.} By combining conditions (T$_{h}$) and (S$_{h}$) in Definition 2.1.1 of \cite{B-D 1} it is easy to see that, for any integer $h, u_{o}$  is $h$-transverse and not an $h$-singularity for $F$ \textit{iff }
\begin{equation*}
(\text{TS}_{h})\quad \begin{cases}
J_{o}(\varphi,\psi)(u_{o})= \ldots = J_{h}(\varphi,\psi)(u_{o})= 0,\\
I_{1}(\varphi,\psi)(u_{o}),\ldots, I_{h}(\varphi,\psi)(u_{o}) \text{  are l.i. }.
\end{cases}
\end{equation*}
By hypothesis $u_{o}$ is $k$-transverse and thus condition (T$_{k}$) of Definition 2.1.1 in \cite{B-D 1} is satisfied. Therefore (TS$_{h}$) holds for any integer $h$ such that $1\leq h\leq k-1$.\\
\indent Let us now assume that $u_{o}$ is $k$-transverse and not a $k$-singularity. Hence condition (TS$_{k}$) is satisfied. A comparison between (TS$_{k}$) and (T$_{k+1}$) in Definition 2.1.1 of \cite{B-D 1} yields that
\begin{equation*}
u_{o}\text{ is }(k+1)\text{-transverse  }\textit{  iff }\; I_{1}(\varphi,\psi)(u_{o}),\ldots,I_{k}(\varphi,\psi)(u_{o}), I_{k+1}(\varphi,\psi)(u_{o})\text{ are l.i. .}
\end{equation*}
Because of the linear independence of $I_{1}(\varphi,\psi)(u_{o}),\ldots,I_{k}(\varphi,\psi)(u_{o})$ and the Algebraic Lemma (cf. Lemma 2.2.1, a) $\Leftrightarrow$ f), of \cite{B-D 1}), the previous equivalence amounts to saying that there exists $v\in \cap^{k}_{j=1}N(I_{j}(u_{o}))$ such that $I_{k+1}(u_{o})v \neq 0$, as we had to show. Finally, (TS$_{k}$) and condition (S$_{k+1}$) of  Definition 2.1.1 in \cite{B-D 1} imply that 
\begin{equation*}
u_{o}\text{ is a }(k+1)\text{-singularity  }\textit{  iff }\; J_{k+1}(u_{o})\neq 0. \;  \bracevert
\end{equation*}

\begin{proposition}\label{Pro233}
\textit{Let} $U,V$ \textit{be open subsets in the} $B$\textit{-spaces }$X,Y$\textit{ and }$F:U\subseteq X \rightarrow V \subseteq Y$\textit{ a }$C^{\infty}$ 0-\textit{Fredholm map. Let }$u_{o}\in S_{1}(F)$\textit{ with }$(\varphi,\psi) \in \mathcal{P}(F,u_{o})$\textit{ and let us assume that, for some integer }$k\geq 1, u_{o}$\textit{ is }$k$-transverse\textit{ for }$F$.\textit{ Then the following equality holds:}
\begin{equation*}
T_{u_{o}}S_{1_{k}}(F) = \cap^{k}_{j=1}N(I_{j}(u_{o})) =\{v \in X:\tau_{j}(u_{o})v \in R(F'(u_{o})),j=1,\ldots,k \}.
\end{equation*}
\end{proposition} 
\indent 
\textbf{Proof.} Let us first set $k=1$. We recall that, by formula (\ref{233}) in the proof of  Theorem \ref{Teo231}, we know that
\begin{equation*}
I_{1}(u)v = \psi(u)\tau_{1}(u)v \,,\; u \in S_{1}(F),\, v\in X.
\end{equation*}
Since by definition $\psi(u_{o})\in R(F'(u_{o}))^{\perp}$, it follows that
\begin{equation*}
N(I_{1}(u_{o})) =\{v \in X:\tau_{1}(u_{o})v \in R(F'(u_{o})) \}.
\end{equation*}
By hypothesis $u_{o}$ is 1-transverse and hence $T_{u_{o}}S_{1}(F)= N(I_{1}(u_{o})) $, cf. Remark 2.4.6, d), in \cite{B-D 1}: this concludes the proof for the case $k=1$.\\
\indent Let us now suppose that the result holds for a given integer $k\geq 1$ and let us prove it for the integer $k+1$. Let $u_{o}$ be a $(k+1)$-transverse singularity. By the preceding Lemma \ref{Lem232} we know, in particular, that $u_{o}$ is $k$-transverse and not a $k$-singularity and thus, by inductive hypothesis, we have that
\begin{equation}\label{2316}
T_{u_{o}}S_{1_{k}}(F)=\{v \in X:\tau_{j}(u_{o})v \in R(F'(u_{o})),j=1,\ldots,k \}. 
\end{equation}
Moreover, from Remark 2.4.6, d), of \cite{B-D 1}, we get that $u_{o} \in S_{1_{k+1}}(F)$ and
\begin{equation}\label{2317}
\begin{split}
T_{u_{o}}S_{1_{k+1}}(F)&=\cap^{k+1}_{j=1}N(I_{j}(u_{o}))=\{v \in \cap^{k}_{j=1}N(I_{j}(u_{o})):I_{k+1}(u_{o})v= 0\}=\\
&=\{v \in T_{u_{o}}S_{1_{k}}(F):I_{k+1}(u_{o})v= 0 \}.
\end{split}
\end{equation}
But, as we already know, $u_{o}$ is $k$-transverse and not a $k$-singularity. Thus, thanks to  identity $(jj_{k})$ of Theorem \ref{Teo231},  we obtain 
\begin{equation*}
I_{k+1}(u_{o})v=\psi(u_{o})\tau_{k+1}(u_{o})v \,,\; v\in T_{u_{o}}S_{1_{k}}(F).
\end{equation*}
Since (from the definition of $\psi(u_{o})$) we know that
\begin{equation*}
I_{k+1}(u_{o})v=0 \quad\textit{iff}\quad\tau_{k+1}(u_{o})v \in R(F'(u_{o}))\, ,
\end{equation*}
equalities (\ref{2316}) and (\ref{2317}) allow us to conclude that 
\begin{equation*}
T_{u_{o}}S_{1_{k+1}}(F) = \{v \in X:\tau_{j}(u_{o})v \in R(F'(u_{o})),j=1,\ldots,k+1 \} ,
\end{equation*}
as we had to show.$\bracevert$
\begin{remark}\label{Rem234}
By means of the above Proposition one can rewrite formula $(j_{h})$ of Theorem \ref{Teo231} in a way that will be useful in the next section. Namely, let $u_{o}$ be a $k$-transverse singularity for $F$. Then $u_{o}\in S_{1_{h}}(F)$ for all integers $h$ such that $h=1,\ldots,k$. Thus, for a given $v \in X$, we have that
 \begin{equation*}
(j^{\,\prime}_{h})\qquad\qquad\tau_{j}(u_{o})v \in R(F'(u_{o})),j=1,\ldots,h \Rightarrow F'(u_{o})\varphi^{\prime}_{h-1}(u_{o})v = -\tau_{h}(u_{o})v . 
\end{equation*}
This is true because, by Proposition \ref{Pro233}, we know that $\tau_{j}(u_{o})v \in R(F'(u_{o})),j=1,\ldots,h,\;\textit{iff }\; v\in T_{u_{o}}S_{1_{h}}(F)$ and this allows to apply $(j_{h})$ and obtain $(j^{\,\prime}_{h})$.
\end{remark}

We are now able to prove Theorem \ref{Teo221}.\\
\par \textbf{Proof of Theorem 2.2.1.} This demonstration is in three steps. In the first part of the proof we obtain equivalences $(LT_{1}),(LS_{1})$; in the second one we show that equivalences $(LT_{k+1}),(LS_{k+1})$ for $k\geq 1$ are true, and in the third part we prove identity $(i_{h})$.\\
\par \textit{Part 1.} We assume that $u_{o}$ is a simple singularity. We saw in the proof of Proposition \ref{Pro233} that
\begin{equation*}
N(I_{1}(u_{o})) =\{v \in X:\tau_{1}(u_{o})v \in R(F'(u_{o})) \}.
\end{equation*}
Here we first have to show that 
\begin{equation*}
u_{o} \text{ is a 1-transverse singularity } \Leftrightarrow \;(LT_{1})\!:\;\exists \;\,v \in X \text{ such that }\tau_{1}(u_{o})v \notin R(F'(u_{o})).
\end{equation*}
By the very definition of 1-transverse singularity (cf. Definition 2.1.1 of \cite{B-D 1}), $u_{o}$ is 1-transverse \textit{iff}  $I_{1}(u_{o})$ is l.i., that is $I_{1}(u_{o})\neq 0$: hence the above statement is true. Finally we need to prove that
\begin{equation*}
u_{o} \text{ is a 1-singularity } \Leftrightarrow \;(LS_{1})\!:\; \sigma_{1}(u_{o}) \notin R(F'(u_{o})).
\end{equation*}
By Definition 2.1.1 of \cite{B-D 1} one has that $u_{o}$ is a 1-singularity \textit{iff} $J_{1}(u_{o})\neq 0$. From the very definition of $J_{1}$ we have that  $J_{1}(u_{o}) = I_{1}(u_{o})\varphi(u_{o})$ (cf. Definition 1.1.5 of \cite{B-D 1}), hence $J_{1}(u_{o})\neq 0$ \textit{ iff }$\varphi(u_{o})\notin N(I_{1}(u_{o}))\textit{ iff }\,\tau_{1}(u_{o})\varphi(u_{o}) \notin R(F'(u_{o}))$. By Definition \ref{De218} we have that $\sigma_{1}(u_{o}) = \tau_{1}(u_{o})\varphi(u_{o})$, therefore the above statement is proved.\\
\par 
\textit{Part 2}. Let us now suppose that, for some integer $k\geq 1,\; u_{o}$ is $k$-transverse and not a $k$-singularity. We wish to show that
\begin{equation*}
u_{o}\textit{ is a }(k+1)\text{-transverse singularity}  \Leftrightarrow (LT_{k+1}):
\begin{cases}
\exists \; v\in X:\tau_{1}(u_{o})v \in R(F'(u_{o})),\\
\quad \quad , \ldots,\tau_{k}(u_{o})v \in R(F'(u_{o})),\\
\qquad \quad \; \tau_{k+1}(u_{o})v \notin R(F'(u_{o})).
\end{cases}
\end{equation*}
The hypotheses on $u_{o}$ imply that $u_{o}\in S_{1_{k+1}}(F)$ (cf. Theorem 2.5.4, c), in \cite{B-D 1}). Moreover, from Lemma \ref{Lem232} we know that
\begin{equation*}
u_{o} \text{ is } (k+1)\text{-transverse } \Leftrightarrow \;\exists \;\; v \in \cap^{k}_{j=1}N(I_{j}(u_{o})):I_{k+1}(u_{o})v \neq 0
\end{equation*}
and, by Proposition \ref{Pro233}, we have
\begin{equation*}
T_{u_{o}}S_{1_{k}}(F)= \cap^{k}_{j=1}N(I_{j}(u_{o})) =\{v \in X:\tau_{j}(u_{o})v \in R(F'(u_{o})),j=1,\ldots,k \}. 
\end{equation*}
Thus $u_{o}$ is $(k+1)$-transverse \textit{iff} there exists $v \in X$ such that $\tau_{j}(u_{o})v \in R(F'(u_{o})),j=1,\ldots,k$, and $I_{k+1}(u_{o})v \neq 0$. Since $u_{o}\in S_{1_{k+1}}(F)$ it follows from Theorem \ref{Teo231}, $(jj_{k})$, that
 \begin{equation}\label{2318}
I_{k+1}(u_{o})v = \psi(u_{o})\tau_{k+1}(u_{o})v\, , \; v\in T_{u_{o}}S_{1_{k}}(F).
\end{equation}
In order to conclude the proof of $(LT_{k+1})$ it then suffices to recall that $\psi(u_{o})\tau_{k+1}(u_{o})v \neq 0$ \textit{ iff } $\tau_{k+1}(u_{o})v \notin R(F'(u_{o}))$. We then have to prove that
\begin{equation*}
u_{o}\text{ is a }(k+1)\text{-singularity} \;  \Leftrightarrow \; (LS_{k+1})\!:\, \sigma_{k+1}(u_{o})\notin R(F'(u_{o})).
\end{equation*}
From Lemma \ref{Lem232} we know that
\begin{equation*}
u_{o}\text{ is a }(k+1)\text{-singularity} \;  \Leftrightarrow \; J_{k+1}(u_{o})\neq 0.
\end{equation*}
Given that $u_{o} \in S_{1_{k+1}}(F)$, $u_{o} \in S_{1_{k}}(F)$ and $\varphi(u_{o})\in T_{u_{o}}S_{1_{k}}(F)$ (cf. Theorem 2.5.4, c), in \cite{B-D 1}). By the very definition of $J_{k+1}$ (cf. Definition 1.1.5 in \cite{B-D 1}) and  formula (\ref{2318}) we get 
\begin{equation*}
J_{k+1}(u_{o})= I_{k+1}(u_{o})\varphi (u_{o})=\psi(u_{o})\tau_{k+1}(u_{o})\varphi(u_{o}).
\end{equation*}
Hence $J_{k+1}(u_{o})\neq 0\;\textit{ iff }\; \tau_{k+1}(u_{o})\varphi(u_{o}) \notin R(F'(u_{o}))$. By  Definition \ref{De218} $\sigma_{k+1}(u_{o}) = \tau_{k+1}(u_{o})\varphi(u_{o}) $, and equivalence $(LS_{k+1})$ is thus proved.\\
\par
\textit{Part 3.} We will show that for $h=1,\ldots, k$ the elements $\varphi_{h}(u_{o})$ satisfy the equations
\begin{equation*}
(i_{h}) \hspace{70pt} F'(u_{o})\varphi_{h}(u_{o}) = -\sigma_{h}(u_{o}).
\end{equation*}
As we know from Part 2, the hypotheses \textquotedblleft$u_{o}$ $k$-transverse\textquotedblright and \textquotedblleft$u_{o}$ not a $k$-singularity\textquotedblright $\,$ imply that $u_{o}\in S_{1_{k+1}}(F)$ and this, in turn, is equivalent to saying that $u_{o}\in S_{1_{k}}(F)$ and $\varphi(u_{o})\in T_{u_{o}}S_{1_{k}}(F)$. From the definition of singular strata (cf. Definition 2.4.1 in \cite{B-D 1}) we know that  $S_{1_{k+1}}(F)\subseteq S_{1_{h}}(F)$ for $h=1,\ldots,k$. Moreover, from d) of Remark 2.4.6 in \cite{B-D 1} we also have that $S_{1_{h}}(F), h=1,\ldots,k$, are manifolds and, of course, $T_{u_{o}}S_{1_{k}}(F)\subseteq T_{u_{o}}S_{1_{h}}(F)$ for $h=1,\ldots,k$. Thus we can use equality $(j_{h})$ of Theorem \ref{Teo231} with $v = \varphi(u_{o})$, that is we can write
\begin{equation*}
F'(u_{o})\varphi^{\prime}_{h-1}(u_{o})\varphi(u_{o}) = -\tau_{h}(u_{o})\varphi(u_{o}), h=1,\ldots,k. 
\end{equation*}
By Definitions \ref{De217} and \ref{De218}, $\varphi^{\prime}_{h-1}(u_{o})\varphi(u_{o}) = \varphi_{h}(u_{o})$ and $\tau_{h}(u_{o})\varphi(u_{o}) = \sigma_{h}(u_{o})$ and thus equality $(i_{h})$ follows at once. $\bracevert$\\
\vspace{2pt}
\subsection{Proof of the Combinatorial Lemma} \label{ss24}
\quad   We now wish to prove the equivalences stated in Lemma 2.2.3, i.e. we have to show that
\begin{align*}
(s_{h}) \quad &\varSigma_{h+1}[\varphi(u_{o}),\ldots,\varphi_{h}(u_{o})]\in R(F'(u_{o})) \Leftrightarrow \varSigma_{h+1}[n_{o},\ldots,n_{h}]\in R(F'(u_{o})),\\
\\
(t_{h}) \quad &T_{h+1}[\varphi(u_{o}),\ldots,\varphi_{h}(u_{o}),v_{o},\ldots,\varphi'_{h-1}(u_{o})v_{o}]\in R(F'(u_{o}))\\
& \hspace{175pt}\Leftrightarrow T_{h+1}[n_{o},\ldots,n_{h},v_{o},\ldots,v_{h}]\in R(F'(u_{o}))
\end{align*}
are satisfied for any $h = 0,\ldots, k$, where $k = 0, 1, 2, 3$, under the assumptions that $u_{o}$ is a $k$-transverse singularity which is not a $k$-singularity (or just a simple singularity for $k=0$) and that there exist $n_{o},\ldots,n_{k},v_{o},\ldots,v_{k}$ which satisfy conditions (\ref{211}), (\ref{214}), (\ref{217}) and (\ref{2110}) of Theorems \ref{Teo211}, \ref{Teo212}, \ref{Teo213} and \ref{Teo214} respectively.\\
\indent Given a fibering pair $(\varphi,\psi)$ for $F$ near $u_{o}$, in proving each step $(s_{k}), k = 0, 1, 2, 3$, the following inductive pattern shall be used: \\
a) we write $\varphi_{k}(u_{o})$ as a linear combination of the vectors $n_{h},h=0,\ldots,k$. This is possible because of formula $(i_{k})$ of Theorem \ref{Teo221} and because a similar relation for $\varphi_{k-1}(u_{o})$ and $n_{h}, h=0,\ldots,k-1$, is found when proving the previous step;\\
b) we write $\varSigma_{k+1}[\varphi(u_{o}),\ldots,\varphi_{k}(u_{o})]$ as a linear combination of  $\varSigma_{h+1}[n_{o},\ldots,n_{h}], h=0,\ldots,k$. This is achieved by using all the identities between $\varphi_{h}(u_{o})$ and $n_{j}, j = 0,\ldots, h$, for $h = 0,\ldots, k$.\\
\indent A similar procedure will give us the proof of the steps $(t_{k})$ for $k = 0, 1, 2, 3$. In this case we will also use relations $(j^{\,\prime}_{h})$ of Remark \ref{Rem234}.\\
\par
\textbf{Proof of Lemma \ref{Lem223} for k = 0.} We want to show that $(s_{o}), (t_{o})$ hold if $u_{o}\in S_{1}(F)$.\\
\indent \textit{Proof of }$(s_{o})$. We recall that dim\,$N(F'(u_{o}))=1$ and $\varphi(u_{o})$ spans $N(F'(u_{o}))$. Thus, if we choose $n_{o}\in N(F'(u_{o}))\setminus \{0\}$, there exists a unique non-zero real number $\alpha_{o}$ such that   
\begin{equation}\label{241}
\varphi(u_{o})=\alpha_{o}n_{o}.
\end{equation}
By definition,  
\begin{equation}\label{242}
\varSigma_{1}[\varphi(u_{o})]=F''(u_{o})[\varphi(u_{o}),\varphi(u_{o})]=\alpha^{2}_{o}F''(u_{o})[n_{o},n_{o}]=\alpha^{2}_{o}\varSigma_{1}[n_{o}]
\end{equation}
hence $(s_{o})$ is satisfied.\\
\indent \textit{Proof of }$(t_{o})$. Given  $n_{o}\in N(F'(u_{o}))\setminus \{0\}$, we have that $\varphi(u_{o})=\alpha_{o}n_{o}$ with $\alpha_{o}\in \mathbb{R}\setminus \{0\}$. Then if $v_{o}\in X$, i.e. the pair $(n_{o},v_{o})$ satisfies condition (\ref{211}) of Theorem \ref{Teo211}, we have by definition of the map $T_{1}$ that
\begin{equation}\label{243}
T_{1}[\varphi(u_{o}),v_{o}]=F''(u_{o})[\varphi(u_{o}),v_{o}]=\alpha_{o}F''(u_{o})[n_{o},v_{o}]=\alpha_{o}T_{1}[n_{o},v_{o}]
\end{equation}
and this concludes the proof.  $\bracevert$\\
\par
\textbf{Proof of Lemma \ref{Lem223} for k = 1}. We have to show that $(s_{o}), (t_{o}),(s_{1}),(t_{1})$ hold if $u_{o}$ is a 1-transverse singularity and is not a 1-singularity. Since $u_{o}\in S_{1}(F)$, from the previous step we already know that $(s_{o}), (t_{o})$ hold along with formulas (\ref{241}), (\ref{242}) and (\ref{243}). Hence it suffices to prove that $(s_{1}), (t_{1})$ are satisfied.\\
\indent \textit{Proof of }$(s_{1})$. Since $u_{o}$ is a 1-transverse singularity which is not a 1-singularity we know, by $(i_{1})$ of Theorem \ref{Teo221}, that $F'(u_{o})\varphi_{1}(u_{o})=-\sigma_{1}(u_{o})$. This is equivalent to $F'(u_{o})\varphi_{1}(u_{o})=-\varSigma_{1}[\varphi(u_{o})]$ thanks to Remark \ref{Rem2111}. Let $n_{o}, n_{1}$ satisfy condition (\ref{214}) of Theorem \ref{Teo212} , i.e. $n_{o}\in N(F'(u_{o}))\setminus \{0\}$ and $
F'(u_{o})n_{1}=-\varSigma_{1}[n_{o}]$. The relation between $\varphi_{1}(u_{o})$ and $n_{1}$ easily follows by recalling that $\varphi(u_{o})=\alpha_{o}n_{o},\alpha_{o} \in \mathbb{R}\setminus \{0\}$, and
\begin{equation*}
F'(u_{o})\varphi_{1}(u_{o})=-\varSigma_{1}[\varphi(u_{o})]=-\alpha^{2}_{o}\varSigma_{1}[n_{o}]=\alpha^{2}_{o}F'(u_{o})n_{1}=F'(u_{o})(\alpha^{2}_{o}n_{1})
\end{equation*}
where the second equality is just (\ref{242}). Then $\varphi_{1}(u_{o})-\alpha^{2}_{o}n_{1}\in N(F'(u_{o}))$ and thus there exists a unique $\alpha_{1}\in \mathbb{R}$ such that 
\begin{equation}\label{244}
\varphi_{1}(u_{o})=\alpha^{2}_{o}n_{1}+\alpha_{1}n_{o}.
\end{equation}
From this relation and thanks to the definition of the map $\varSigma_{2}$ we obtain
 \begin{align}\label{245}
\begin{split}
\varSigma_{2}[\varphi(u_{o}),\varphi_{1}(u_{o})]\!&=\!F^{(3)}(u_{o})[\varphi(u_{o}),\varphi(u_{o}),\varphi(u_{o})]+3F''(u_{o})[\varphi_{1}(u_{o}),\varphi(u_{o})]=\\
&=\!\alpha^{3}_{o}(F^{(3)}(u_{o})[n_{o},n_{o},n_{o}]+3F''(u_{o})[n_{1},n_{o}])+3\alpha_{1}\alpha_{o}F''(u_{o})[n_{o},n_{o}]\!=\\
&=\!\alpha^{3}_{o}\varSigma_{2}[n_{o},n_{1}]+3\alpha_{1}\alpha_{o}\varSigma_{1}[n_{o}].
\end{split}
\end{align}
Since $\alpha_{o}\neq 0$ and $\varSigma_{1}[n_{o}]\in R(F'(u_{o}))$ we easily get that $\varSigma_{2}[\varphi(u_{o}),\varphi_{1}(u_{o})]\in R(F'(u_{o}))$ \textit{ iff }$\varSigma_{2}[n_{o},n_{1}]\in R(F'(u_{o}))$.\\
\indent \textit{Proof of }$(t_{1})$. Let $n_{o},n_{1},v_{o},v_{1}$ satisfy condition (\ref{214}) of Theorem \ref{Teo212}. This means that $n_{o}\in N(F'(u_{o}))\setminus \{0\},F'(u_{o})n_{1}= -\varSigma_{1}[n_{o}], v_{o} \in X$ and $F'(u_{o})v_{1}=-T_{1}[n_{o},v_{o}]$. Since $(t_{o})$ is satisfied, the condition $T_{1}[n_{o},v_{o}]\in R(F'(u_{o}))$ is equivalent to $T_{1}[\varphi(u_{o}),v_{o}]\in R(F'(u_{o}))$ and thus, thanks to Remark \ref{Rem2111}, we have that $\tau_{1}(u_{o})v_{o}\in R(F'(u_{o}))$. Then, from formula $(j^{\,\prime}_{1})$ of Remark \ref{Rem234}, we obtain that $F'(u_{o})\varphi'(u_{o})v_{o}=-\tau_{1}(u_{o})v_{o}=-T_{1}[\varphi(u_{o}),v_{o}]$. This provides the relation between $\varphi'(u_{o})v_{o}$ and $v_{1}$. Indeed, by using (\ref{241}) and (\ref{243}) we find 
 \begin{equation*}
F'(u_{o})\varphi'(u_{o})v_{o}=-T_{1}[\varphi(u_{o}),v_{o}]=-\alpha_{o}T_{1}[n_{o},v_{o}]=-\alpha_{o}F'(u_{o})v_{1}=-F'(u_{o})(\alpha_{o}v_{1}).
\end{equation*}
Hence there exists a unique $\beta_{1}\in \mathbb{R}$ such that
\begin{equation}\label{246}
\varphi'(u_{o})v_{o} = \alpha_{o}v_{1} + \beta_{1}n_{o}.
\end{equation}
From this relation, using (\ref{241}) and (\ref{244}), we obtain 
\begin{align}\label{247}
\begin{split}
&\quad \;
T_{2}[\varphi(u_{o}),\varphi_{1}(u_{o}),v_{o},\varphi'(u_{o})v_{o}]=\\
&=F^{(3)}(u_{o})[\varphi(u_{o}),\varphi(u_{o}),v_{o}]+2F''(u_{o})[\varphi'(u_{o})v_{o},\varphi(u_{o})]+F''(u_{o})[\varphi_{1}(u_{o}),v_{o}]=\\
&=\alpha^{2}_{o}\big(F^{(3)}(u_{o})[n_{o},n_{o,}v_{o}]+2F''(u_{o})[v_{1},n_{o}]+F''(u_{o})[n_{1}+n_{o},v_{o}]\big)\,+\\
&\quad+2\beta_{1}\alpha_{o}F''(u_{o})[n_{o},n_{o}]+\alpha_{1}F''(u_{o})[n_{o},v_{o}]=\\
&=\alpha^{2}_{o}T_{2}[n_{o},n_{1},v_{o},v_{1}]+2\beta_{1}\alpha_{o}\varSigma_{1}[n_{o}]+\alpha_{1}T_{1}[n_{o},v_{o}].
\end{split}
\end{align}
Since $\alpha_{o}\neq 0$ and $\varSigma_{1}[n_{o}],T_{1}[n_{o},v_{o}] \in R(F'(u_{o}))$ we can conclude that  \\

$T_{2}[\varphi(u_{o}),\varphi_{1}(u_{o}),v_{o},\varphi'(u_{o})v_{o}] \in R(F'(u_{o}))$ \textit{ iff } $T_{2}[n_{o},n_{1},v_{o},v_{1}] \in R(F'(u_{o})). \; \bracevert$\\
\par
\textbf{Proof of Lemma \ref{Lem223} for k = 2}. Our goal is to prove that $(s_{o}), (t_{o}),(s_{1}),(t_{1})$ and $(s_{2}),(t_{2})$ hold when $u_{o}$ is a 2-transverse singularity and is not a 2-singularity. From Lemma \ref{Lem232} we know that $u_{o}$ is a 1-transverse singularity but not a 1-singularity and thus, from the previous step, we get that $(s_{o}), (t_{o}),(s_{1}),(t_{1})$ and formulas (\ref{241}), ..., (\ref{247}) are true.\\
\indent \textit{ Proof of }$(s_{2})$. If $u_{o}$ is a 2-transverse singularity which is not a 2-singularity, the equality $F'(u_{o})\varphi_{2}(u_{o})=-\sigma_{2}(u_{o})$ follows from $(i_{2})$ of Theorem \ref{Teo221}. Thanks to Remark \ref{Rem2111} this amounts to saying that $F'(u_{o})\varphi_{2}(u_{o})=-\varSigma_{2}[\varphi(u_{o}),\varphi_{1}(u_{o})]$ . Let $n_{o}, n_{1}, n_{2}$ satisfy condition (\ref{217}) of Theorem \ref{Teo213}, i.e. $n_{o}\in N(F'(u_{o}))\setminus \{0\}, F'(u_{o})n_{1}=-\varSigma_{1}[n_{o}]$ and $F'(u_{o})n_{2}=-\varSigma_{2}[n_{o},n_{1}]$. As in the previous step, we look for a relation between $\varphi_{2}(u_{o})$ and $n_{2}$. We already know that
\begin{equation*}
\varphi(u_{o})=\alpha_{o}n_{o} \,,\, \varphi_{1}(u_{o})=\alpha^{2}_{o}n_{1} + \alpha_{1}n_{o}\, 
,\end{equation*}
for suitable $\alpha_{o} \in \mathbb{R}\setminus \{0\}$ and $\alpha_{1} \in \mathbb{R}$ (see (\ref{241}) and (\ref{244})). From (\ref{245}) we have that
\begin{align*}
\begin{split}
F'(u_{o})\varphi_{2}(u_{o}) &=-\varSigma_{2}[\varphi(u_{o}),\varphi_{1}(u_{o})]=\alpha^{3}_{o}(-\varSigma_{2}[n_{o},n_{1}])+3\alpha_{1}\alpha_{o}(-\varSigma_{1}[n_{o}])=\\
&=\alpha^{3}_{o}(F'(u_{o})n_{2})+3\alpha_{1}\alpha_{o}(F'(u_{o})n_{1})=F'(u_{o})(\alpha^{3}_{o}n_{2}+3\alpha_{1}\alpha_{o}n_{1}).
\end{split}
\end{align*}
Hence there exists a unique $\alpha_{2}\in \mathbb{R}$ such that  
\begin{equation}\label{248}
\varphi_{2}(u_{o})=\alpha^{3}_{o}n_{2}+3\alpha_{1}\alpha_{o}n_{1}+\alpha_{2}n_{o}.
\end{equation}
This implies, from the definition of $\varSigma_{3}$ and by using (\ref{241}), (\ref{244}), that 
\begin{align}\label{249}
\begin{split}
&\quad \;
\varSigma_{3}[\varphi(u_{o}),\varphi_{1}(u_{o}),\varphi_{2}(u_{o})]=\\
&=F^{(4)}(u_{o})[\varphi(u_{o}),\varphi(u_{o}),\varphi(u_{o}),\varphi(u_{o})]+6F^{(3)}(u_{o})[\varphi_{1}(u_{o}),\varphi(u_{o}),\varphi(u_{o})]\,+\\
&\quad+4F''(u_{o})[\varphi_{2}(u_{o}),\varphi(u_{o})]+3F''(u_{o})[\varphi_{1}(u_{o}),\varphi_{1}(u_{o})]=\\
&=\alpha^{4}_{o}\big(F^{(4)}(u_{o})[n_{o},n_{o},n_{o},n_{o}]+6F^{(3)}(u_{o})[n_{1},n_{o},n_{o}]+4F''(u_{o})[n_{2},n_{o}]\,+\\
&\quad+3F''(u_{o})[n_{1},n_{1}]\big)+6\alpha_{1}\alpha^{2}_{o}\big(F^{(3)}(u_{o})[n_{o},n_{o},n_{o}]+3F''(u_{o})[n_{1},n_{o}]\big)\,+\\
&\quad+(4\alpha_{2}\alpha_{o}+3\alpha^{2}_{1})F''(u_{o})[n_{o},n_{o}]=\\
&=\alpha^{4}_{o}\varSigma_{3}[n_{o},n_{1},n_{2}]+6\alpha_{1}\alpha^{2}_{o}\varSigma_{2}[n_{o},n_{1}]+(4\alpha_{2}\alpha_{o}+3\alpha^{2}_{1})\varSigma_{1}[n_{o}].
\end{split}
\end{align}
Since $\varSigma_{2}[n_{o},n_{1}]$ and $\varSigma_{1}[n_{o}]$ belong to $R(F'(u_{o}))$ and $\alpha_{o}\in \mathbb{R}\setminus \{0\}$ we easily get that $\varSigma_{3}[\varphi(u_{o}),\varphi_{1}(u_{o}),\varphi_{2}(u_{o})]\in R(F'(u_{o}))\; \textit{ iff }\;\varSigma_{3}[n_{o},n_{1},n_{2}] \in R(F'(u_{o}))$.\\
\indent\textit{ Proof of }$(t_{2})$. Let $n_{o},n_{1},n_{2},v_{o},v_{1},v_{2}$  satisfy condition (\ref{217}) of Theorem \ref{Teo213}, i.e. $n_{o}\in N(F'(u_{o}))\setminus \{0\}, F'(u_{o})n_{1}=-\varSigma_{1}[n_{o}], F'(u_{o})n_{2}=-\varSigma_{2}[n_{o},n_{1}], v_{o} \in X, F'(u_{o})v_{1}=-T_{1}[n_{o},v_{o}]$ and $F'(u_{o})v_{2} = -T_{2}[n_{o},n_{1},v_{o},v_{1}]$. When we proved the previous equivalences we already established equalities (\ref{241}), (\ref{244}), (\ref{246}), (\ref{248}), i.e.
\begin{align*}
&\varphi(u_{o})\; = \alpha_{o}n_{o},\\
&\varphi_{1}(u_{o}) = \alpha^{2}_{o}n_{1}+\alpha_{1}n_{o},\\
&\varphi_{2}(u_{o}) = \alpha^{3}_{o}n_{2}+3\alpha_{1}\alpha_{o}n_{1}+\alpha_{2}n_{o},\\
&\varphi'(u_{o})v_{o}=\alpha_{o}v_{1}+\beta_{1}n_{o},
\end{align*}
for suitable $\alpha_{o},\alpha_{1},\alpha_{2},\beta_{1} \in \mathbb{R},\alpha_{o} \neq 0$. An analogous formula is to be found for $\varphi^{\prime}_{1}(u_{o})v_{o}$ and $v_{2}$. We claim that 
\begin{equation}\label{2410}
F'(u_{o})\varphi^{\prime}_{1}(u_{o})v_{o}=-T_{2}[\varphi(u_{o}),\varphi_{1}(u_{o}),v_{o},\varphi'(u_{o})v_{o}].
\end{equation} 
This claim follows, as seen above, from relation $(j^{\,\prime}_{2})$ of Remark \ref{Rem234} and from Remark \ref{Rem2111}. In order to apply $(j^{\,\prime}_{2})$ we need to know that $\tau_{1}(u_{o})v_{o},\tau_{2}(u_{o})v_{o}\in R(F'(u_{o}))$. This is true \textit{iff} $T_{1}[\varphi(u_{o}),v_{o}]$ and $T_{2}[\varphi(u_{o}),\varphi_{1}(u_{o}),v_{o},\varphi'(u_{o})v_{o}]$ belong to $R(F'(u_{o}))$, thanks to Remark \ref{Rem2111}. In turn, this holds since the equivalences $(t_{o})$ and $(t_{1})$, proved above, are satisfied and moreover $T_{1}[n_{o},v_{o}],\,T_{2}[n_{o},n_{1},v_{o},v_{1}]\in R(F'(u_{o}))$ by hypothesis. Then, from relation $(j^{\,\prime}_{2})$ of Remark \ref{Rem234}, we obtain that $F'(u_{o})\varphi^{\prime}_{1}(u_{o})v_{o}\!=-\tau_{2}(u_{o})v_{o}$ and (\ref{2410}) is thus proved thanks to Remark \ref{Rem2111}. From (\ref{2410}) and (\ref{247})  we get 
\begin{align*}
F'&(u_{o})\varphi^{\prime}_{1}(u_{o})v_{o}=-(\alpha^{2}_{o}
T_{2}[n_{o},n_{1},v_{o},v_{1}]+2\beta_{1}\alpha_{o}\varSigma_{1}[n_{o}]+\alpha_{1}T_{1}[n_{o},v_{o}])=\\
&=\alpha^{2}_{o}F'(u_{o})v_{2}+2\beta_{1}\alpha_{o}F'(u_{o})n_{1}+\alpha_{1}F'(u_{o})v_{1}=F'(u_{o})(\alpha^{2}_{o}v_{2}+2\beta_{1}\alpha_{o}n_{1}+\alpha_{1}v_{1}),
\end{align*}
hence there exists a unique $\beta_{2}\in \mathbb{R}$ such that
\begin{equation}\label{2411}
\varphi^{\prime}_{1}(u_{o})v_{o}=\alpha^{2}_{o}v_{2} + 2 \beta_{1}\alpha_{o}n_{1}+\alpha_{1}v_{1}+\beta_{2}n_{o}.
\end{equation}
Therefore,
 \begin{align}\label{2412}
\begin{split}
&\quad \; T_{3}[\varphi(u_{o}),\varphi_{1}(u_{o}),\varphi_{2}(u_{o}),v_{o},\varphi'(u_{o})v_{o},\varphi^{\prime}_{1}(u_{o})v_{o}]=\\
&=F^{(4)}(u_{o})[\varphi(u_{o}),\varphi(u_{o}),\varphi(u_{o}),v_{o}]+3F^{(3)}(u_{o})[\varphi(u_{o}),\varphi(u_{o}),\varphi'(u_{o})v_{o}]\,+\\
&\quad+3F^{(3)}(u_{o})[\varphi_{1}(u_{o}),\varphi(u_{o}),v_{o}]+3F''(u_{o})[\varphi(u_{o}),\varphi^{\prime}_{1}(u_{o})v_{o}]\,+\\
&\quad+3F''(u_{o})[\varphi_{1}(u_{o}),\varphi'(u_{o})v_{o}]+F''(u_{o})[\varphi_{2}(u_{o}),v_{o}]=\\
&=\alpha^{3}_{o}\big(F^{(4)}(u_{o})[n_{o},n_{o},n_{o},v_{o}]+3F^{(3)}(u_{o})[n_{o},n_{o},v_{1}]+3F^{(3)}(u_{o})[n_{1},n_{o},v_{o}]\,+\\
&\quad+3F''(u_{o})[n_{o},v_{2}]+3F''(u_{o})[n_{1},v_{1}]+F''(u_{o})[n_{2},v_{o}]\big)\,+\\
&\quad+3\alpha^{2}_{o}\beta_{1}\big(F^{(3)}(u_{o})[n_{o},n_{o},n_{o}]+2F''(u_{o})[n_{o},n_{1}]+F''(u_{o})[n_{1},n_{o}]\big)\,+\\
&\quad+3\alpha_{1}\alpha_{o}\big(F^{3}(u_{o})[n_{o},n_{o},v_{o}]+2F''(u_{o})[n_{o},v_{1}]+F''(u_{o})[n_{1},v_{o}]\big)\,+\\
&\quad+3\alpha_{o}\beta_{2}F''(u_{o})[n_{o},n_{o}]+3\alpha_{1}\beta_{1}F''(u_{o})[n_{o},n_{o}]+\alpha_{2}F''(u_{o})[n_{o},v_{o}]=\\
&=\alpha^{3}_{o}T_{3}[n_{o},n_{1},n_{2},v_{o},v_{1},v_{2}]+3\alpha^{2}_{o}\beta_{1}\varSigma_{2}[n_{o},n_{1}]+3\alpha_{1}\alpha_{o}T_{2}[n_{o},n_{1},v_{o},v_{1}]\,+\\
&\quad+3\alpha_{o}\beta_{2}\varSigma_{1}[n_{o}]+3\alpha_{1}\beta_{1}\varSigma_{1}[n_{o}]+\alpha_{2}T_{1}[n_{o},v_{o}].
\end{split}
\end{align}
Since $\alpha_{o}\neq 0$ and $\varSigma_{1}[n_{o}],\varSigma_{2}[n_{o},n_{1}],T_{1}[n_{o},v_{o}],T_{2}[n_{o},n_{1},v_{o},v_{1}]\in R(F'(u_{o}))$ we conclude that $(t_{2})$ is satisfied. $\bracevert$\\
\par
\textbf{Proof of Lemma \ref{Lem223} for k = 3}. The scheme of this proof is again based on the points a) and b) given at the beginning of this section. For this reason we only sketch the proof of equivalence $(s_{3})$ and we fully omit the proof of $(t_{3})$. Here we have to show that $(s_{o}), (s_{1}), (s_{2})$ and $(s_{3})$ hold if $u_{o}$ is a 3-transverse singularity which is not a 3-singularity. Once again, by Lemma \ref{Lem232} one has that $u_{o}$ is a 2-transverse singularity which is not a 2-singularity and hence we already know, from the previous step, that $(s_{o}), (s_{1}), (s_{2})$ and formulas (\ref{241}), ..., (\ref{249}) hold.\\
\indent
\textit{Proof of }$(s_{3})$. Since $u_{o}$ is a 3-transverse singularity which is not a 3-singularity then, from equality $(i_{3})$ of Theorem \ref{Teo221},  $F'(u_{o})\varphi_{3}(u_{o})=-\sigma_{3}(u_{o})$; also, we have that $F'(u_{o})\varphi_{3}(u_{o})=-\varSigma_{3}[\varphi(u_{o}),\varphi_{1}(u_{o}),\varphi_{2}(u_{o})]$ from Remark \ref{Rem2111}\\
Let now $n_{o},n_{1},n_{2},n_{3}$ satisfy condition (\ref{2110}) of Theorem \ref{Teo214}. Once more, the relation between $\varphi_{3}(u_{o})$ and $n_{3}$ is found by using
\begin{align*}
&\varphi(u_{o})\; = \alpha_{o}n_{o},\\
&\varphi_{1}(u_{o}) = \alpha^{2}_{o}n_{1}+\alpha_{1}n_{o},\\
&\varphi_{2}(u_{o}) = \alpha^{3}_{o}n_{2}+3\alpha_{1}\alpha_{o}n_{1}+\alpha_{2}n_{o}\, ,
\end{align*}
where $\alpha_{o}\neq 0, \alpha_{1}, \alpha_{2}$ are suitable real numbers, see (\ref{241}), (\ref{244}) and (\ref{248}). An easy computation gives that 
\begin{align*}
\begin{split}
F'(u_{o})\varphi_{3}(u_{o})\!&=-\varSigma_{3}[\varphi(u_{o}),\varphi_{1}(u_{o}),\varphi_{2}(u_{o})]=\\
&=\alpha^{4}_{o}(-\varSigma_{3}[n_{o},n_{1},n_{2}])\!+\!6\alpha_{1}\alpha^{2}_{o}(-\varSigma_{2}[n_{o},n_{1}])\!+\!(4\alpha_{2}\alpha_{o}+3\alpha^{2}_{1})(-\varSigma_{1}[n_{o}])\!=\\
&=\alpha^{4}_{o}(F'(u_{o})n_{3})+6\alpha_{1}\alpha^{2}_{o}(F'(u_{o})n_{2})+(4\alpha_{2}\alpha_{o}+3\alpha^{2}_{1})(F'(u_{o})n_{1})=\\
&=F'(u_{o})\big(\alpha^{4}_{o}n_{3}+6\alpha_{1}\alpha^{2}_{o}n_{2}+(4\alpha_{2}\alpha_{o}+3\alpha^{2}_{1})n_{1}\big)
\end{split}
\end{align*}
where in the second equality we used (\ref{249}). From this it follows that there exists a suitable $\alpha_{3}\in \mathbb{R}$ such that
\begin{equation}\label{2413}
\varphi_{3}(u_{o})=\alpha^{4}_{o}n_{3}+6\alpha_{1}\alpha^{2}_{o}n_{2}+(4\alpha_{2}\alpha_{o}+3\alpha^{2}_{1})n_{1}+\alpha_{3}n_{o}.
\end{equation}
Finally, by using the very definition of $\varSigma_{4}$ and formulas (\ref{241}), (\ref{244}), (\ref{248}) and (\ref{2413}) we get 
\begin{align*}
\begin{split}
&\quad \;
\varSigma_{4}[\varphi(u_{o}),\varphi_{1}(u_{o}),\varphi_{2}(u_{o}),\varphi_{3}(u_{o})]=F^{(5)}(u_{o})[\varphi(u_{o}),\varphi(u_{o}),\varphi(u_{o}),\varphi(u_{o}),\varphi(u_{o})]\,+\\
&\quad+10F^{(4)}(u_{o})[\varphi_{1}(u_{o}),\varphi(u_{o}),\varphi(u_{o}),\varphi(u_{o})]+10F^{(3)}(u_{o})[\varphi_{2}(u_{o}),\varphi(u_{o}),\varphi(u_{o})]\,+\\
&\quad+15F^{(3)}(u_{o})[\varphi_{1}(u_{o}),\varphi_{1}(u_{o}),\varphi(u_{o})]+5F''(u_{o})[\varphi_{3}(u_{o}),\varphi(u_{o})]\,+\\
&\quad+10F''(u_{o})[\varphi_{2}(u_{o}),\varphi_{1}(u_{o})]=\\
&=\alpha^{5}_{o}\big(F^{(5)}(u_{o})[n_{o},n_{o},n_{o},n_{o},n_{o}]+10F^{(4)}(u_{o})[n_{1},n_{o},n_{o},n_{o}]+10F^{(3)}(u_{o})[n_{2},n_{o},n_{o}]\,+\\
&\quad+15F^{(3)}(u_{o})[n_{1},n_{1},n_{o}]+5F''(u_{o})[n_{3},n_{o}]+10F''(u_{o})[n_{2},n_{1}]\big)\,+\\
&\quad+10\alpha_{1}\alpha^{3}_{o}\big(F^{(4)}(u_{o})[n_{o},n_{o},n_{o},n_{o}]+6F^{(3)}(u_{o})[n_{1},n_{o},n_{o}]+4F''(u_{o})[n_{2},n_{o}]\,+\\
&\quad+3F''(u_{o})[n_{1},n_{1}]\big)+10\alpha_{2}\alpha^{2}_{o}\big(F^{(3)}(u_{o})[n_{o},n_{o},n_{o}]+3F''(u_{o})[n_{1},n_{o}]\big)\,+\\
&\quad+15\alpha^{2}_{1}\alpha_{o}\big(F^{(3)}(u_{o})[n_{o},n_{o},n_{o}]+3F''(u_{o})[n_{1},n_{o}]\big)\,+\\
&\quad+5\alpha_{4}\alpha_{o}F''(u_{o})[n_{o},n_{o}]+10\alpha_{2}\alpha_{1}F''(u_{o})[n_{o},n_{o}]=\\
&=\alpha^{5}_{o}\varSigma_{4}[n_{o},n_{1},n_{2},n_{3}]+10\alpha_{1}\alpha^{3}_{o}\varSigma_{3}[n_{o},n_{1},n_{2}]+(10\alpha_{2}\alpha^{2}_{o}+15\alpha^{2}_{1}\alpha_{o})\varSigma_{2}[n_{o},n_{1}]\,+\\
&\quad+(5\alpha_{3}\alpha_{o}+10\alpha_{2}\alpha_{1})\varSigma_{1}[n_{o}]
\end{split}
\end{align*}
and this concludes the proof. $\bracevert$\\
\vspace{2pt}
\subsection{An Example of Swallow's Tail Singularity for a Class of Differential Problems}\label{ss25}
\quad As an application of the pointwise characterization of low-order singularities studied in the previous sections, here we wish to give a local existence and multiplicity result for the Neumann problem of a Liénard equation, i.e.  
\begin{equation*}
\text{(P)}\; \begin{cases}
u''+f(u)u'+g(u)=h \qquad \text{in  }(0,\pi)\\
u'(0)= 0
=u'(\pi),
\end{cases}
\end{equation*}
where $f,g \in C^{\infty}(\mathbb{R})$ and, for a given right-hand side $h \in C^{0}([0,\pi])$, we are interested in solutions $u \in C^{2}_{N}([0,\pi])\equiv\{u \in C^{2}([0,\pi]):u'(0)=0=u'(\pi)\}$. We study problem (P) by considering the map $F:X \rightarrow Y$, where $X:=C^{2}_{N}([0,\pi])$ and $Y:=C^{0}([0,\pi])$, defined as $F(u):=u''+f(u)u'+g(u),u\in X$. We shall see that, under suitable assumptions for the nonlinearities $f$ and $g$, the function $u\equiv 0$ is a \textit{swallow's tail} singularity for $F$. We can then apply the local behaviour results of Chapter \ref{s1} in order to obtain a full description of the number of solutions near $u\equiv 0$ to problem (P) for $h$ in a suitable neighbourhood of $F(0)$; particularly, we can say that there exists $h$ near $F(0)$, in fact an open set of such $h$'s, with exactly four solutions near $u\equiv 0$. \\
\indent To our knowledge this provides an explicit example of a swallow's tail singularity for a class of differential problems; the other known examples, see \cite{Da} and \cite{Ruf2}, are of implicit type though in slightly different ways. In \cite{Da} the author studied the Dirichlet problem for a semilinear equation where the nonlinearity is a fourth degree polynomial with non-constant coefficients which depend on two parameters and it is shown that for some values of these parameters there exists a coefficient such that the function $u\equiv 0$ is a swallow's tail singularity for the associated map. On the other hand, the Neumann problem for an elliptic semilinear equation with a cubic nonlinearity depending on a real parameter is studied in \cite{Ruf2} and the author proves that, for suitable values of the parameter, there exists a butterfly singularity for the associated map. Hence the stratification of singularities (see e.g. Remark 2.4.6 and Theorem 2.5.4 in \cite{B-D 1}) allows to saying that there exist swallow's tail singularities in a suitable neighbourhood of the butterfly singularity. We also point out a theoretical result in \cite{C-Ti} on the existence of $k$-singularities, for all $k$: suitable coefficients (functions) can be determined in abstract for a nonlinear elliptic boundary value prototype problem with Dirichlet conditions so that the associated map is a \textit{global} $k$-singularity.\\
\indent As a matter of notation, we set $0_{X}:=0\in X$ and $0_{Y}:=0\in Y$ for the origins of $X, Y$ respectively. In order to show that the singular point $0_{X}$ is a swallow's tail we shall use the pointwise conditions for low-order singularities stated in Theorems \ref{Teo211}, \ref{Teo212}, \ref{Teo213}. It is worthwhile to mention that, thanks to the symmetric nature of problem (P) at $0_{X}$, such conditions can be simplified and, in practice, the number of required calculations can be significantly reduced.\\
\indent We note that the technique used in this section can be suitably modified to study the Dirichlet problem for the same class of equations. Indeed, it can be shown that, under quite weak assumptions on $f$ and $g$, all the low-order singularities considered in this chapter can be found for both Neumann and Dirichlet problems, \cite{B-D 4}. However, when studying more general boundary value problems where all kinds of singular points can occur, the pointwise approach is no longer convenient. In order to deal with these cases in \cite{B-D 3} we develop a new approach to identify the simple singularities, which we call the \textit{dual characterization} of singularities. Several results of multiplicity of solutions to boundary value problems are obtained if the dual characterization introduced in \cite{B-D 3} is combined with the study of the local behaviour near a simple singularity given here in Chapter 1, \cite{B-D 5}.\\
\par Since the nonlinearities $f$ and $g$ in problem (P) are $C^{\infty}$ functions it is a well-known fact that $F$ is a $C^{\infty}$ map and it is straightforward to show that $F$ is a 0-Fredholm map. Then, we have: 

\begin{theorem}\label{Teo251}
\textit{ Let }$f, g\in C^{\infty}(\mathbb{R})$\textit{ satisfy the following conditions:}
\begin{equation*}
f(0)= 0 \,,\; f'(0) = 0 \,,\; g(0)= 0 \,,\; g'(0)= 1 \,,\; g''(0) \neq 0\,; \qquad\qquad\text{(I)}
\end{equation*}
\begin{equation*}
g'''(0) = \frac{5}{3}g''(0)^{2}\,,\;f'''(0)\neq \frac{11}{3}f''(0)g''(0).   \quad\qquad\qquad\qquad\qquad\;\text{(II)}
\end{equation*}
\textit{Then}\\
i)\quad $0_{X}$\textit{ is a swallow's tail singularity for }$F$;\\
ii) \textit{there exist neighbourhoods }$U$\textit{ and }$V$\textit{ of }$0_{X}$\textit{ and }$0_{Y}$\textit{ respectively such that for any }$h\in V$\textit{ problem }(P\textit{) has at most four solutions }$u\in U$. \textit{Moreover, there exists }$\tilde{h}\in V$\textit{ such that problem }(P)\textit{ has exactly four distinct solutions in }$U$.
\end{theorem}
\begin{remark}\label{Rem252}
It is quite easy to give examples of classes of functions $f$ and $g$ satisfying the above conditions (I) and (II). For instance, we can consider the following cases:\\
1) $f$ and $g$ are asymptotically linear:
\begin{align*}
f(u)&=log(cos h(au)), a \in \mathbb{R}\setminus \{0\},\\
g(u)&=(1+10\alpha^{2})u +\frac{3\alpha u^{2}-10\alpha^{2}u}{u^{2}+1}, \alpha \in \mathbb{R}\setminus \{0\}\,;
\end{align*}
2) $f$ is a polynomial and $g$ is asymptotically linear:
\begin{align*}
&f(u)=au^{2}+bu^{3}+\varSigma^{N}_{n=4}c_{n}u^{n},\, a,b,c_{n}   \in \mathbb{R},\, b \neq \frac{22}{3} a\alpha \,,\\
&g(u)=(1+10\alpha^{2})u +\frac{3\alpha u^{2}-10\alpha^{2}u}{u^{2}+1}, \alpha \in \mathbb{R}\setminus\{0\}\,;
\end{align*}
3) $f$ and $g$ are polynomials:
\begin{align*}
&f(u)=au^{2}+bu^{3}+\varSigma^{N}_{n=4}c_{n}u^{n},\, a,b,c_{n}   \in \mathbb{R},\, b \neq \frac{22}{3} a\alpha \,,\\
&g(u)= u+3\alpha u^{2} + 10\alpha^{2}u^{3}+\varSigma^{M}_{m=4}\gamma_{m}u^{m},\alpha\in \mathbb{R} \setminus\{0\},\gamma_{m}\in \mathbb{R}\,.
\end{align*}
\end{remark}
Part i) of Theorem\ref{Teo251} will be proved later on, using Theorem \ref{Teo213}. On the other hand, if we know that $0_{X}$ is a swallow's tail singularity then we can use the Local Multiplicity Theorem ( Theorem \ref{Teo115}) to give a direct proof of part ii).\\
\par \textbf{Proof of Theorem \ref{Teo251}, ii)}. Let us suppose that we have already shown that $0_{X}$ is a swallow's tail singularity (i.e. a 3-singularity) for $F$. Then from Theorem \ref{Teo115} we know that there exist neighbourhoods $U$ and $V$ of $0_{X}$ and $0_{Y}=F(0_{X})$ respectively such that for any $h\in V$ the equation $F(u)=h$ has at most four solutions $u\in U$, hence the same is true for problem (P). Moreover, there exists $\tilde{h}\in V$ such that $F(u)=\tilde{h}$, i.e. problem (P) with right-hand side $\tilde{h}$, has exactly four distinct solutions in $U$.\\
\par 
Before proving the first part of Theorem \ref{Teo251} it is convenient to list some properties of the map $F$ and recall some elementary results from calculus. As a matter of notation, from now on we will simply write 0 instead of $0_{X}$ when $0_{X}$ is a variable, e.g. $F'(0)$ instead of  $F'(0_{X})$.

\begin{remark}\label{Rem253} a) For any $u\in X$ and $k = 2, 3,\ldots,$ the following formulas for the derivatives of $F$ can be easily proved by direct computation and induction:\\
- if $v\in X$ then $F'(u):X\rightarrow Y$ is given by
\begin{equation*}
F'(u)v =v'' + f(u)v' + \big(f'(u)u' + g'(u)\big)v ;
\end{equation*}
- if $v_{i}\in X, i=1.\ldots,k,$ then $F^{k}(u):\underbrace{X \times \ldots\times X }_{k}\rightarrow Y$ has the form
\begin{equation*}
F^{k}(u)[v_{1},v_{2},\ldots,v_{k}]= f^{(k-1)}(u)(v_{1}\cdotp v_{2}\cdotp \ldots \cdotp v_{k})' + \big(f^{(k)}(u)u' + g^{k}(u)\big)v_{1}\cdotp v_{2} \cdotp \ldots \cdotp v_{k}.
\end{equation*}
In particular, evaluating at 0 and taking into account the conditions in (I) we have that
\begin{align}\label{251}
F'(0)v =& \: v'' + v ;\\ \label{252}
\text{for }k=2 \qquad\qquad\qquad\, F''(0)[v_{1}, v_{2}] =& \:g''(0)v_{1}v_{2} ;\\ \label{253}
\text{for }k=3 \qquad\qquad\; F^{(3)}(0)[v_{1}, v_{2},v_{3}] =&\: f''(0)(v_{1}v_{2}v_{3})' + g'''(0)v_{1}v_{2}v_{3} ;\\ \label{254}
\text{for }k=4 \qquad\quad F^{(4)}(0)[v_{1}, v_{2},v_{3},v_{4}] =&\: f'''(0)(v_{1}v_{2}v_{3}v_{4})' + g^{IV}(0)v_{1}v_{2}v_{3}v_{4}.
\end{align}
\indent b) From formula (\ref{251}) we get that $N(F'(0))= \{v\in X:F'(0)v = v''+ v =0_{Y} \}$ has dimension 1 and, precisely, that $N(F'(0))$ is spanned by the function $n(t):=\cos(t)$.\\
\indent  c) It is not difficult to see that $R(F'(0))=\{w \in Y: \int_{0}^{\pi} w(t)\cos(t)dt =0\}$. Let us consider the $L^{2}$-scalar product $<\cdotp,\cdotp>$, i.e.\;$<\!w,z\!>\;:= \int_{0}^{\pi} w(t)z(t)dt$, for $w,z \in L^{2}(0,\pi)$. Then, in terms of this product, we have that $w \in R(F'(0))$ \textit{ iff }$<\!w,n\!>\;=0$.\\
\indent d) In the following we will use the fact that the linear operator $F'(0)$ is \textit{symmetric} with respect to the scalar product $<\cdotp ,\cdotp>$, where this means that for any $v_{1}, v_{2}\in X$ 
 \begin{equation}\label{255}
<\!F'(0)v_{1},v_{2}\!>\; =\; <\!v_{1}, F'(0)v_{2}\!>\,.
\end{equation}
\indent e) We recall the following formula from calculus: for any $m \in \mathbb{N}$, 
\begin{equation*}
I_{2m}:= \int_{0}^{\pi/2} \cos^{2m}(t)dt = \frac{(2m-1)!!}{(2m)!!}\frac{\pi}{2}.
\end{equation*}
From this we have that 
\begin{equation}\label{256}
\int_{0}^{\pi} \cos^{2m}(t)dt = 2I_{2m} =\frac{(2m-1)!!}{(2m)!!}\pi , m \in \mathbb{N}.
\end{equation}
Finally, for any $m,\mu \in \mathbb{N}$
\begin{equation}\label{257}
\int_{0}^{\pi}\cos^{2m}(t)\sin(t)dt =\frac{2}{2m+1}\,;
\end{equation}
\begin{equation}\label{258}
\int_{0}^{\pi}\cos^{2m+1}(t)\sin^{\mu}(t)dt = 0.
\end{equation}
Formula (\ref{258}) is proved by setting $\cos^{2}(t)=1-\sin^{2}(t)$ and then integrating with respect to $\tau=\sin(t)$.
\end{remark}
\textbf{Proof of Theorem \ref{Teo251}, i)}. If the function $g$ satisfies the hypotheses in (I), then the point $0_{X}$ is a simple singularity for the map $F$ and $N(F'(0))$ is spanned by the function $n(t)=\cos(t)$ as seen in point b) of the above Remark. By using the pointwise conditions stated in Theorems \ref{Teo211}, \ref{Teo212}, \ref{Teo213}, we shall now prove the first part of Theorem 2.5.1. We use the hypotheses in (I) to show that $0_{X}$ is not a fold but is a 1-transverse singularity, while the conditions in (II) allow us to show that $0_{X}$ is not a cusp but is a 2-transverse singularity as well as a swallow's tail. Actually, since the derivative of the map $F$ satisfies formula (\ref{255}), simpler equivalent conditions can be used to classify the singularity $0_{X}$. Hence the full proof is given by the following steps:\\
\indent \textit{Step 1:} $0_{X}$ is not a fold singularity.\\
\indent \textit{Step 2:} $0_{X}$ is a 1-transverse singularity.\\
\indent \textit{Step 3:} $0_{X}$ is not a cusp singularity.\\
\indent \textit{Step 4:} $0_{X}$ is a 2-transverse singularity.\\
\indent \textit{Step 5:} $0_{X}$ is a swallow's tail singularity.\\
\indent \textit{Step 6:} The conditions used in Steps 4 and 5 are equivalent to those stated in Theorems \ref{Teo212} and \ref{Teo213}.\\
\par 
\textit{Proof of Step 1}. We use condition (\ref{213}) of Theorem \ref{Teo211}: for a given $n_{o} \in N(F'(0))\setminus \{0\}$, the point $0_{X}$ is a fold singularity \textit{iff}
\begin{equation*}
F''(0)[n_{o}, n_{o}]\notin R(F'(0)).
\end{equation*}
Hence $0_{X}$ is not a fold singularity \textit{iff}
\begin{equation}\label{259}
F''(0)[n_{o}, n_{o}]\in R(F'(0)).
\end{equation}
For the sake of simplicity, throughout the rest of the proof we will write $n$ instead of $n_{o}$ and always take 
 \begin{equation}\label{2510}
n(t)=\cos(t).
\end{equation}
From Remark \ref{Rem253}, c), condition (\ref{259}) is equivalent to
\begin{equation}\label{2511}
<\!F''(0)[n,n],n\!>\; =0 \, ,
\end{equation}
where $<\cdotp,\cdotp>$ is the $L^{2}$-scalar product. Given the explicit form of $F''$ (see \ref{252}) and the choice of $n(t)$, (\ref{2511}) becomes
\begin{equation}\label{2512}
g''(0) \int_{0}^{\pi} \cos^{3}(t)dt = 0 \, .
\end{equation}
This is true thanks to formula (\ref{258}) with $m=1$ and $\mu=0$.\\
\par
\textit{Proof of Step 2}. In order to show that $0_{X}$ is a 1-transverse singularity we use condition (\ref{212}) of Theorem \ref{Teo211}. For $n(t)$ as set above we have to find a vector  $v_{o}\in X$ such that 
\begin{equation}\label{2513}
F''(0)[n, v_{o}]\notin R(F'(0)).
\end{equation}
As in the previous step, this is equivalent to 
\begin{equation}\label{2514}
g''(0) \int_{0}^{\pi} \cos^{2}(t)v_{o}(t)dt \neq 0 \, .
\end{equation}
Since $g''(0) \neq 0 $ we only have to find $v_{o} \in X$ such that $\int_{0}^{\pi} \cos^{2}(t)v_{o}(t)dt \neq 0 $. To this end it suffices to take $v_{o}(t)=\cos^{2}(t)$ which belongs to $X$ and satisfies the condition. Hence $0_{X}$ is a 1-transverse singularity.\\
\par
\textit{Proof of Step 3.} The previous steps allow us to use Theorem \ref{Teo212} in order to show that $0_{X}$ is not a cusp singularity for $F$. Since $0_{X}$ is a 1-transverse singularity which is not a fold the conditions in (\ref{214}) of Theorem \ref{Teo212} are satisfied; hence $0_{X}$ is not a cusp singularity \textit{iff} there exist $n, n_{1}\in X$, as in (\ref{214}), such that
\begin{equation}\label{2515}
F^{(3)}(0)[n,n,n]+3F''(0)[n_{1},n] \in R(F'(0)) \,.
\end{equation}
From formulas (\ref{252}) and (\ref{253}) for $F''(0)$ and $F^{(3)}(0)$, (\ref{2515}) turns out to be equivalent to 
\begin{equation}\label{2516}
<\!f''(0)(n^{3})'+ g'''(0)n^{3},n\!> +\;3\!<\!g''(0)n_{1},n,n\!>\;= 0 \,.
\end{equation}
From the conditions in (\ref{214}), we have to find an element $n_{1}\in X$ such that 
\begin{equation*}
F'(0)n_{1} = -F''(0)[n,n] = -g''(0)n^{2} \,.
\end{equation*}
This amounts to choosing $n_{1}$ as a particular solution to the problem 
\begin{equation*}
\text{(P$_{1}$)}\, \begin{cases}
v''(t)+v(t)=-g''(0)\cos^{2}(t) \qquad (0,\pi)\\
v'(0)= 0 = v'(\pi).
\end{cases}
\end{equation*}
It is not difficult to verify that 
\begin{equation}\label{2517}
n_{1}(t):=\frac{g''(0)}{3}(\cos^{2}(t)-2)
\end{equation}
solves (P$_{1}$). Then, with the above choice for $n_{1}$, a simple calculation shows that (\ref{2516}) can be rewritten as 
\begin{equation*}
-3f''(0) \int_{0}^{\pi}\cos^{3}(t)\sin(t)dt + g'''(0) \int_{0}^{\pi}  \cos^{4}(t)dt + g''(0)^{2}\int_{0}^{\pi}\cos^{2}(t)(\cos^{2}(t)-2)dt=0.
\end{equation*}
By (\ref{258}) the first integral is zero and so the above condition can be reformulated as
\begin{equation}\label{2518}
(g'''(0)+ g''(0)^{2})\int_{0}^{\pi}\cos^{4}(t)dt - 2 g''(0)^{2}\int_{0}^{\pi}\cos^{2}(t)dt=0\,.
\end{equation}
On the other hand, by using formula (\ref{256}) we have that 
\begin{equation*}
\int_{0}^{\pi}\cos^{2}(t)dt = 2I_{2}= \frac{1}{2}\pi \; ,\; \int_{0}^{\pi}\cos^{4}(t)dt = 2I_{4}=\frac{3}{8}\pi \, .
\end{equation*}
Therefore condition (\ref{2518}) becomes 
\begin{equation*}
3g'''(0) -5g''(0)^{2} = 0 \, .
\end{equation*}
This last condition is verified because of the first hypothesis in (II). We have thus proved that $0_{X}$ is not a cusp.\\
\par
\textit{Proof of Step 4.} Since the point $0_{X}$ is a 1-transverse singularity which is not a fold, we use condition (\ref{215}) of Theorem \ref{Teo212} to prove that $0_{X}$ is 2-transverse. Precisely, $0_{X}$ is a 2-transverse singularity \textit{iff} there exist $n, n_{1}, v_{o}, v_{1} \in X$ which satisfy conditions (\ref{214}) of Theorem \ref{Teo212} and with $v_{o}$ such that
\begin{equation}\label{2519}
F^{(3)}(0)[n,n,v_{o}]+2F''(0)[n,v_{1}]+F''(0)[n_{1},v_{o}]\notin R(F'(0)) \,,
\end{equation}
i.e.
\begin{equation}\label{2520}
<\!F^{(3)}(0)[n,n,v_{o}],n\!> +\;2\!<\!F''(0)[n,v_{1}],n\!> +<\!F''(0)[n_{1},v_{o}],n\!>\;\neq 0\,.
\end{equation}
We claim that the above inequality can be rewritten as 
\begin{equation}\label{2521}
<\!F^{(3)}(0)[n,n,v_{o}],n\!> +\;3\!<\!F''(0)[n_{1},v_{o}],n\!>\; \neq 0\,.
\end{equation}
For this particular map $F$, this means that in order to show that $0_{X}$ is 2-transverse we need not know the element $v_{1}$. The equivalence between formulas (\ref{2520}) and (\ref{2521}) will be proved in Step 6. At this stage we assume that (\ref{2521}) holds and by choosing $n$ and $n_{1}$ as in the previous step we obtain that (\ref{2521}) becomes
\begin{align}\label{2522} 
\begin{split}
<\!\!f''(0)(n^{2}v_{o})' &+ g'''(0)n^{2}v_{o},n\!\!>\! +\, 3\!<\!\!g''(0)n_{1}v_{o},n\!\!>\, = \!f''(0)\int_{0}^{\pi}(\cos^{2}(t)v_{o}(t))'\cos(t)dt+\\
&+g'''(0)\int_{0}^{\pi}\cos^{3}(t)v_{o}(t)dt + g''(0)^{2}\int_{0}^{\pi}\cos(t)(\cos^{2}(t)-2)v_{o}(t)dt\neq 0 \, , 
\end{split}
\end{align}
where we used again formulas (\ref{252}) and (\ref{253}) and $v_{o}\in X$ is still to be determined; we claim that it suffices to take $v_{o}(t)=cos^{3}(t)$ in order to verify (\ref{2522}).\\
Before proving the last point, we have to show that the fourth condition in (\ref{214}) of Theorem \ref{Teo212} is satisfied for such a choice of $v_{o}$. Indeed, 
\begin{equation*}
F''(0)[n, v_{o}]\in R(F'(0)) \Leftrightarrow \; <\!F''(0)[n, v_{o}],n\!>\; = 0\, ,
\end{equation*}
i.e.
\begin{equation*}
g''(0)\int_{0}^{\pi}n^{2}v_{o}= g''(0)\int_{0}^{\pi}\cos^{5}(t)dt = 0 \,,
\end{equation*}
which is true by formula (\ref{258}).\\
For $v_{o}=cos^{3}(t)$, (\ref{2522}) becomes
\begin{align}\label{2523}
\begin{split}
&-5f''(0)\int_{0}^{\pi}\cos^{5}(t)\sin(t)dt\! + \!(g'''(0)+ g''(0)^{2}) \int_{0}^{\pi}\cos^{6}(t)dt
\!-\!2g''(0)^{2}\int_{0}^{\pi}\cos^{4}(t)dt=\\
&=(g'''(0)+g''(0)^{2})\int_{0}^{\pi}\cos^{6}(t)dt -2g''(0)^{2}\int_{0}^{\pi}\cos^{4}(t)dt \neq 0
\end{split}
\end{align}
as the first integral is zero thanks to formula (\ref{258}).\\
From (\ref{256}) we have that 
\begin{equation*}
\int_{0}^{\pi}\cos^{4}(t)dt = 2I_{4}= \frac{3}{8}\pi \, , \;\int_{0}^{\pi}\cos^{6}(t)dt = 2I_{6}= \frac{5}{16}\pi\, .
\end{equation*}
Thus, from (\ref{2523}) we get the equivalent condition
\begin{equation}\label{2524}
5 g'''(0) - 7 g'' (0)^{2} \neq 0 \, ,
\end{equation}
which is true from the hypotheses in (I) and (II), i.e. $g''(0) \neq 0$ and $g'''(0)=\frac{5}{3} g''(0)^{2}$. Therefore (\ref{2524}) is verified, hence $0_{X}$ is a 2-transverse singularity.\\
\par
\textit{Proof of Step 5.} Finally, we show that $0_{X}$ is a swallow's tail singularity for $F$. To this end, since the point $0_{X}$ is 2-transverse but is not a cusp (as shown in Steps 3 and 4) we use Theorem \ref{Teo213}. Given $n, n_{1}$ and $n_{2}\in X$ satisfying conditions (\ref{217}), we have to prove that (\ref{219}) is verified. By proceeding as in the above steps, we rewrite this condition as
\begin{equation}\label{2525}
<\!F^{(4)}(0)[n,n,n,n]+6F^{(3)}(0)[n_{1},n,n]+3F''(0)[n_{1},n_{1}]+4F''(0)[n_{2},n],n\!>\; \neq 0\,.
\end{equation}
We claim that (\ref{2525}) is equivalent to showing that
\begin{align}\label{2526}
\begin{split}
<\!F^{(4)}(0)[n,n,n,n],n\!>&+6\!<\!F^{(3)}(0)[n_{1},n,n],n\!>+\;4\!<\!F^{(3)}(0)[n,n,n],n_{1}\!>+\\
&+15<\!F''(0)[n_{1},n_{1}],n\!>\; \neq 0\,.
\end{split}
\end{align}
Hence in order to prove that $0_{X}$ is a swallow's tail we do not need to find the element $n_{2}$. The equivalence of (\ref{2525}) and (\ref{2526}) will be proved in Step 6. By assuming that such an equivalence holds, we can express the above formula in the following way:
 \begin{align}\label{2527}
\begin{split}
<\!f'''(0)(n^{4})' &+ g^{IV}(0)n^{4},n\!>+\;6\!<\!f''(0)(n_{1},n^{2})' + g'''(0)n_{1}n^{2},n\!> +\\
& +4\!<\!f''(0)(n^{3})' + g'''(0)n^{3},n_{1}\!> +\; 15\!<\!g''(0)n^{2}_{1},n\!> \;\neq 0\,,
\end{split}
\end{align}
where we used the form of the derivatives of $F$ at $0_{X}$, see Remark \ref{Rem253}, a).\\
If we choose $n$ and $n_{1}$ as in the previous steps, i.e. as in (\ref{2510}) and (\ref{2517}), formula (\ref{2527}) can be rewritten as
 \begin{align}\label{2528}
\begin{split}
&f'''(0)\int_{0}^{\pi}(\cos^{4}(t))'\cos(t)dt + g^{IV}(0)\int_{0}^{\pi}\cos^{5}(t)dt + \\
&+2f''(0)g''(0)\int_{0}^{\pi}\big(\cos^{2}(t)(\cos^{2}(t)-2)\big)'cos(t)dt+\\
&+2g'''(0)g''(0)\int_{0}^{\pi}(\cos^{2}(t)-2)\cos^{3}(t)dt+\frac{4}{3}f''(0)g''(0)\int_{0}^{\pi}(\cos^{3}(t))'(\cos^{2}(t)-2)dt+\\
&+\frac{4}{3}g'''(0)g''(0)\int_{0}^{\pi}(\cos^{2}(t)-2)\cos^{3}(t)dt+\frac{5}{3}(g''(0))^{3}\int_{0}^{\pi}(\cos^{2}(t)-2)^{2}\cos(t)dt\neq  0\,.
\end{split}
\end{align}
We notice that, thanks to (\ref{258}), four of the previous integrals are equal to zero; after a few  computations (\ref{2528}) can thus be reformulated as 
\begin{align*}
\begin{split}
&-f'''(0)\int_{0}^{\pi}4\cos^{4}(t)\sin(t)dt - 2f''(0)g''(0)\int_{0}^{\pi}(4\cos^{4}(t)-4\cos^{2}(t))\sin(t)dt+\\
&-4f''(0)g''(0)\int_{0}^{\pi}(\cos^{4}(t)-2\cos^{2}(t))\sin(t)dt \neq 0\,.
\end{split}
\end{align*}
From (\ref{257}), after some calculations, we get
\begin{equation*}
-3f'''(0)+11f''(0)g''(0) \neq 0\, ,
\end{equation*}
which is verified because of the second hypothesis in (II). Hence we proved that $0_{X}$ is a swallow's tail. \\
\par
\textit{Step 6.}\\
\textit{Proof of (\ref{2520}) $\Leftrightarrow$ (\ref{2521})}. It suffices to show that 
\begin{center}
$<\!F''(0)[n,v_{1}],n\!>$ = $<\!F''(0)[n_{1},v_{o}],n\!> .$
\end{center} 
By formula (\ref{252}) we have that 
 \begin{align*}
&<\!F''(0)[n,v_{1}],n\!>\; =g''(0)\int_{0}^{\pi}n(t)v_{1}(t)n(t)dt =\\
&=\;<\!g''(0)nn,v_{1}\!>\; =\; <\!F''(0)[n,n],v_{1}\!> \; =\; <\!-F'(0)n_{1},v_{1}\!>\, ,
\end{align*}
where the last equality follows from the third condition in \ref{214} of Theorem \ref{Teo212}. Since $F'(0)$ is symmetric, see (\ref{255}), then 
\begin{equation*}
<\!-F'(0)n_{1},v_{1}\!>\; =\;<\!n_{1}, -F'(0)v_{1}\!>\; =\;<\!n_{1},F''(0)[n,v_{o}]\!> \, ,
\end{equation*}
where the last equality follows from the fifth condition in (\ref{214}) of Theorem \ref{Teo212}. Hence 
\begin{equation*}
<\!F''(0)[n,v_{1}],n\!>\; =\;<\!n_{1},F''(0)[n,v_{o}]\!>\; =\;<\!F''(0)[n_{1},v_{o}],n\!> \, ,
\end{equation*}
where the last equality follows from the very definition of $F''(0)$.\\
\textit{Proof of (\ref{2525}) $\Leftrightarrow$ (\ref{2526})}. In this case we focus on the addendum $4\!<\!F''(0)[n_{2},n],n\!>$ in (\ref{2525}). Arguing as above, we obtain
\begin{equation*}
<\!F''(0)[n_{2},n],n\!>\; = F''(0)[n,n],n_{2}>\; =\; <\!-F'(0)n_{1},n_{2}\!> 
\end{equation*}
where the last equality is given by the third condition in (\ref{217}) of Theorem \ref{Teo213}. Since $F'(0)$ is symmetric we have that
 \begin{equation*}
<\!-F'(0)n_{1},n_{2}\!>\; =\;<\!n_{1}, -F'(0)n_{2}\!>\; =\;<\!n_{1},F^{(3)}(0)[n,n,n]+3F''(0)[n_{1},n]\!> 
\end{equation*}
thanks to the fifth condition in (\ref{217}) of Theorem \ref{Teo213}. Therefore
\begin{align*}
4\!<\!F''(0)[n_{2},n],n\!> &= 4\!<\!n_{1},F^{(3)}(0)[n,n,n]+ 3F''(0)[n_{1},n]\!>\; =\\
&=4\!<\!F^{(3)}(0)[n,n,n],n_{1}\!>+\; 12\!<\!F''(0)[n_{1},n_{1}],n\!>\,,
\end{align*}
and this suffices to conclude.$\bracevert$\\ \\

\end{document}